\documentclass[a4paper,11pt]{amsart}
\usepackage{graphicx}
\usepackage{amsmath}
\usepackage{amssymb}
\usepackage{amsthm}
\usepackage{amscd}
\usepackage{fullpage}
\usepackage[all]{xy}
\usepackage{tikz-cd}
\usepackage{latexsym}
\usepackage{mathrsfs}
\usepackage{dsfont} 
\usepackage{hyperref}

\newtheorem{theorem}{Theorem}[section]
\newtheorem{lemma}[theorem]{Lemma}
\newtheorem{corollary}[theorem]{Corollary}
\newtheorem{proposition}[theorem]{Proposition}
\newtheorem{introtheorem}{Theorem}

\theoremstyle{definition}
\newtheorem{remark}[theorem]{Remark}
\newtheorem{example}[theorem]{Example}
\newtheorem{definition}[theorem]{Definition}

\newcommand{\longmapsfrom}{\longleftarrow\!\shortmid}
\renewcommand{\emptyset}{\varnothing}
\newcommand{\Z}{\mathbb{Z}}
\newcommand{\Zc}{\mathcal{Z}}
\newcommand{\C}{\mathcal{C}}
\newcommand{\D}{\mathcal{D}}
\newcommand{\X}{\mathcal{X}}
\newcommand{\Y}{\mathcal{Y}}
\renewcommand{\P}{\mathcal{P}}
\newcommand{\Q}{\mathcal{Q}}
\newcommand{\T}{\mathcal{T}}
\newcommand{\U}{\mathcal{U}}
\newcommand{\V}{\mathcal{V}}
\newcommand{\A}{\mathcal{A}}
\newcommand{\B}{\mathcal{B}}
\newcommand{\R}{\mathcal{R}}
\newcommand{\K}{\mathbb{K}}
\renewcommand{\a}{\alpha}
\renewcommand{\b}{\beta}
\renewcommand{\l}{\lambda}
\newcommand{\g}{\gamma}
\newcommand{\ovl}{\overline}
\newcommand{\lora}{\longrightarrow}
\renewcommand{\S}{\Sigma}
\newcommand{\wt}{\widetilde}

\DeclareMathOperator{\Hom}{Hom}
\DeclareMathOperator{\Ext}{Ext}
\DeclareMathOperator{\id}{id}
\DeclareMathOperator{\add}{add}
\DeclareMathOperator{\ind}{ind}
\DeclareMathOperator{\PC}{PC}
\DeclareMathOperator{\PE}{PE}
\DeclareMathOperator{\PT}{PT}

\hypersetup{colorlinks,
linkcolor = {red!50!black},
citecolor = {blue!50!black},
urlcolor = {blue!80!black}
}

\newcommand{\harxiv}[1]{ \href{http://arxiv.org/abs/#1}{\texttt{arXiv:#1}}}

\begin{document}

\title{Torsion pairs, t-structures, and co-t-structures for completions of discrete cluster categories}

\author{Sofia Franchini} 
\address{Sofia Franchini, School of Mathematical Sciences, Lancaster University, UK.}
\email{s.franchini1@lancaster.ac.uk}

\begin{abstract}
We give a classification of torsion pairs, t-structures, and co-t-structures in the Paquette--Y\i ld\i r\i m completion of the Igusa--Todorov discrete cluster category. We prove that the aisles of t-structures and co-t-structures are in bijection with non-crossing partitions enriched with some additional data. We also observe that recollements exist in the completion and we classify them.
\end{abstract}

\maketitle
{\small
\tableofcontents}

\addtocontents{toc}{\protect\setcounter{tocdepth}{1}}

\section{Introduction}

Given a positive integer $m$ and a field $\K$, Igusa and Todorov defined in \cite{IT} a cluster category $\C_m$ which generalises the classical cluster category $\C(A_n)$ of type $A_n$ introduced by Buan, Marsh, Reineke, Reiten, and Todorov in \cite{BMRRT} for finite-dimensional hereditary algebras. The category $\C(A_n)$ has a nice  geometric model in terms of an $(n+3)$-gon, introduced by Caldero, Chapoton, and Schiffler in \cite{CCS}. When $m = 1$ or $m = 2$, the category $\C_m$ can be regarded as the orbit category of the infinite quiver $A_{\infty}$ or $A_{\infty}^{\infty}$, respectively, studied by Liu and Paquette in \cite{LP}, in analogy with the finite-rank case. When $m = 1$, the category $\C_m$ is equivalent to the Holm--J\o rgensen category defined in \cite{HJ} as the finite derived category of $\K[T]$ viewed as a graded algebra, or can be obtained by stabilising a certain subcategory of a Grassmannian category of infinite rank, see \cite{ACFGS}. When $m = 1$, $\C_m$ is also the unique algebraic triangulated category generated by a 2-spherical object, up to triangle equivalence. In particular, these categories come up in many different contexts.

The category $\C_m$ has many nice properties, for instance it is a $\Hom$-finite, $\K$-linear, Krull--Schmidt triangulated category, and has a geometric model which allows us to use combinatorial tools to classify some important classes of subcategories. The indecomposable objects of $\C_m$ can be regarded as the arcs of an $\infty$-gon, $\Zc_m$, having $m$ two-sided accumulation points. The $\Hom$-spaces are at most one-dimensional and can be understood in terms of crossings of arcs. Moreover, $\C_m$ is 2-Calabi--Yau,  
i.e. $\Ext^1(a,b) \cong D\Ext^1(b,a)$ for each pair of objects $a$ and $b$. The cluster-tilting subcategories of $\C_m$ were classified in \cite{GHJ}, and in \cite{HJ} and \cite{LP} for $m = 1$ and $m = 2$, as certain triangulations of the infinity-gon $\Zc_m$.

We can also be interested in studying a completion $\ovl\C_m$ of $\C_m$. Cummings and Gratz studied in \cite{CG} the Neeman's completion of $\C_m$. We work with the Paquette--Y\i ld\i r\i m completion of $\C_m$, defined in \cite{PY}, which was obtained by taking the Verdier quotient of $\C_{2m}$ with respect to a specific thick subcategory. In this article, by completion we mean the Paquette--Y\i ld\i r\i m completion of $\C_m$. This was first defined by Fisher in \cite{F} for the case $m = 1$, by closing the category $\C_1$ under certain homotopy colimits. August, Cheung, Faber, Gratz, and Schroll proved in \cite{ACFGS} that the category $\ovl\C_1$ is also equivalent to a stable Grassmannian category of infinite rank. The completion inherits many properties from $\C_m$, for instance, $\ovl\C_m$ is still a $\Hom$-finite, $\K$-linear, Krull--Schmidt triangulated category and has also a geometric model similar to the one for $\C_m$. The indecomposable objects of $\ovl\C_m$ are in bijection with arcs, or limits of arcs, of $\Zc_m$, and the $\Hom$-spaces are still at most one-dimensional. Moreover, $\ovl\C_m$ also has cluster-tilting subcategories which are in correspondence with some triangulations of the $\infty$-gon $\Zc_m$. Çanakçi, Kalck, and Pressland endowed $\ovl\C_m$ with an extriangulated structure $\mathbb{E}$ and classified the cluster-tilting subcategories with respect to $\mathbb{E}$ in terms of a larger class of triangulations of the $\infty$-gon $\Zc_m$, see \cite{CKP}. 

Despite $\C_m$ and $\ovl\C_m$ having many similarities, these two categories also have relevant differences. One remarkable difference is that $\ovl\C_m$ is not 2-Calabi--Yau, although it is ``weakly 2-Calabi--Yau" with respect to the extriangulated structure of \cite{CKP}. Therefore, classifying subcategories in $\ovl\C_m$ is a way to develop intuition for a more general setting, which is not necessarely 2-Calabi--Yau.  

The geometric models of $\C_m$ and $\ovl\C_m$ allow one to classify some important classes of subcategories using arc combinatorics. Torsion pairs of a triangulated category play an important role in representation theory, as they provide a decomposition of the category into smaller subcategories. By a result of Iyama and Yoshino in \cite{IY}, the torsion pairs in a ``small" triangulated category are completely determined by their torsion classes, which are characterised as extesion-closed precovering subcategories. In order to classify torsion pairs, we classify the precovering subcategories and the extension-closed subcategories of $\ovl\C_m$. We keep those properties separate and independent from each other.

\begin{introtheorem}[{Theorem \ref{theorem classification precovering subcategories}}]
Let $\X$ be an additive full subcategory of $\ovl\C_m$. Then $\X$ is a precovering subcategory of $\ovl\C_m$ if and only if the set of arcs corresponding to the indecomposable objects of $\X$ satisfies the completed precovering condition, i.e. it is closed under certain configurations of converging sequences of arcs of $\X$.
\end{introtheorem}

\begin{introtheorem}[{Proposition \ref{proposition extension closed subcategories completion}}]	
Let $\X$ be an additive full subcategory of $\ovl\C_m$. Then $\X$ is an extension-closed subcategory of $\ovl\C_m$ if and only if the set of arcs corresponding to the indecomposable objects of $\X$ satisfies the completed Ptolemy condition, i.e. it is closed under taking Ptolemy arcs for each pair of crossing arcs of $\X$.
\end{introtheorem}

The torsion pairs in $\C_m$ were classified by Gratz, Holm and J\o rgensen in \cite{GHJ} generalising the classifications of Ng in \cite{N} and of Chang, Zhou, and Zhu in \cite{CZZ} for the cases $m = 1$ and $m = 2$ respectively. By combining the two results above, we classify the torsion pairs in $\ovl\C_m$.

\begin{introtheorem}[{Theorem \ref{theorem classification of torsion pairs}}]
Let $\X$ be an additive full subcategory of $\ovl\C_m$. Then $\X$ is a torsion class in $\ovl\C_m$ if and only if the set of arcs corresponding to the indecomposable objects of $\X$ satisfies the completed precovering condition and the completed Ptolemy condition.
\end{introtheorem}

Particular kinds of torsion pairs are t-structures and co-t-structures, which are similar concepts but have important differences. For instance, t-structures are related to a notion of homology, and the heart of a t-structure is abelian, while the co-heart of a co-t-structure is presilting. The t-structures in $\C_m$ were first classified for $m = 1$ or $m = 2$ in \cite{N} and \cite{CZZ}. Gratz and Zvonareva classified the t-structures for $m\geq 1$ using decorated non-crossing partitions. These consist of a non-crossing partition of the finite set $\{1,\dots,m\}$ together with some additional data consisting of elements of the closure of the $\infty$-gon $\Zc_m$. In \cite{GZ} the authors also proved that the t-structures in $\C_m$ form a lattice under inclusions of aisles. This reflects the fact that the non-crossing partitions form a lattice under refinement. In this paper we classify the t-structures of $\ovl\C_m$ using similar combinatorial objects.

\begin{introtheorem}[{Theorem \ref{theorem t-structures}, Theorem \ref{theorem lattice structures t-structures and co-t-structures}}]
There is a bijection between the aisles of t-structures in $\ovl\C_m$ and the half-decorated non-crossing partitions of $\{1,\dots,2m\}$. Moreover, the t-structures in $\ovl\C_m$ form a lattice under inclusion of aisles.
\end{introtheorem}

We also classify the co-t-structures in $\ovl\C_m$ using related combinatorial objects. In $\C_m$ the only co-t-structures are the trivial ones $(0,\C_m)$ and $(\C_m,0)$. This marks a further difference with the completion, where non-trivial co-t-structures exist.

\begin{introtheorem}[{Theorem \ref{theorem classification aisles co-t-structures}, Theorem \ref{theorem lattice structures t-structures and co-t-structures}}]
There is a bijection between the aisles of the co-t-structures in $\ovl\C_m$ and the alternating non-crossing partitions of $\{1,\dots,2m\}$. Moreover, the co-t-structures in $\ovl\C_m$ form a lattice under inclusion of aisles.
\end{introtheorem}

Another interesting aspect of triangulated categories are their recollements. These can be thought as exact sequences of triangulated categories. Recollements are in bijection with torsion-torsion free triples, i.e. triples $(\X,\Y,\Zc)$ where $(\X,\Y)$ and $(\Y,\Zc)$ are t-structures, see for instance \cite{NS}. In our context, where the triangulated categories are $\Hom$-finite, $\K$-linear, and Krull--Schmidt, recollements are also in bijection with the functorially finite thick subcategories. Gratz and Zvonareva classified the thick subcategories of $\C_m$ in \cite{GZ}, and Murphy classified the thick subcategories of $\ovl\C_m$ in \cite{M}. The category $\C_m$ can be thought as ``triangulated simple", as its only functorially finite thick subcategories are $0$ and $\C_m$. The completion has different behaviour, indeed in $\ovl\C_m$ more functorially finite thick subcategories exist. We have the following classification of the functorially finite thick subcategories. 

\begin{introtheorem}[{Corollary \ref{corollary functorially finite thick subcategories}, Corollary \ref{corollary functorially finite thick subcategories lattice structure}}]
There is a bijection between the functorially finite thick subcategories of $\ovl\C_m$ and certain alternating non-crossing partitions of $\{1,\dots,2m\}$. Moreover, the functorially finite thick subcategories of $\ovl\C_m$ form a lattice under inclusion.
\end{introtheorem}

 \subsection*{Acknowledgments} The author thanks her supervisor David Pauksztello for useful conversation and suggestions and for carefully reading previous versions of this paper, and Raquel Coelho Sim\~{o}es for useful technical support and advice. The author thanks the referee for carefully reading this paper and for the useful comments. The author acknowledges support by the EPSRC through a mathematical sciences studentship and the grant EP/V050524/1.
  
\section{Background}\label{section background}

Throughout this section $\T$ will be a $\Hom$-finite, $\K$-linear, Krull--Schmidt triangulated category with shift functor $\S\colon\T\to\T$, unless otherwise stated. We denote by $\ind\T$ the class of indecomposable objects of $\T$. 

Any subcategory $\X$ of $\T$ is assumed to be full, and we sometimes write $\X\subseteq\T$. We say that $\X$ is an \emph{additive subcategory} if $\X$ is closed under direct sums, isomorphisms, direct summands, and contains the zero object. Given $\X$ and $\Y$ subcategories of $\T$, we write 
\[
\X\ast\Y = \left\{t\in\T\middle|\text{ there exists a triangle $x\lora t\lora y\lora \S x$ for some $x\in\X$ and $y\in\Y$}\right\}.
\]
A subcategory $\X$ of $\T$ is called
\begin{itemize}
\item \emph{extension-closed} if $\X\ast\X = \X$;
\item \emph{suspended} if it is extension-closed, additive, and $\S\X\subseteq\X$;
\item \emph{co-suspended} if it is extension-closed, additive, and $\S^{-1}\X\subseteq\X$;
\item \emph{thick} if it is suspended and co-suspended.
\end{itemize}
Given $\X$ and $\Y$ additive subcategories of $\T$, we write $\Hom(\X,\Y) = 0$ if $\Hom(x,y) = 0$ for each $x\in\X$ and $y\in\Y$. We denote 
\[
\X^\perp = \{t\in\T\mid \Hom_{\D}(\X,t) = 0\}\quad\text{and}\quad ^\perp\X = \{t\in\T\mid \Hom_{\D}(t,\X) = 0\}.
\] 
Let $\X$ be an additive subcategory of $\T$ and let $t\in\T$. We say that a morphism $f\colon x\to t$ is an \emph{$\X$-precover} of $t$ if $x\in\X$ and any $g\colon x'\to t$ with $x'\in \X$ factors through $f$. We say that an $\X$-precover $f\colon x\to t$ is an \emph{$\X$-cover} if additionally it is \emph{right minimal}, i.e. for any $g\colon x\to x$ if $fg = f$ then $g$ is an isomorphism. Covers are unique up to isomorphism, while precovers are not. We say $\X$ is \emph{precovering} if any $t\in \T$ admits an $\X$-precover.
The notions of \emph{preenvelope}, \emph{envelope}, and \emph{preenveloping subcategory} are dual.
If $\X$ is precovering and preenveloping, we say that $\X$ is \emph{functorially finite}.

\begin{remark}\label{remark precovering at the level of indecomposables}
In our context being precovering can be checked at the level of the indecomposable objects. More precisely, $\X$ is a precovering subcategory of $\T$ if and only if for any $t\in\ind\T$ there exist $x\in\X$ and $f\colon x\to t$ such that any $g\colon x'\to t$ with $x'\in\ind\X$ factors through $f$, cf. \cite[p. 81]{AS}.
\end{remark}

\subsection{Torsion pairs}

Let $\X$ and $\Y$ be additive subcategories of $\T$. The pair $(\X,\Y)$ is called
\begin{itemize}
\item \emph{torsion pair} if $\Hom(\X,\Y) = 0$ and $\T = \X\ast\Y$, see \cite{IY};
\item \emph{t-structure} if it is a torsion pair and $\S\X\subseteq\X$, see \cite{BBDG};
\item \emph{co-t-structure} if it is a torsion pair and  $\S^{-1}\X\subseteq\X$, see \cite{P} and \cite{B} where they are called \emph{weight structures}. 
\end{itemize}
Let $(\X,\Y)$ be a torsion pair, then $\X$ is called \emph{torsion class} and $\Y$ is called \emph{torsion-free class}. If $(\X,\Y)$ is a t-structure or a co-t-structure, $\X$ is called \emph{aisle} and $\Y$ is called \emph{co-aisle}. The \emph{heart} of a t-structure $(\X,\Y)$ is $\X\cap\S\Y$. The \emph{co-heart} of a co-t-structure  $(\X,\Y)$ is $\X\cap\S^{-1}\Y$. Let $(\X,\Y)$ be a t-structure or a co-t-structure, we say that $(\X,\Y)$ is
\begin{itemize}
\item \emph{left bounded}, or \emph{right bounded}, if $\T = \bigcup_{n\in\Z}\S^n\X$, or $\T = \bigcup_{n\in\Z}\S^n\Y$ respectively;
\item \emph{bounded} if it is left bounded and right bounded;
\item \emph{left non-degenerate}, or \emph{right non-degenerate}, if  $\bigcap_{n\in\Z}\S^n\X = 0$, or $\bigcap_{n\in\Z}\S^n\Y = 0$ respectively;
\item \emph{non-degenerate} if it is left non-degenerate and right non-degenerate.
\end{itemize}
It is straightforward to check that if $(\X,\Y)$ is left bounded then it is right non-degenerate, and if it is right bounded then it is left non-degenerate. 

\begin{proposition}[{\cite[Proposition 2.3]{IY}}]\label{proposition iyama-yoshino torsion pairs}
Let $\X,\Y$ be additive subcategories of $\T$. Then $(\X,\Y)$ is a torsion pair if and only if $\X$ is extension-closed and precovering, and $\Y = \X^{\perp}$.
\end{proposition} 

We recall the following notion from \cite{B}. Let $(\X,\Y)$ be a co-t-structure.
\begin{itemize}
\item If $\X$ is functorially finite, then $(^\perp\X,\X)$ is called its \emph{left-adjacent t-structure}.
\item If $\Y$ is functorially finite, then $(\Y,\Y^{\perp})$ is called its \emph{right-adjacent t-structure}. 
\end{itemize} 

\subsection{Decomposition of triangles}\label{section decomposition of triangles}

In this section we provide a decomposition of certain triangles of $\T$ (Proposition \ref{proposition decomposition of a triangle in triangulated setting}) and we characterise the extension-closed subcategories of $\T$ (Proposition \ref{proposition decomposition of a triangle in triangulated setting}), under the assumption that $\T$ has at most one dimensional $\Hom$-spaces. These results will be useful in Section \ref{section extension-closed subcategories} for describing the extension-closed subcategories in the Igusa--Todorov discrete cluster categories. We need the following definition and lemma.

\begin{definition}[{\cite[Definition 4.2]{CP}}]
The morphisms $f_1\colon x_1\to y$ and $f_2\colon x_2\to y$ are \emph{factorization free} if there is no $g\colon x_1\to x_2$ such that $f_1 = f_2g$, and there is no $h\colon x_2\to x_1$ such that $f_2 = f_1h$.
\end{definition}

\begin{lemma}\label{lemma factorization free morphisms}
Let $x_1,\dots,x_n,y\in\ind\T$ and $f = (f_1,\dots,f_n)\colon \bigoplus_{i = 1}^n x_i\to y$ be a morphism of $\T$. If $f$ is right minimal then $f_1,\dots,f_n$ are pairwise factorization free.
\end{lemma}
\begin{proof}
Assume that there exists a morphism $\varphi\colon x_j\to x_i$ such that $f_i = f_j\varphi$ for some $i\neq j$. We define the morphism $\g = (\g)_{s,t}\colon \bigoplus_{i = 1}^n x_i\to \bigoplus_{i = 1}^n x_i$ as
\[
(\g)_{s,t} = 
\begin{cases}
1 & \text{if $s = t\neq i$,}\\
0 & \text{if $s = t = i$,}\\
\varphi & \text{if $s = j$ and $t = i$,}\\
0 & \text{otherwise.}
\end{cases}
\]
It is straightforward to check that $f\g = f$. Since $\g$ is not an isomorphism, we have a contradiction with the fact that $f$ is right minimal. We conclude that $f_1,\dots,f_n$ are pairwise factorization free.
\end{proof}

\begin{proposition}\label{proposition decomposition of a triangle in triangulated setting}
Assume that the $\Hom$-spaces of $\T$ are at most one dimensional. Let $a\lora e \lora b\overset{h}{\lora} \S a$ be a triangle in $\T$ with $a,b_1,\dots,b_n\in\ind\T$, $b = \bigoplus_{i = 1}^n b_i$, and $h = (h_1,\dots,h_n)$. Then there exist $b_1',\dots,b_k'\in\ind\T$ and a morphism $h' = (h_1',\dots,h_k')\colon b' = \bigoplus_{i = 1}^k b_i'\to \S a$ such that $b'$ is a direct summand of $b$, $h_1',\dots,h_k'$ are pairwise factorization free, and there is the following isomorphism of triangles.
\[
\begin{tikzcd}
a\ar[r]\ar[d,"1"] & e'\oplus b''\ar[r]\ar[d,"\wr"] & b'\oplus b''\ar[r, "{(h', 0)}"]\ar[d,"\wr"]  & \S a\ar[d,"1"]\\
a\ar[r] & e\ar[r]& b\ar[r,"h"]& \S a	
\end{tikzcd}
\]
\end{proposition}
\begin{proof}
Without loss of generality, we can assume that $h_1,\dots,h_n\neq 0$, see \cite[Lemma 3.1]{CP}. Since the $\Hom$-spaces are at most one dimensional, it is straightforward to check that $h\colon b\to \S a$ is an $\add\{b\}$-precover of $\S a$. Thus, by \cite[Lemma 4.1]{J} there exists $b',b''\in\add\{b\}$ and an isomorphism $\a\colon b'\oplus b''\lora b$ such that the composition
\[
b'\overset{\left(\begin{smallmatrix}1\\0\end{smallmatrix}\right)}{\lora}b'\oplus b''\overset{\a}{\lora} b\overset{h}{\lora}\S a
\] 
is an $\add\{b\}$-cover of $\S a$, which we denote by $h'\colon b'\to \S a$. We denote $h\a = (h',h'')\colon b'\oplus b''\to \S a$. Since $h'\colon b'\to \S a$ is an $\add\{b\}$-cover of $\S a$, there exists $\b\colon b''\to b'$ such that $h'\b = h''$, and then
\[
\left(\begin{matrix}h' & 0\end{matrix}\right)\left(\begin{matrix} 1 & \b \\ 0 & 1\end{matrix}\right) = \left(\begin{matrix}h' & h'' \end{matrix}\right) = h\a.
\]
As a consequence, we obtain the isomorphism of triangles in the claim. Since $h'$ is right minimal, by Lemma \ref{lemma factorization free morphisms} , $h_1',\dots,h_n'$ are pairwise factorization free.
\end{proof}

\begin{remark}\label{remark consequence lemma decomposition of a triangle in triangulated setting}
Keeping the notation of Lemma \ref{proposition decomposition of a triangle in triangulated setting}, since $h_1',\dots,h_n'$ are pairwise factorization free, we have that $b_i'\not\cong\S a$, $h_i'\neq 0$, and $b_i'\not\cong b_j'$ for each $i\neq j$.
\end{remark}

The following proposition provides a sufficient condition for checking that a subcategory is extension-closed.

\begin{proposition}\label{proposition closure under extensions}
Assume that the $\Hom$-spaces of $\T$ are at most one dimensional, and let $\U$ be an additive subcategory of $\T$. Assume that $\U$ is closed under extensions of the form $a\lora e\lora b \overset{h}{\lora} \S a$ with $b = \bigoplus_{i = 1}^n b_i$, $a,b_1,\dots,b_n\in\ind\T$, $h = (h_1,\dots,h_n)$, and $h_1,\dots,h_n$ pairwise factorization free. Then $\U$ is closed under extensions.
\end{proposition}
\begin{proof}
We divide the proof into claims.
	
\emph{Claim 1.} The subcategory $\U$ is closed under extensions of the form $a\lora e\lora b\lora \S a$ with $a\in\ind\U$ and $b\in\U$.
	
Consider a triangle $a\lora e\lora b\lora \S a$ with $a\in\ind\U$ and $b\in\U$. This is isomorphic to a triangle of the form $a\lora e'\oplus b''\lora b'\oplus b''\overset{(h',0)}{\lora} \S a$ where $h' = (h_1',\dots,h_k')\colon b' = \bigoplus_{i = 1}^k b_i'\to \S a$ is as in Proposition \ref{proposition decomposition of a triangle in triangulated setting}. The triangle $a\lora e'\lora b'\overset{h'}{\lora} \S a$ satisfies the assumptions of our statement and therefore $e'\in\U$. Moreover, since $e\cong e'\oplus b''$ and $b''$ is a direct summand of $b$, we have that $e\in\U$.
	
\emph{Claim 2.} The subcategory $\U$ is closed under extensions.
	
Consider a non-split extension $a\lora e\lora b\lora \S a$ in $\C_m$ with $\bigoplus_{i = 1}^k a_i$, $b = \bigoplus_{j = 1}^n b_j$, and $a_i,b_j\in\ind\U$ for each $i\in\{1,\dots,k\}$ and $j\in\{1,\dots,n\}$. We proceed by induction on $k$. If $k = 1$ then we have the statement by Claim 1. Assume that $k\geq 2$, and consider the following Octahedral Axiom diagram.
\[
\begin{tikzcd}
& a_k\ar[r,"1"]\ar[d] & a_k\ar[d,dashed] & \\
\S^{-1}b\ar[r]\ar[d,"1"'] & a\ar[r]\ar[d] & e\ar[r]\ar[d,dashed] & b\ar[d,"1"'] \\
\S^{-1}b\ar[r] & \bigoplus_{i = 1}^{k-1}a_i\ar[r]\ar[d,"0"'] & x\ar[r]\ar[d,dashed] & b\\
& \S a_k\ar[r,"1"] & \S a_k
\end{tikzcd}
\]
Consider the triangle $\bigoplus_{i = 1}^{k-1}a_i\lora x\lora b\lora \bigoplus_{i = 1}^{k-1}\S a_i$. By the induction hypothesis we obtain that $x\in\U$. Now consider the triangle $a_k\lora e\lora x\lora \S a_k$. Since $a_k\in\ind\U$ and $x\in\U$, by Claim 1 we conclude that $e\in\U$.
\end{proof}

\subsection{Verdier quotients}

Let $\D$ be a thick subcategory of $\T$, we recall how to obtain the Verdier quotient $\T/\D$. We refer to \cite[Section 4]{K} for a detailed description.

We consider $S$ the class of morphisms $f\colon t_1\to t_2$ of $\T$ which extend to triangles of the form $t_1\overset{f}{\lora} t_2\lora d\lora \S t_1$ with $d\in\D$. The category $\T/\D$ has
\begin{itemize}
\item as objects exactly the same objects of $\T$;
\item as morphisms the equivalence classes of left fractions, see \cite[Section 3.1]{K};
\item the quotient functor $Q\colon \T\to \T/\D$ which acts as the identity on objects and makes the morphisms in $S$ invertible, and is universal with this property.
\end{itemize}

The category $\T/\D$ has a triangulated structure which consists of
\begin{itemize}
\item the shift functor $\S\colon\T/\D\to\T/\D$ induced by the shift functor $\S\colon \T\to \T$; 
\item triangles given by isomorphic copies of the images of the triangles of $\T$ after $Q$. 
\end{itemize}
With this triangulated structure, the quotient functor $Q$ is a triangulated functor, i.e. it is an additive functor commuting with $\S$ and sending triangles of $\T$ to triangles of $\T/\D$, see \cite[Lemma 4.3.1]{K}. From \cite[Proposition 4.6.2]{K} we recall that
\begin{itemize}
\item a morphism $f\colon t_1\to t_2$ in $\T$ is such that $Q(f) = 0$ in $\T/\D$ if and only if $f = hg$ for some $g\colon t_1\to d$ and $h\colon d\to t_2$ with $d\in\D$; 
\item an object $t\in\T$ is such that $Q(t)\cong 0$ in $\T/\U$ if and only if $t\in\D$.
\end{itemize}

Let $\U\subseteq\T$ and $\X\subseteq\T/\D$. The \emph{essential image} of $\U$ after $Q$, and the \emph{preimage} of $\X$ after $Q$, are respectively
\begin{align*}
Q(\U) &  = \{x\in\T/\D\mid x\cong Q(u) \text{ in $\T/\D$ for some $u\in\U$}\}\text{ and}\\
Q^{-1}(\X) & = \{t\in\T\mid Q(t)\cong x \text{ in $\T/\D$ for some $x\in\X$}\}.
\end{align*}

We have the following generalisation of \cite[Proposition 2.3.1]{V}.

\begin{proposition}\label{proposition bijection between extension closed subcategories in the verdier quotient}
Let $\T$ be a triangulated category and $\D$ be a thick subcategory. The following is an inclusion preserving bijection.
\begin{align*}
\left\{ \parbox{7cm}{\centering Extension-closed additive subcategories $\U\subseteq\T$ such that $\D\subseteq\U$} \right\} 
& \longleftrightarrow\left\{ \parbox{4.75cm}{\centering Extension-closed additive subcategories of $\T/\D$} \right\}\\
\U & \longmapsto Q(\U)\\
Q^{-1}(\X) & \longmapsfrom \X
\end{align*}
\end{proposition}

The argument in \cite{V} is in part not applicable with our assumptions when checking that the maps are well defined. Therefore, we provide an argument for this statement. Before doing so, we have the following lemma, which is included in the argument of \cite[Proposition 2.3.1]{V}. Our assumptions are more general than those of \cite{V}, but the argument still applies.

\begin{lemma}\label{lemma Q(X) closed under isomorphisms}
Let $\D$ be a thick subcategory of $\T$ and $\U$ be an extension-closed additive subcategory of $\T$ containing $\D$. If $t\in \T$ and $u\in\U$ are such that $Q(t)\cong Q(u)$ in $\T/\D$, then $u\in\U$. 
\end{lemma}

\begin{proof}[Proof of Proposition \ref{proposition bijection between extension closed subcategories in the verdier quotient}]
We check that the maps are well defined. To show that the two maps are mutually inverse we can proceed as in the argument of \cite[Proposition 2.3.1]{V}. Let $\U$ be an extension-closed additive subcategory of $\T$ containing $\D$. It is straightforward to see that $Q(\U)$ is closed under isomorphism, $0\in Q(\U)$, and that $Q(\U)$ is closed under direct sums. Moreover, by Lemma \ref{lemma Q(X) closed under isomorphisms}, it is straightforward to check that $Q(\U)$ is closed under direct summands. 

Now we show that $Q(\U)$ is extension-closed. Consider a triangle in $\T/\D$
\[
\text{(T)} \quad x_1\lora y\lora x_2\lora \S x_1
\]
with $x_1,x_2\in Q(\U)$. Then there is a triangle $a\lora e\lora b\lora \S a$ in $\T$ whose image under $Q$ is isomorphic to the triangle (T) in $\T/\D$, see the proof of \cite[Lemma 4.3.1]{K}. Thus, in $\T/\D$ we have the isomorphisms $Q(a)\cong x_1$, $Q(b)\cong x_2$ and $Q(e)\cong y$. Since $x_1,x_2\in Q(\U)$, we have that there exist $u_1,u_2\in\U$ such that $x_1\cong Q(u_1)$ and $x_2\cong Q(u_2)$. Then, by Lemma \ref{lemma Q(X) closed under isomorphisms}, we have that $a,b\in\U$. Since $\U$ is extension-closed, we obtain that  $e\in\U$ and as a consequence $y\cong Q(e)\in Q(\U)$. Thus, the map $\U\mapsto Q(\U)$ is well defined.
	
Let $\X$ be an extension-closed additive subcategory of $\T/\D$. We check that $Q^{-1}(\X)$ is an additive subcategory of $\T$. It is straightforward to see that $0\in Q^{-1}(\X)$, $Q^{-1}(\X)$ is closed under isomorphisms, direct sums, direct summands, and that $\D\subseteq Q^{-1}(\X)$. Now we show that $Q^{-1}(\X)$ is extension-closed. Consider a triangle $a\lora e\lora b\lora \S a$ in $\T$ with $a,b\in Q^{-1}(\X)$. Then its image under $Q$ is a triangle $Q(a)\lora Q(e)\lora Q(b)\lora \S Q(a)$ in $\T/\D$ with $Q(a), Q(b)\in\X$. Since $\X$ is extension-closed, then $Q(e)\in\X$. As a consequence $e\in Q^{-1}(\X)$. Hence, the map $\Y\mapsto Q^{-1}(\Y)$ is well defined.

Finally, from the definitions of $Q(\U)$ and $Q^{-1}(\X)$, it is straightforward to check that the maps $\U\mapsto Q(\U)$ and $\X\mapsto Q^{-1}(\X)$ preserve inclusion.
\end{proof}

\subsection{Non-crossing partitions}\label{section non-crossing partitions}

Let $k$ be a positive integer. Consider the unit circle $S^1$ with anticlockwise orientation, and a finite set  of elements of $S^1$, which we label as $\{1,\dots,k\} = [k]$, with the cyclic order $1<2<\cdots<k<1$.

A \emph{non-crossing partition} of $[k]$ is a partition $\P$ of $[k]$ such that for any $i_1,i_2,j_1,j_2\in[k]$ which are in cyclic order $i_1<j_1<i_2<j_2<i_1$, if $i_1,i_2\in B$ and $j_1,j_2\in C$ for some $B,C\in\P$, then $B = C$. If $\P$ is a non-crossing partition, its elements are called \emph{blocks}.

\begin{figure}[ht]
\centering
\includegraphics[height = 4cm]{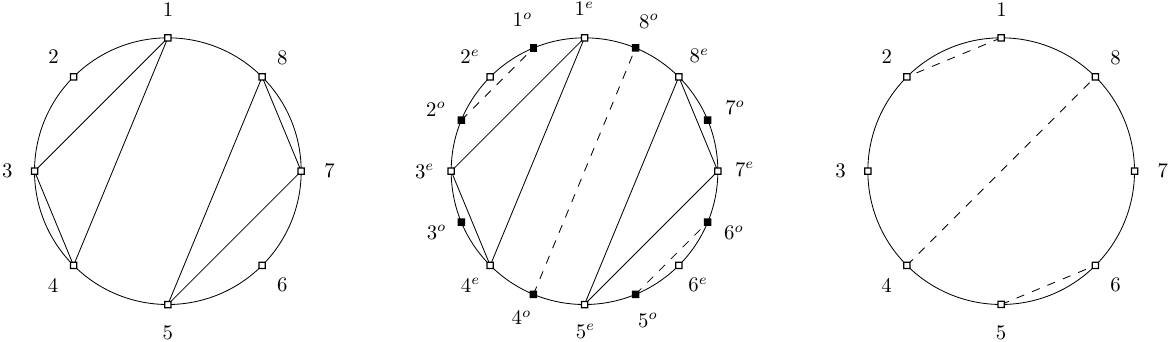}
\caption{On the left $\P = \{\{1,3,4\}, \{2\}, \{5,7,8\},\{6\}\}$ is a non-crossing partition of $[8]$, and $\P^c = \{\{1,2\},\{3\},\{4,8\},\{5,6\},\{7\}\}$ on the right is its Kreweras complement.}
\label{figure kreweras complement}
\end{figure}

The \emph{Kreweras complement}, $\P^c$, of a non-crossing partition $\P$ of $[k]$ is obtained as follows, see Figure \ref{figure kreweras complement} for an illustration.
\begin{enumerate}
\item Double the elements of $[k]$ to get the set $[k^e]\cup[k^o] = \{1^e,1^o,\dots,k^e,k^o\}$ with cyclic order $1^e<1^o<\dots<k^e<k^o<1^e$.
\item Define $\P^e$ as the non-crossing partition of $[k^e]$ which consists on $\P$. 
\item Complete $\P^e$ to a \emph{serrée} (\emph{dense}) non-crossing partition $\P^{e}\cup\P^o$ of $[k^e]\cup[k^o]$, see \cite[p. 338]{K1}.
\item Define $\P^c$ as $\P^o$ and relabel the elements of $[k^o]$ as $1,\dots,k$. 
\end{enumerate}

There is a partial order on the set of non-crossing partitions of $[k]$ given as follows: $\P\leq \P'$ if for each block $B$ of $\P$ there exists a block $B'$ of $\P'$ such that $B\subseteq B'$. The non-crossing partitions of $[k]$ form a \emph{lattice}: each pair of non-crossing partitions has a least upper bound $\P\vee\P'$, called \emph{join}, and a greatest lower bound $\P\wedge\P'$, called \emph{meet}. We refer to \cite[Section 2]{K1} for more details.

\begin{remark}\label{remark join and meet kreweras complement}
Let $\P$ and $\P'$ be non-crossing partitions of $[k]$. Then $(\P\vee\P')^c = \P^c\wedge\P'^c$. Indeed, it is straightforward to check that $\P\leq\P'$ if and only if $\P'^c\leq\P^c$. Therefore, the Kreweras complement of the least upper bound of $\P$ and $\P'$ is equal to the greatest lower bound of the Kreweras complements of $\P$ and $\P'$.
\end{remark}

\section{The categories $\mathcal{C}_m$ and $\overline{\mathcal{C}}_m$}\label{section the categories Cm and ovlCm}

In this section we recall the Igusa--Todorov discrete cluster category $\C_m$, introduced in \cite{IT}, the Paquette--Y\i ld\i r\i m completion $\ovl\C_m$, introduced in \cite{PY}, and their geometric models.

\subsection{The $\infty$-gons $\mathcal{Z}_m$, $\mathcal{Z}_{2m}$, and $\overline{\mathcal{Z}}_m$}\label{section infinity-gons}

\begin{figure}[ht]
\centering
\includegraphics[height = 4cm]{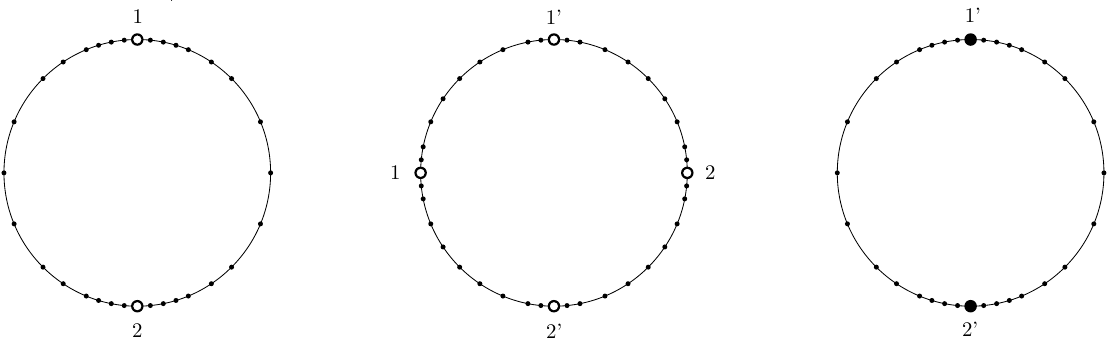}
\caption{On the left the $\infty$-gon $\Zc_2$, in the centre $\Zc_{4}$, and on the right $\ovl\Zc_2$. The white circles denote the accumulation points of $\Zc_2$ and $\Zc_{4}$, the black circles denote the accumulation points of $\ovl\Zc_2$.}
\label{figure infinity-gons}
\end{figure}

We consider the unit circle $S^1$ with anticlockwise orientation, endowed with the usual topology. Given a positive integer $m$, the $\infty$-gon $\Zc_m$ is an infinite discrete subset of $S^1$ consisting of $m$ copies of $\Z$ embedded in $S^1$ with $m$ two-sided accumulation points, see Figure \ref{figure infinity-gons}. We denote the accumulation points of $\Zc_m$ by $\{1,\dots,m\} = [m]$. Given $p\in[m]$, we denote by $\Z^{(p)}$ all the elements of $\Zc_m$ which belong to the $p$-th copy of $\Z$. The accumulation points are in cyclic order $1<\cdots<m<1$. If $p\in[m]$ is an accumulation point, we denote the successor and the predecessor of $p$ with respect to the cyclic order by $p^+$ and $p^-$. We also regard $[m]$ as a totally ordered set $1<\cdots<m$. This total order induces a total order $\leq$ on $\Zc_m\cup[m]$. 

We can define intervals in $\Zc_m$. Given $x,y\in\Zc_m\cup[m]$ we denote
\[
[x,y) = 
\begin{cases}
\{z\in\Zc_m\mid x\leq z < y\} & \text{if $x\leq y$, and}\\
\{z\in\Zc_m\mid z\leq x \text{ or } z>y\} & \text{otherwise.}
\end{cases}
\]
Similarly, we can define the intervals $(x,y]$, $(x,y)$, and $[x,y]$. Since the set $\Zc_m$ is discrete, for each $z\in\Zc_m$ there exists a predecessor $z-1$ and a successor $z+1$. 

\begin{definition}
A pair $x = (x_1,x_2)$ of elements of $\Zc_m$ is called \emph{arc} if $x_2\geq x_1+2$, and in that case $x_1$ and $x_2$ are called \emph{endpoints} or \emph{coordinates} of $x$. Given two arcs $x = (x_1,x_2)$ and $y = (y_1,y_2)$ of $\Zc_m$, we say that $x$ and $y$ \emph{cross} if $x_1<y_1<x_2<y_2$ or $y_1<x_1<y_2<x_2$. Given $p,q\in[m]$ with $p\leq q$, we define 
\[
\Z^{(p,q)} = \left\{(x_1,x_2) \text{ is an arc of $\Zc_m$ }\middle| \text{ $x_1\in\Z^{(p)}$ and $x_2\in\Z^{(q)}$}\right\}.
\] 
\end{definition}

From $\Zc_m$ we define another $\infty$-gon $\ovl\Zc_m$. To this end, we take an intermediate step by considering the $\infty$-gon $\Zc_{2m}$. We re-label the accumulation points of $\Zc_{2m}$ as $1',1,\dots,m',m$, see Figure \ref{figure infinity-gons}. The set of accumulation points $[m']\cup[m]$ has cyclic order $1'<1<\cdots<m'<m<1'$ and a total order $1'<1<\cdots<m'<m$, which induces a total order on $\Zc_{2m}$. The notions of interval, successor, predecessor, arc, are the same as for the set $\Zc_m$. 

On $\Zc_{2m}$ we define an equivalence relation $\sim$ as follows. For each $x,y\in\Zc_{2m}$ we have that 
\[
x\sim y\text{ if and only if } x = y\text{ or } x,y\in\Z^{(p)} \text{ for some } p\in[m'].
\] 
Consider $x\in\Zc_{2m}$, we sometimes denote the equivalence class of $x$ by $\ovl x$. If $x\in\Z^{(p)}$ for some $p\in[m']$, we identify $\ovl x = p$ with an abuse of notation. We define the set $\ovl\Zc_m = \Zc_{2m}/\sim$ and we observe that $\ovl\Zc_m$ can be regarded as the set $\Zc_m\cup[m']$. The total order on $\Zc_{2m}$ induces a total order on $\ovl\Zc_m$.

Given a point $z\in\ovl\Zc_m = \Zc_m\cup[m']$, we define the successor $z+1$ of $z$ as
\[
z+1 = 
\begin{cases}
\text{the successor of $z$ in $\Zc_m$} & \text{ if $z\in \Zc_m$,}\\
z & \text{ if $z\in[m']$.}
\end{cases}
\] 
We can define $z-1$ analogously. The notions of arc of $\ovl\Zc_m$ and of crossing arcs are the same as those for $\Zc_m$. Given $p,q\in[m']\cup[m]$ we define the following sets
\[
C^{(p)} = 
\begin{cases}
\Z^{(p)} & \text{ if $p\in[m]$,}\\
\{p\} & \text{ if $p\in[m']$}
\end{cases}\text{ and } C^{(p,q)} = \left\{(x_1,x_2) \text{ is an arc of $\ovl\Zc_m$ }\middle| \text{ $x_1\in C^{(p)}$ and $x_2\in C^{(q)}$}\right\}.
\] 
The following notation will be useful later. Given $x$ and $y$ both elements of $\Zc_m$, or both elements of $\ovl\Zc_m$, with $x_2\geq x_1+2$ or $x_1\geq x_2+2$, we denote
\[
|x_1,x_2| =
\begin{cases}
(x_1,x_2) & \text{if $x_1<x_2$,}\\
(x_2,x_1) & \text{if $x_2<x_1$.}
\end{cases}
\]

\subsection{Geometric models}\label{section geometric models}

Given the $\infty$-gon $\Zc_m$ and a field $\K$, the category $\C_m$ was defined in \cite{IT}. This is a $\K$-linear, $\Hom$-finite, Krull--Schmidt triangulated category. We denote its shift functor by $\S\colon\C_m\to\C_m$. Moreover, $\C_m$ is 2-Calabi--Yau, i.e. $\S^2$ is a Serre functor. We recall some properties of $\C_m$.

\begin{itemize}
\item There is a bijection between the isoclasses of indecomposable objects of $\C_m$ and the arcs of $\Zc_m$. We regard the indecomposable objects of $\C_m$ as arcs of $\Zc_m$, see \cite[Section 2.4.1]{IT}.
\item  Given $x = (x_1,x_2)\in\ind\C_m$ we have that $\S x = (x_1-1,x_2-1)$ by \cite[Lemma 2.4.3]{IT}.
\item Given $x,y\in\ind\C_m$, by \cite[Lemma 2.4.4]{IT} we have that
\[
\Hom_{\C_m}(x,y)\cong
\begin{cases}
\K & \text{if $x$ and $\S^{-1}y$ cross,}\\
0 & \text{otherwise.}
\end{cases}
\]
\end{itemize}

The completion $\ovl\C_m$ of $\C_m$ was defined in \cite{PY}, we recall its construction. Consider the set $\Zc_{2m}$ defined in Section \ref{section infinity-gons} and the associated category $\C_{2m}$. Define $\D$ as
\[
\D = \add\left\{\bigcup_{p\in[m']}\Z^{(p,p)}\right\}.
\]
The category $\D$ is a thick subcategory of $\C_{2m}$ and $\ovl\C_m$ is defined as the Verdier quotient $\C_{2m}/\D$. This is a $\K$-linear, $\Hom$-finite, Krull--Schmidt triangulated category. We denote the quotient functor as $\pi\colon \C_{2m}\to \C_{2m}/\D = \ovl\C_m$ and its shift functor by $\S\colon \ovl\C_m\to\ovl\C_m$ as for $\C_m$. We recall the following properties.

\begin{itemize}
\item The isoclasses of indecomposable objects of $\ovl\C_m$ are in bijection with the arcs of $\ovl\Zc_m$, see \cite[Corollary 3.11]{PY}. 
\item For any $x = (x_1,x_2)\in\ind\C_{2m}\setminus\ind\D$ the object $\pi x\in\ovl\C_m$ is indecomposable by \cite[Proposition 3.10]{PY} and can be regarded as the arc $\left(\ovl x_1,\ovl x_2\right)$ of $\ovl\Zc_m$.
\item Let $x\in\ind\ovl\C_m$, then there exists $x'\in\ind\C_{2m}$ such that $\pi x'\cong x$. Indeed, if $x = (x_1,x_2)$  with $x_1,x_2\in\ovl\Zc_m$, we can take $x_1',x_2'\in\Zc_{2m}$ such that $\ovl x_1' = x_1$ and $\ovl x_2' = x_2$, we define $x' = (x_1',x_2')\in\ind\C_{2m}$. 
\item Given $x = (x_1,x_2)\in\ind\ovl\C_m$ we have that $\S x = (x_1-1,x_2-1)$.
\item The Hom-spaces of $\ovl\C_m$ between indecomposable objects are at most one-dimensional. More precisely, we have the following proposition.
\end{itemize}

\begin{proposition}[{\cite[Proposition 3.14]{PY}}]\label{proposition Hom sets completion}
Let $x, y\in\ind\ovl{\C}_m$. Then $\Hom_{\ovl{\C}_m}(x,\S y)\cong \K$ if and only if one of the following statements holds.
\begin{itemize}
\item The arcs $x$ and $y$ cross.
\item The arcs $x$ and $y$ share exactly one endpoint $z\in[m']$, and we can reach $y$ by rotating $x$ in the anticlockwise direction about $z$.
\item The arcs $x$ and $y$ share both endpoints $z_1,z_2\in[m']$.
\end{itemize}
Otherwise $\Hom_{\ovl{\C}_m}(x,\S y) = 0$.
\end{proposition}

\begin{remark}\label{remark additive full subcategories of IT category are sets of arcs}
For both categories $\C_m$ and $\ovl\C_m$, we identify the indecomposable objects with arcs of $\Zc_m$ or $\ovl\Zc_m$, and the full additive subcategories with sets of arcs.
\end{remark}

From now on any subcategory of $\C_m$ or $\ovl\C_m$ is assumed to be additive and full.

\section{The AR quiver of $\overline{\mathcal{C}}_m$} \label{section quiver}

In this section we describe the AR quiver of $\ovl\C_m$. It is well known that $\ovl\C_m$ does not have a Serre functor, and therefore no almost split triangles in general. In this setting, by AR quiver we mean the quiver having as vertices the isoclasses of indecomposable objects of $\ovl\C_m$ and as arrows the irreducible morphisms between them.

\subsection{The coordinate system}\label{section coordinate systems}

Recall from \cite[Theorem 2.4.13]{IT} that the AR quiver of $\C_m$ consists of
\begin{itemize}
\item $m$ components of type $\Z A_{\infty}$, corresponding to the arcs of $\Z^{(p,p)}$ for $p\in[m]$,
\item $\binom{m}{2}$ components of type $\Z A_{\infty}^{\infty}$, corresponding to the arcs of $\Z^{(p,q)}$ for $p,q\in[m]$ with $p\neq q$.
\end{itemize}

\begin{figure}[ht]
\centering
\includegraphics[height = 3cm]{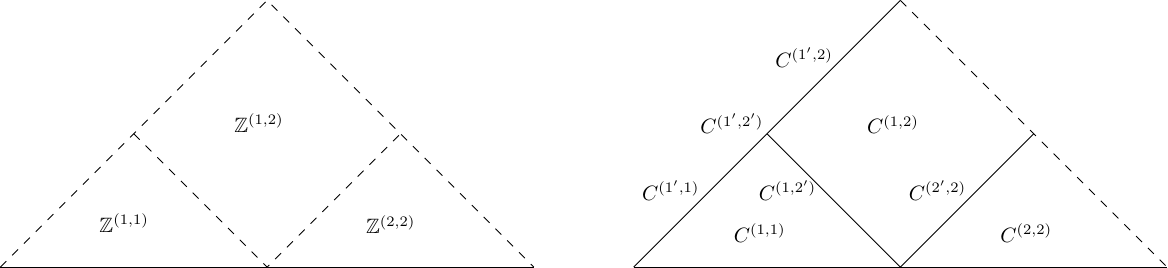}
\caption{On the left the AR quiver of $\C_2$, on the right the AR quiver of $\ovl\C_2$.}\label{figure quivers}
\end{figure}

For the category $\ovl\C_m$, we can arrange the isoclasses of the indecomposable objects of $\ovl\C_m$ into a coordinate system having
\begin{itemize}
\item $m$ components of type $\Z A_{\infty}$, corresponding to the arcs of $C^{(p,p)}$ for $p\in[m]$,
\item $\binom{m}{2}$ components of type $\Z A_{\infty}^{\infty}$, corresponding to the arcs of $C^{(p,q)}$ for $p,q\in[m]$ with $p\neq q$,
\item $\binom{m}{2}$ components of type $A_1$, corresponding to the arcs of $C^{(p,q)}$ for $p,q\in[m']$ with $p\neq q$,
\item $m^2$ components of type $A_{\infty}^{\infty}$, corresponding to the arcs of $C^{(p,q)}$ for $p,q\in[m']\cup[m]$ such that either $p\in[m']$ and $q\in[m]$, or $p\in[m]$ and $q\in [m']$.
\end{itemize}
Figure \ref{figure quivers} illustrates this coordinate system. In Proposition \ref{proposition irreducible morphisms} we describe the irreducible morphisms of $\ovl\C_m$ and thus show that the above describes the AR quiver of $\ovl\C_m$.

\subsection{Hom-hammocks}

Before describing the $\Hom$-hammocks of $\ovl\C_m$, we extend the definition of Hom-hammocks of $\C_m$, see \cite[Definition 2.1]{HJ}, from $m = 1$ to the general case $m\geq 1$.

\begin{definition}\label{definition hom hammocks IT category}
Let $a = (a_1,a_2)\in\ind\C_m$. We define 
\begin{align*}
H^+(a) & = \{(x_1,x_2)\in\ind\C_m\mid a_1\leq x_1\leq a_2-2 \text{ and } x_2\geq a_2\}\text{ and}\\
H^-(a) & = \{(x_1,x_2)\in\ind\C_m\mid x_1\leq a_1 \text{ and } a_1+2\leq x_2\leq a_2\}.
\end{align*}
\end{definition}

\begin{remark}
For $a,b\in\ind\C_m$, by \cite[Lemma 2.4.2]{IT} it follows that $\Hom_{\C_m}(a,b)\cong\K$ if and only if $b\in H^+(a)\cup H^-(\S^2 a)$, or equivalently, $a\in H^+(\S^{-2} b)\cup H^-(b)$.
\end{remark} 

We define the Hom-hammocks for the category $\ovl\C_m$ analogously to $\C_m$; Figure \ref{figure hammocks} provides an illustration.

\begin{definition}\label{definition hom-hammocks completion}
Let $a = (a_1,a_2)\in\ind\ovl\C_m$, and let $p,q\in[m']\cup[m]$ be such that $a\in C^{(p,q)}$. We define the \emph{$\Hom$-hammocks} $\ovl H^+(a)$ and $\ovl H^-(a)$ as follows.
\[
\ovl H^+(a) = 
\begin{cases}
\{(x_1,x_2)\in\ind\ovl\C_m\mid a_1\leq x_1\leq a_2-2 \text{ and } x_2\geq a_2\} & \text{ if } q\in[m],\\
\{(x_1,x_2)\in\ind\ovl\C_m\mid a_1\leq x_1 < a_2 \text{ and } x_2\geq a_2\} & \text{ if } q\in[m'].
\end{cases}
\]
\[
\ovl H^-(a) = 
\begin{cases}
\{(x_1,x_2)\in\ind\ovl\C_m\mid 1'\leq x_1\leq a_1 \text{ and } a_1+2\leq x_2\leq a_2\} & \text{ if } p,q\in[m],\\
\{(x_1,x_2)\in\ind\ovl\C_m\mid 1'\leq x_1\leq a_1 \text{ and } a_1+2\leq x_2< a_2\} &  \text{ if $p\in[m]$ and $q\in[m']$,}\\
\{(x_1,x_2)\in\ind\ovl\C_m\mid 1'\leq x_1 < a_1 \text{ and } a_1\leq x_2\leq a_2\} & \text{ if } p\in[m']\text{ and } q\in[m],\\
\{(x_1,x_2)\in\ind\ovl\C_m\mid 1'\leq x_1< a_1 \text{ and } a_1\leq x_2< a_2\} & \text{ if } p,q\in[m'].\\
\end{cases}
\]
\end{definition}

We can reformulate Proposition \ref{proposition Hom sets completion} as follows.

\begin{proposition}
Let $a,b\in\ind\ovl\C_m$. Then $\Hom_{\ovl\C_m}(a,b)\cong \K$ if and only if $b\in \ovl H^+(a)\cup \ovl H^-(\S^{2}a)$.
\end{proposition}

Since $\ovl\C_m$ is not 2-Calabi-Yau, for $a,b\in\ind\ovl\C_m$ in general $b\in \ovl H^+(a)\cup\ovl H^-(\S^2 a)$ is not equivalent to $a\in \ovl H^+(\S^{-2}b)\cup \ovl H^-(b)$. Therefore, we also define the \emph{reverse $\Hom$-hammocks} $\ovl I^+$ and $\ovl I^-$, for which $b\in \ovl H^+(a)\cup\ovl H^-(\S^2 a)$ if and only if $a\in \ovl I^+(\S^{-2}b)\cup \ovl I^-(b)$.

\begin{definition}
Let $a\in\ind\ovl\C_m$, and let $p,q\in[m']\cup[m]$ be such that $a\in C^{(p,q)}$. We define the \emph{reverse $\Hom$-hammocks} $\ovl I^+(a)$ and $\ovl I^-(a)$ as follows.
\[
\ovl I^+(a) = 
\begin{cases}
\{(x_1,x_2)\in\ind\ovl\C_m\mid a_1\leq x_1\leq a_2-2 \text{ and } x_2\geq a_2\} & \text{ if } p,q\in[m],\\
\{(x_1,x_2)\in\ind\ovl\C_m\mid a_1\leq x_1\leq a_2 \text{ and } x_2 > a_2\} & \text{ if } p\in[m]\text{ and } q\in[m'],\\
\{(x_1,x_2)\in\ind\ovl\C_m\mid a_1< x_1\leq a_2-2 \text{ and } x_2\geq a_2\} & \text{ if } p\in[m']\text{ and } q\in[m],\\
\{(x_1,x_2)\in\ind\ovl\C_m\mid a_1< x_1\leq a_2 \text{ and } x_2> a_2\} & \text{ if } p,q\in[m'].\\
\end{cases}
\]
\[
\ovl I^-(a) = 
\begin{cases}
\{(x_1,x_2)\in\ind\ovl\C_m\mid 1'\leq x_1\leq a_1 \text{ and } a_1+2\leq x_2\leq a_2\} & \text{ if } p\in[m],\\
\{(x_1,x_2)\in\ind\ovl\C_m\mid 1'\leq x_1\leq a_1 \text{ and } a_1< x_2\leq a_2\} & \text{ if } p\in[m'].
\end{cases}
\]
\end{definition}

Note that in general $\ovl H^+(a) \neq \ovl I^+(a)$ and $\ovl H^-(a) \neq \ovl I^-(a)$. We can reformulate Proposition \ref{proposition Hom sets completion} as follows.

\begin{proposition}
Let $a,b\in\ind\ovl\C_m$. Then $\Hom_{\ovl\C_m}(a,b)\cong \K$ if and only if $a\in \ovl I^+(\S^{-2}b)\cup \ovl I^-(b)$.
\end{proposition}

\subsection{Factorization properties}

Here we study the factorization properties of the morphisms of $\ovl\C_m$. First we recall the factorization properties of $\C_m$ which will be useful later. We say that a morphism $f\colon a\to b$ in $\C_m$ \emph{factors through an object} $c\in\ind\C_{m}$ if there exist $g\colon a\to c$ and $h\colon d\to c$ such that $f = hg$. We say that a morphism $f\colon a\to b$ in $\C_m$ \emph{factors through $\D$} if it factors through some $d\in\ind\D$. 

\begin{lemma}[{\cite[Lemma 2.4.2]{IT}}]\label{lemma factorization properties in Cm}
Let $a,b,c\in\ind\C_m$. Assume that one of the following statements holds.
\begin{enumerate}
\item $b\in H^+(a)$ and $c\in H^+(a)\cap H^+(b)$.
\item $b\in H^+(a)$ and $c\in H^-(\S^2 a)\cap H^-(\S^2 b)$.
\item $b\in H^-(\S^2 a)$ and $c\in H^-(\S^2 a)\cap H^+(b)$.
\end{enumerate}
Then any morphism $a\to c$ in $\C_m$ factors through $b$.
\end{lemma} 

The following lemma will be useful for proving Proposition \ref{proposition factorization properties}.

\begin{figure}[ht]
\centering
\includegraphics[width = 15cm]{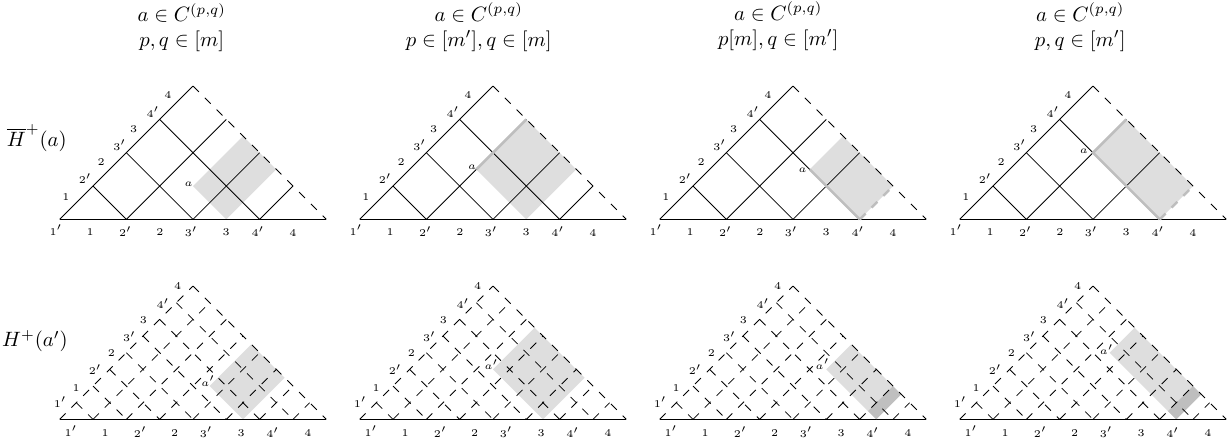}\\[0.5cm]
\includegraphics[width = 15cm]{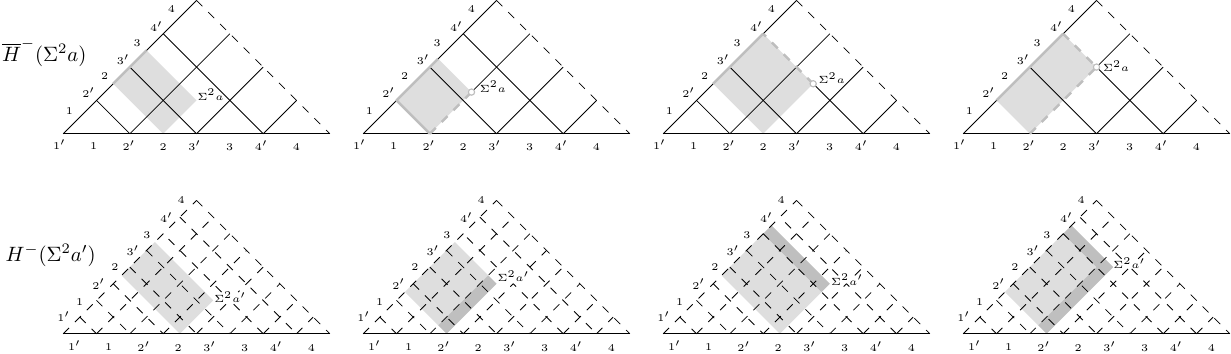}
\caption{Illustration of Lemma \ref{lemma hammocks}. The darker area in $H^+(a')$ or in $H^+(\S^2 a')$, if present, denotes the objects $b'\in H^+(a')\cup H^-(\S^2 a')$ such that any non-zero morphism $a'\to b'$ factors through $\D$. Whenever the darker area is not present there are no such objects in $H^+(a')\cup H^-(\S^2 a')$.}
\label{figure hammocks}
\end{figure}

\begin{lemma}\label{lemma hammocks}
Let $a,b\in\ind\ovl\C_m$ be such that $\Hom_{\ovl\C_m}(a,b)\cong\K$, and let $a'\in\ind\C_m$ be such that $\pi a'\cong a$. The following statements hold.
\begin{enumerate}
\item If $b\in\ovl H^+(a)$, then there exists $b'\in\ind\C_{2m}$ such that $\pi b'\cong b$, $b'\in  H^+(a')$, and any non-zero morphism $a'\to b'$ in $\C_{2m}$ does not factor through $\D$.
\item If $b\in\ovl H^-(\S^2 a)$, then there exists $b'\in\ind\C_{2m}$ such that $\pi b'\cong b$, $b'\in H^-(\S^2 a')$, and any non-zero morphism $f'\colon a'\to b'$ in $\C_{2m}$ does not factor through $\D$.
\end{enumerate}
\end{lemma}
\begin{proof}
We show statement (1), statement (2) is analogous. Assume that $a = (a_1,a_2)\in C^{(p,q)}$ with $p\in[m]$ and $q\in[m']$, the other cases are similar. We write $a' = (a_1',a_2')\in\Z^{(p,q)}$. Since $b = (b_1,b_2)\in \ovl H^+(a)$, we have that $a_1\leq b_1 < a_2$ and $b_2\geq a_2$, see Figure \ref{figure intervals}.

\begin{figure}[ht]
\centering
\includegraphics[height = 6cm]{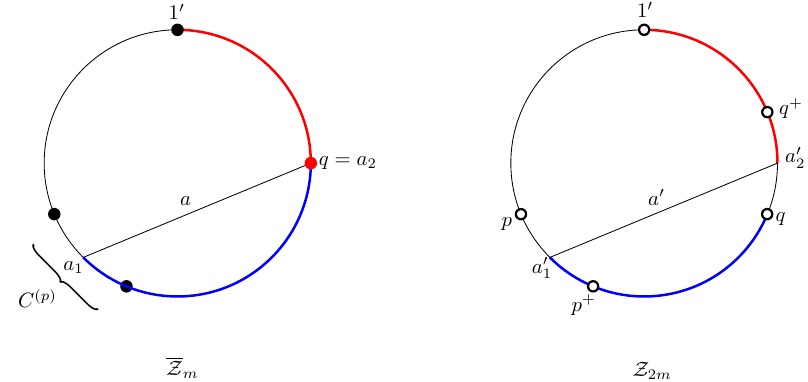}
\caption{The element $b_1\in\ovl\Zc_m$ is such that $a_1\leq b_1< a_2$ and the element $b_2\in\ovl\Zc_m$ is such that $b_2\geq a_2$. We can find $b_1'\in\Zc_{2m}$ such that $\ovl b_1' = b_1$ and $a_1'\leq b_1'$. Similarly for $b_2'$ in $\Zc_{2m}$.}
\label{figure intervals}
\end{figure}

It is straightforward to check that there exists $b_1'$ such that $\ovl b_1' = b_1$ and $a_1'\leq b_1'<q$, and there exists $b_2'$ such that $\ovl b_2' = b_2$ and $b_2'\geq a_2'$. Therefore, $b'\in H^+(a')$ and, since $b_1'\notin \Z^{(q)}$, any non-zero morphism $a'\to b'$ does not factor through $\D$, see Figure \ref{figure hammocks}.
\end{proof}

The following lemma is dual to the lemma above, and will be useful for proving Lemma \ref{lemma A aisle t-structure}.

\begin{lemma}\label{lemma reverse hammocks}
Let $a,b\in\ind\ovl\C_m$ be such that $\Hom_{\ovl\C_m}(b,a)\cong\K$, and let $a'\in\ind\C_m$ be such that $\pi a'\cong a$. The following statements hold.
\begin{enumerate}
\item If $b\in\ovl I^-(a)$, then there exists $b'\in\ind\C_{2m}$ such that $\pi b'\cong b$, $b'\in  H^-(a')$, and any non-zero morphism $b'\to a'$ in $\C_{2m}$ does not factor through $\D$.
\item If $b\in\ovl I^+(\S^{-2} a)$, then there exists $b'\in\ind\C_{2m}$ such that $\pi b'\cong b$, $b'\in H^+(\S^{-2} a')$, and any non-zero morphism $f'\colon b'\to a'$ in $\C_{2m}$ does not factor through $\D$.
\end{enumerate}
\end{lemma}

Now we have the factorization properties of $\ovl\C_m$.

\begin{proposition}\label{proposition factorization properties}
Let $a,b,c\in\ind\ovl\C_m$. Assume that one of the following statements holds.
\begin{enumerate}
\item $b\in \ovl H^+(a)$ and $c\in \ovl H^+(a)\cap \ovl H^+(b)$.
\item $b\in \ovl H^+(a)$ and $c\in \ovl H^-(\S^{2}a)\cap \ovl H^-(\S^{2}b)$.
\item $b\in \ovl H^-(\S^{2}a)$ and $c\in \ovl H^-(\S^{2}a)\cap \ovl H^+(b)$.
\end{enumerate}	
Then any morphism $a\to c$ in $\ovl\C_m$ factors through $b$. 
\end{proposition}
\begin{proof}
We prove statement (1), statements (2) and (3) are analogous. Fix $a'\in\ind\C_{2m}$ such that $\pi a'\cong a$. Since $b\in\ovl H^+(a)$, by Lemma \ref{lemma hammocks} there exists $b'\in\ind\C_{2m}$ such that $\pi b'\cong b$, $b'\in H^+(a')$, and any non-zero morphism $a'\to b'$ does not factor through $\D$. Fix such $b'$, since $c\in \ovl H^+(b)$, then there exists $c'\in\ind\C_{2m}$ such that $\pi c'\cong c$, $c'\in H^+(b')$, and any non-zero morphism $b'\to c'$ does not factor through $\D$. We show that $c'\in H^+(a')\cap H^+(b')$.
	
We denote $a = (a_1,a_2)$, $a' = (a_1',a_2')$, $c = (c_1,c_2)$, and $c' = (c_1',c_2')$. Assume that $c'\notin H^+(a')$, then $c'_1\geq a_2'-1$. It is straightforward to check that as a consequence $c_1\geq a_2-1$. Then $c\not\in\ovl H^+(a)$ and we have a contradiction. Therefore $c'\in H^+(a')\cap H^+(b')$. Now, if there exists a non-zero morphism $f'\colon a'\to c'$ which factors through $\D$, then $a_2',c_1' \in\Z^{(q)}$, see Figure \ref{figure hammocks}. This implies that $c_1 = q = a_2$, and then $c\notin\ovl H^+(a)$ giving a contradiction. Then any non-zero morphism $f'\colon a'\to c'$ does not factor through $\D$.

Since $b'\in H^+(a')$ and $c'\in H^+(a')\cap H^+(b')$, by Lemma \ref{lemma factorization properties in Cm} there exist $h'\colon a'\to b'$ and $g'\colon b'\to c'$ such that $f' = g'h'$, and then $\pi f' = \pi(g')\pi(h')$. Since $f'$ does not factor through $\D$, we have that $\pi f' \neq 0$. Now consider a non-zero morphism $f\colon a\to c$ in $\ovl\C_m$. Since $\Hom_{\ovl\C_m}(a,c)\cong\K$, we have that $f = \l\pi f$ for some $\l\in\K^*$, and then $f = \l\pi f' = \l\pi(g')\pi(h')$. This concludes the argument.
\end{proof}

\subsection{Irreducible morphisms}\label{section irreducible morphisms}

In this section we describe the irreducible morphisms of $\ovl\C_m$. From Section \ref{section coordinate systems} we already know that the isoclasses of indecomposable objects of $\ovl\C_m$ are in bijection with the arcs of $\ovl\Zc_m$ and that they can be arranged in a coordinate system.

We recall that, from \cite[Lemma 2.4.11]{IT}, for an object $(a_1,a_2)\in\ind\C_{2m}$, the irreducible morphisms of $\C_m$ are exactly the non-zero morphisms of the form $(a_1,a_2)\to (a_1,a_2+1)$ or $(a_1,a_2)\to (a_1+1,a_2)$, provided that $(a_1+1,a_2)$ is still an arc, i.e. $a_2\geq a_1+3$.

\begin{proposition}\label{proposition irreducible morphisms}
Let $a = (a_1,a_2),b = (b_1,b_2)\in\ind\ovl\C_m$. Assume that $a,b\in C^{(p,q)}$ for some $p,q\in[m']\cup[m]$ and that one of the following conditions holds.
\begin{enumerate}
\item $p\in[m']$, $q\in[m]$ and $(b_1,b_2) = (a_1,a_2+1)$.
\item $p\in[m]$, $q\in[m']$ and $(b_1,b_2) = (a_1+1,a_2)$.
\item $p,q\in[m]$ and $(b_1,b_2) = (a_1,a_2+1)$ or $(b_1,b_2) = (a_1+1,a_2)$.	
\end{enumerate}
Then any non-zero morphism $a\to b$ is irreducible. Moreover, there are no other irreducible morphisms in $\ovl\C_m$ between indecomposable objects.
\end{proposition}
\begin{proof}
First we show that if any of the conditions (1), (2) and (3) holds, then any non-zero morphism $f\colon a \to b$ is irreducible. Assume that condition (1) holds, for the other cases we can proceed analogously. Consider a non-zero morphism $f\colon a\to b$ and note that, since $(a_1,a_2+1)\not\cong (a_1,a_2)$, $f$ is not a split monomorphism nor a split epimorphism. Assume that $f = hg$ for some $g\colon a\to c$, $h\colon c\to b$, and $c\in\ovl\C_m$. Since the $\Hom$-spaces are one dimensional, we can assume that $c\in\ind\ovl\C_m$. We show that $g$ is a split monomorphism or $h$ is a split epimorphism. Note that 
\[c\in \left(\ovl H^+(a)\cup\ovl H^-(\S^2 a)\right)\cap\left( \ovl I^-(b)\cup \ovl I^+(\S^{-2}b)\right).
\]
Assume that $c\not\cong a$ and $c\not\cong b$, then $g\colon a\to c$ factors as $g = lf$ with $l\colon b\to c$, see Figure \ref{figure illustration argument proposition irreducible morphisms}. Since $0\neq f = hg = hlf$, it follows that $hl\colon b\to b$ is non-zero and $hl = \l 1_b$ for some $\l\in\K^*$. This implies that $b\cong c$, which gives a contradiction with our assumption. We conclude that $c\cong a$ or $c\cong b$, i.e. $f\colon a\to b$ is irreducible.

\begin{figure}[ht]
\centering
\includegraphics[height = 4cm]{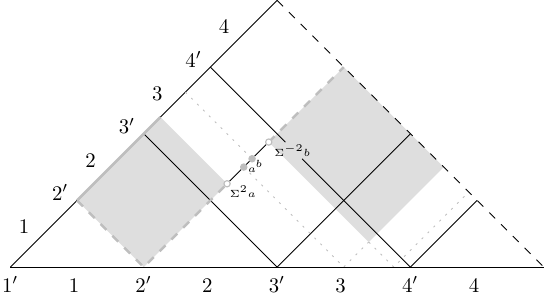}	
\caption{The object $c$ belongs to the grey area.}
\label{figure illustration argument proposition irreducible morphisms}
\end{figure}

Now, consider $a = (a_1,a_2),b = (b_1,b_2)\in\ind\ovl\C_m$ and a non-zero morphism $f\colon a\to b$. We show that if $f$ is irreducible then it has to be of the form listed in the statement. Let $p,q\in[m']\cup[m]$ be such that $a\in C^{(p,q)}$. Assume that $p\in[m']$ and $q\in[m]$, the other cases are analogous. Note that if $b_2\neq a_2$ then, from Proposition \ref{proposition factorization properties}, $f$ factors through the irreducible morphism $a\to (a_1,a_2+1)$, and then $f$ is not irreducible unless $(b_1,b_2) = (a_1,a_2+1)$. If $b_2 = a_2$, then consider the object $c = (b_2-1,a_2)$ and the non-zero morphisms $g\colon a\to c$ and $h\colon c\to b$. From Proposition \ref{proposition factorization properties} we have that $f = hg$, and then $f$ is not irreducible. We can conclude that if $f$ is irreducible then $(b_1,b_2) = (a_1,a_2+1)$.
\end{proof}

From Proposition \ref{proposition irreducible morphisms} we obtain that $\ovl\C_m$ does not have a Serre functor. This is a well known fact, but we could not find an argument in the literature. We provide the argument below for the convenience of the reader.

\begin{corollary}
The category $\ovl\C_m$ does not have a Serre functor.
\end{corollary}
\begin{proof}
By \cite[Proposition I.2.3]{RV}, it is equivalent to show that $\ovl\C_m$ does not always have almost split triangles. Let $a = (a_1,a_2)\in\ind\ovl\C_m$ with $a_1,a_2\in[m']$.  By Proposition \ref{proposition irreducible morphisms}, there are no irreducible morphisms from and to $a$, and therefore there are no almost split triangles of the form $a\lora e\lora {\tau}^{-1} a\lora \S a$ and $\tau a\lora e\lora a\lora \S\tau a$.
\end{proof}

\section{Precovering and preenveloping subcategories}\label{section precovering and preenveloping subcategories}

In this section we classify the precovering and preenveloping  subcategories of $\ovl\C_m$ using arc combinatorics. We also relate the precovering or preenveloping subcategories in $\ovl\C_m$ to their preimages in $\C_{2m}$ under the localisation functor $\pi\colon \C_{2m}\to \ovl\C_m$. 

In \cite{PY} the authors classified the functorially finite weak cluster-tilting subcategories of $\ovl\C_m$, i.e. the cluster-tilting subcategories, in terms of certain triangulations of the $\infty$-gon $\ovl\Zc_m$. This classification generalises \cite{F} for the case $m = 1$. After endowing $\ovl\C_m$ with a specific extriangulated structure, the cluster-tilting subcategories were also classified in \cite{CKP} in terms of a larger class of triangulations of $\ovl\Zc_m$. Here we classify subcategories of $\ovl\C_m$ which are just precovering or preenveloping.

\subsection{Precovering subcategories of $\C_m$}\label{section precovering subcategories IT category}
We start by recalling the classification of the precovering subcategories of $\C_m$ from \cite{GHJ}.

The following definition corresponds to \cite[Definition 0.7]{GHJ}, but we use a different formulation which is more convenient for our purposes. For the notation $|x_1,x_2|$ we refer to Section \ref{section infinity-gons}. 

\begin{definition}\label{definition pc conditions}
Let $\U$ be a subcategory of $\C_m$. We say that $\U$ satisfies the \emph{precovering conditions} ($\PC$ conditions for short) if it satisfies the following combinatorial conditions.
\begin{enumerate}
\item[($\PC 1$)] If there exists a sequence $\{(x_1^n,x_2^n)\}_n\subseteq \U\cap \Z^{(p,q)}$ for some $p,q\in[m]$ such that $p\neq q$ and the sequences $\{x_1^n\}_n$ and $\{x_2^n\}_n$ are strictly increasing, then there exist strictly decreasing sequences $\{y_1^n\}_n\subseteq\Z^{(p^+)}$ and $\{y_2^n\}_n\subseteq\Z^{(q^+)}$ such that  $\{|y_1^n,y_2^n|\}_n\subseteq \U$.
\item[($\PC 2$)] If there exists a sequence $\{(x_1^n,x_2^n)\}_n\subseteq \U\cap \Z^{(p,q)}$ for some $p,q\in[m]$ such that $p\neq q^+$ and the sequences $\{x_1^n\}_n$  and $\{x_2^n\}_n$ are respectively strictly decreasing and strictly increasing, then there exist strictly decreasing sequences $\{y_1^n\}_n\subseteq\Z^{(p)}$ and $\{y_2^n\}_n\subseteq\Z^{(q^+)}$ such that $\{|y_1^n,y_2^n|\}_n\subseteq \U$.
\item[($\PC 2'$)] If there exists a sequence $\{(x_1^n,x_2^n)\}_n\subseteq \U\cap\Z^{(p,q)}$ for some $p,q\in[m]$ such that $q\neq p^+$, $p\neq q$, and the sequences $\{x_1^n\}_n$ and $\{x_2^n\}_n$ are respectively strictly increasing and strictly decreasing, then there exist strictly decreasing sequences $\{y_1^n\}_n\subseteq\Z^{(p^+)}$ and $\{y_2^n\}_n\subseteq\Z^{(q)}$ such that $\{(y_1^n,y_2^n)\}_n\subseteq \U$.
\item[($\PC 3$)] If there exists a sequence $\{(x_1,x_2^n)\}_n\subseteq \U\cap\Z^{(p,q)}$ for some $p,q\in[m]$ such that the sequence $\{x_2^n\}_n$ is strictly increasing, then there exists a strictly decreasing sequence $\{y_2^n\}_n\subseteq\Z^{(q^+)}$ such that  $\{|x_1,y_2^n|\}_n\subseteq \U$.
\item[($\PC 3'$)] If there exists a sequence $\{(x_1^n,x_2)\}_n\subseteq\U\cap\Z^{(p,q)}$ for some $p,q\in[m]$ such that $p\neq q$ and the sequence $\{x_1^n\}_n$ is strictly increasing, then there exists a strictly decreasing sequence $\{y_1^n\}_n\subseteq\Z^{(p^+)}$ such that $\{(y_1^n,x_2)\}_n\subseteq \U$.
\end{enumerate}
\end{definition}

The conditions $(\PC 1)$, $(\PC 3)$, and $(\PC 3')$ correspond to condition $(\PC 1)$ in \cite[Definition 0.7]{GHJ}, and conditions $(\PC 2)$, $(\PC 2')$, $(\PC 3)$, $(\PC 3')$ correspond to $(\PC 2)$ in \cite[Definition 0.7]{GHJ}. A subcategory of $\C_m$ is precovering if and only if it satisfies the $\PC$ conditions, see \cite[Theorem 3.1]{GHJ}. This characterization generalises \cite{N} for the case $m = 1$.

\subsection{Precovering subcategories of $\overline{\mathcal{C}}_m$}\label{section precovering subcategories of the completion}

Now we classify the precovering subcategories of $\ovl\C_m$. Our approach is to relate the precovering subcategories of $\ovl\C_m$ to some subcategories of $\C_{2m}$ which are ``almost precovering". To do so, we need to introduce an auxiliary subcategory of $\C_{2m}$. Fix $z^0\in\Z$. For each $p\in[m']$ we denote by $z^0_p\in\Zc_{2m}$ the copy of $z^0$ belonging to $\Z^{(p)}$. 

\begin{definition}\label{definition A}
We define the subcategory $\A$ of $\C_{2m}$ as
\[
\A = \add\left\{(a_1,a_2)\in\ind\C_{2m}\middle| a_1,a_2\in\bigcup_{p\in[m]}(p,z^0_{p^+}]\right\}.
\]
\end{definition}

\begin{figure}[ht]
\centering
\includegraphics[height = 4cm]{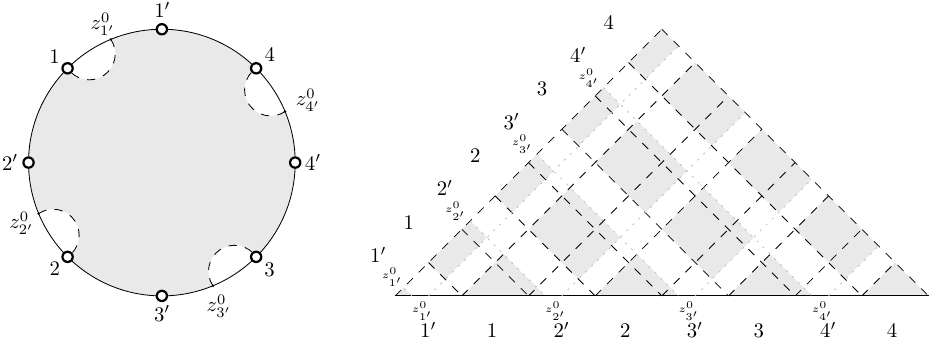}
\caption{The category $\A$.}
\label{figure the category A}
\end{figure}

Figure \ref{figure the category A} illustrates the subcategory $\A$.

Now we define the completed versions of the PC conditions.

\begin{figure}[ht]
\centering
\includegraphics[height = 4cm]{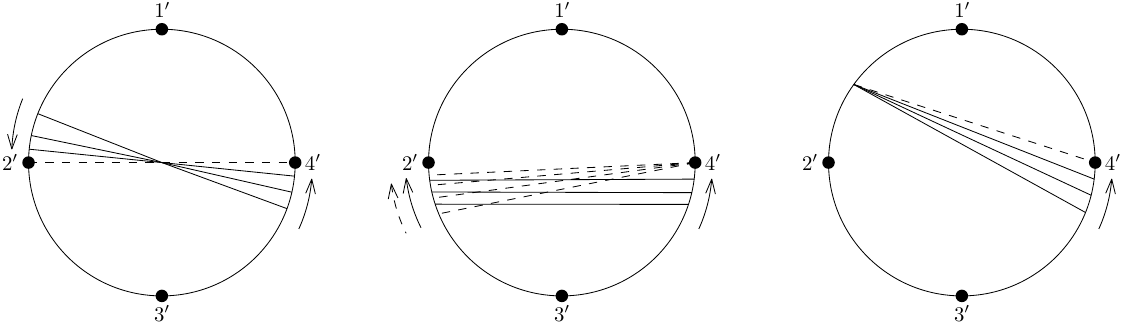}
\caption{On the left $(\ovl\PC 1)$, in the middle $(\ovl\PC 2)$, and on the right $(\ovl\PC 3)$.}
\label{figure pc conditions completion}
\end{figure}

\begin{definition}\label{definition completed pc conditions}
Let $\X$ be a subcategory of $\ovl\C_m$. We say that $\X$ satisfies the \emph{completed precovering conditions} (the $\ovl\PC$ conditions for short) if it satisfies the following combinatorial conditions.
\begin{enumerate}
\item[($\ovl\PC 1$)] If there exists a sequence $\{(x_1^n,x_2^n)\}_n\subseteq \X\cap C^{(p,q)}$ for some $p,q\in[m]$ such that $p\neq q$ and the sequences $\{x_1^n\}_n$ and $\{x_2^n\}_n$ are strictly increasing, then $|p^+,q^+|\in\X$.
\item[($\ovl\PC 2$)] If there exists a sequence $\{(x_1^n,x_2^n)\}_n\subseteq \X\cap C^{(p,q)}$ for some $p,q\in[m]$ and the sequences $\{x_1^n\}_n$  and $\{x_2^n\}_n$ are respectively strictly decreasing and strictly increasing, then there exists a strictly decreasing sequence $\{y_1^n\}_n\subseteq C^{(p)}$ such that $\{|y_1^n,q^+|\}_n\subseteq\X$.
\item[($\ovl\PC 2'$)]If there exists a sequence $\{(x_1^n,x_2^n)\}_n\subseteq \X\cap C^{(p,q)}$ for some $p,q\in[m]$ such that $p\neq q$, and the sequences $\{x_1^n\}_n$ and $\{x_2^n\}_n$ are respectively strictly increasing and strictly decreasing, then there exists a strictly decreasing sequence $\{y_2^n\}_n\subseteq C^{(q)}$ such that $\{(p^+,y_2^n)\}_n\subseteq \X$.
\item[($\ovl\PC 3$)] If there exists a sequence $\{(x_1,x_2^n)\}_n\subseteq \X\cap C^{(p,q)}$ for some $p\in[m']\cup[m]$ and $q\in[m]$ such that $p\neq q^+$ and the sequence $\{x_2^n\}_n$ is strictly increasing, then $|x_1,q^+|\in\X$.
\item[($\ovl\PC 3'$)] If there exists a sequence $\{(x_1^n,x_2)\}_n\subseteq\X\cap C^{(p,q)}$ for some $p\in[m]$ and $q\in[m']\cup[m]$ such that $p\neq q$, $q\neq p^+$, and the sequence $\{x_1^n\}_n$ is strictly increasing, then $(p^+,x_2)\in\X$.
\end{enumerate} 
\end{definition}
Figure \ref{figure pc conditions completion} illustrates some of the $\ovl\PC$ conditions.

The main result of this section is the following.

\begin{theorem}\label{theorem classification precovering subcategories}
Let $\X$ be a subcategory of $\ovl\C_m$. The following statements are equivalent.
\begin{enumerate}
\item $\X$ is precovering in $\ovl\C_m$.
\item $\pi^{-1}\X\cap\A$ is precovering in $\C_{2m}$.
\item $\X$ satisfies the $\ovl\PC$ conditions.
\end{enumerate}
\end{theorem}

The following lemmas will be useful to prove the theorem above.

\begin{lemma}\label{lemma A aisle t-structure}
The following statements hold.
\begin{enumerate}
\item The category $\A$ is the aisle of a t-structure in $\C_{2m}$. 
\item If $x\in\ind\C_{2m}$ and $x\notin\D$, then there exists an $\A$-cover $a\to x$ of $x$ such that $\pi a\cong \pi x$ and $a$ is indecomposable.
\end{enumerate}
\end{lemma}
\begin{proof}
We prove statement (1). From Proposition \ref{proposition iyama-yoshino torsion pairs} it is enough to check that $\A$ is suspended and precovering. It is straightforward to check that $\A$ satisfies the $\PC$ conditions, therefore by \cite[Theorem 3.1]{GHJ} it is precovering. For showing that $\S\A\subseteq\A$, consider $a = (a_1,a_2)\in\ind\A$, then, by Definition \ref{definition A}, $a_1-1,a_2-1\in\bigcup_{p\in[m]}(p,z_{p^+}^0-1]\subseteq\bigcup_{p\in[m]}(p,z_{p^+}^0]$. As a consequence, $\S a = (a_1-1,a_2-1)\in\ind\A$. Now we prove that $\A$ is extension-closed. Let $a\lora e\lora b\lora \S a$ be a triangle of $\C_{2m}$ with $a,b\in\A$. Then all the endpoints of the indecomposable summands of $a$ and $b$ belong to the set $\bigcup_{p\in[m]}(p,z_{p^+}^0]$. Therefore, by \cite[Lemma 3.4]{GZ}, all the endpoints of the indecomposable summands of $e$ belong to the same set, and it follows that $e\in\A$.

Now we prove statement (2). Let $x = (x_1,x_2)\in\ind\C_{2m}$, by \cite[Lemma 4.1]{J} $x$ has an $\A$-cover $a\to x$. We show that $a\in\ind\A$ and that $\pi a\cong \pi x$. Let $p,q\in[m']\cup[m]$ be such that $x\in \Z^{(p,q)}$, and assume that $p \in[m']$ and $q\in[m]$, the other cases are analogous. If $x_1\in (p,z_p^0]$ then $x\in\ind\A$ and $1_x\colon x\to x$ is an $\A$-cover of $x$. If $x_1\notin (p,z_p^0]$ then let $a' = (z_p^0,x_2)\in\ind\C_{2m}$. We have that $a'\in\ind\A$. There exists a non-zero morphism $a'\to x$ and it is straightforward to check that it is an $\A$-precover. Moreover $a'\to x$ is right-mimimal and therefore an $\A$-cover. Since covers are unique up to isomorphism, we have that $a\cong a'$. We also have that $\pi a \cong \pi a' = (\ovl{z}_p^0,\ovl{x}_2) = (\ovl{x}_1,\ovl{x}_2) = \pi x$. This concludes the argument.
\end{proof}

Consider $a,b\in\ind\ovl\C_m$ such that $\Hom_{\ovl\C_m}(b,a) \cong \K$. We recall, from Lemma \ref{lemma reverse hammocks}, that we can fix $a'\in\ind\C_{2m}$ such that $\pi(a') \cong a$, and then there exists $b'\in\ind\C_{2m}$, which depends on the choice of $a'$, such that $\pi(b')\cong b$, $\Hom_{\C_{2m}}(b',a')\cong\K$, and any non-zero morphism $b'\to a'$ does not factor through $\D$.

\begin{lemma}\label{lemma hom sets A-cover}
Let $a,b\in\ind\ovl\C_m$ be such that $\Hom_{\ovl\C_{m}}(b,a) \cong \K$, and let $a'\in\ind\C_{2m}$ be such that $\pi(a')\cong a$. Let $b'\in\ind\C_{2m}$ be such that $\pi(b')\cong b$, $\Hom_{\C_{2m}}(b',a')\cong\K$, and any non-zero morphism $b'\to a'$ does not factor through $\D$. Let $b''\to b'$ be the $\A$-cover of $b'$. Then $\Hom_{\C_{2m}}(b'',a')\cong\K$ and any non-zero morphism $b''\to a'$ does not factor through $\D$.
\end{lemma}
\begin{proof}
By Lemma \ref{lemma reverse hammocks} such $b'$ exists. If $b'\in\A$ then $b''\cong b'$ and we have the statement. Assume that $b'\notin\A$ and let $p,q\in[m']\cup[m]$ be such that $b' = (b_1',b_2')\in\Z^{(p,q)}$. Consider the case $p\in[m']$ and $q\in[m]$, the other cases are analogous. Since $b'\notin\A$ we have that $b_1'\notin (p,z_p^0]$, i.e. $b_1'\in[z_p^0+1,p^+)$. From the argument of Lemma \ref{lemma A aisle t-structure}, $b'' = (z_p^0,b_2')$. We have that $a'\in H^+(b')\cup H^-(\S^2 b')$, we show that $a'\in H^+(b'')\cup H^-(\S^2 b'')$. 

If $a'\in H^+(b')$ then $b_1'\leq a_1'\leq b_2'-2$ and $a_2'\geq b_2'$. Since $b_1'>z_p^0$, then $a'\in H^+(b'')$. Now, if $a'\in H^-(\S^2 b')$ then $a_1'\leq b_1'-2$ and $b_1'\leq a_2'\leq b_2'-2$. Assume that $a_1'\not\leq z_p^0-2$, then $b_1'-2\leq a_1'\leq z_p^0-1$. In particular, $a_1'\in \Z^{(p)}$ and any non-zero morphism $b'\to a'$ factors through $\D$ giving a contradiction, see Figure \ref{figure hammocks}. Therefore $a_1'\leq z_p^0-2$. Moreover, since $b_1'>z_p^0$, from $b_1'\leq a_2'\leq b_2'-2$ we also have that $z_p^0\leq a_2'\leq b_2'-2$ and obtain that $a'\in H^-(\S^2 b'')$. We can conclude that $\Hom_{\C_{2m}}(b'',a')\cong\K$. 

We show that any non-zero morphism $b''\to a'$ does not factor through $\D$. If this is not the case, then $a'\in H^-(\S^2 b'')$ and $a_1'\in\Z^{(p)}$. As a consequence, $\Hom_{\ovl\C_m}(b,a) = 0$ giving a contradiction. We obtain that any non-zero morphism $b''\to a'$ does not factor through $\D$, and this concludes the argument.
\end{proof}

\begin{figure}[ht]	
\centering
\includegraphics[height = 8cm]{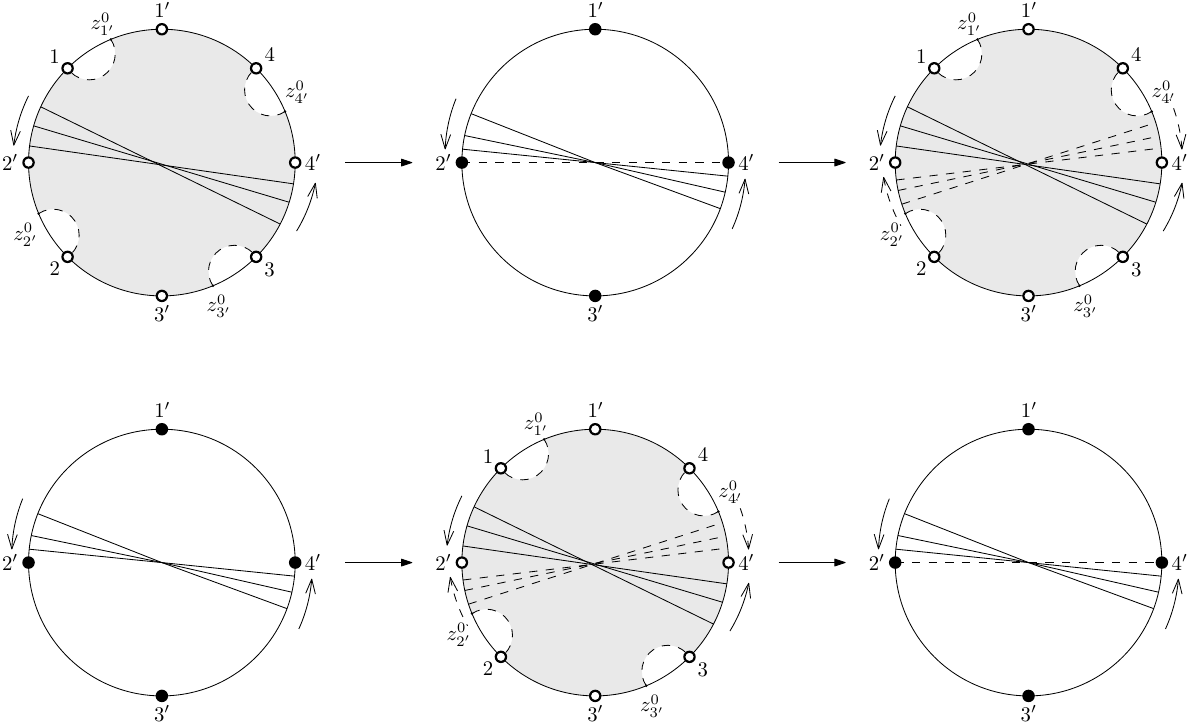}
\caption{Illustration of the argument of Lemma \ref{lemma X completed PC conditions if and only if preimage of X intersected with A PC conditions}.}
\label{figure pc conditions preimage pc conditions}
\end{figure}

\begin{lemma}\label{lemma X completed PC conditions if and only if preimage of X intersected with A PC conditions}
Let $\X$ be a subcategory of $\ovl\C_m$. The subcategory $\X$ satisfies the $\ovl\PC$ conditions if and only if $\pi^{-1}\X\cap \A$ satisfies the $\PC$ conditions.
\end{lemma}
\begin{proof}	
We show that $\X$ satisfies $(\ovl\PC 1)$ if and only if $\pi^{-1}\X\cap\A$ satisfies $(\PC 1)$, we refer to Figure \ref{figure pc conditions preimage pc conditions} for an illustration. Assume that $\X$ satisfies $(\ovl\PC 1)$ and that there exists a sequence $\{(x_1^n,x_2^n)\}_n\subseteq \Z^{(p,q)}\cap \left(\pi^{-1}\X\cap \A\right)$ for some $p,q\in[m']\cup[m]$ such that $p\neq q$ with $\{x_1^n\}_n$ and $\{x_2^n\}_n$ strictly increasing sequences. Note that $p,q\notin[m']$, otherwise for $n$ big enough we have $\{(x_1^n,x_2^n)\}_n\not\subseteq\A$. For each $n$ we define $\pi x^n = y^n = (y_1^n,y_2^n)$, note that $\{y_1^n\}_n$ and $\{y_2^n\}_n$ are still strictly increasing sequences, and consider the sequence $\{y^n\}_n\subseteq \X\cap C^{(p,q)}$. Since $\X$ satisfies $(\ovl\PC 1)$, then $|p^+,q^+|\in\X$ and $\pi^{-1}\X$ contains any arc of $\C_{2m}$ having one endpoint in $\Z^{(p^+)}$ and the other in $\Z^{(q^+)}$. In particular, there exist strictly decreasing sequences $\{z_1^n\}_n\subseteq \Z^{(p^+)}$ and $\{z_2^n\}_n\subseteq \Z^{(q^+)}$ such that $\{|z_1^n,z_2^n|\}_n\subseteq \pi^{-1}\X\cap \A$. This proves that $\pi^{-1}\X\cap \A$ satisfies $(\PC 1)$. 
	
Now assume that $\pi^{-1}\X\cap \A$ satisfies $(\PC 1)$ and that there exists a sequence $\{(x_1^n,x_2^n)\}_n\subseteq \X\cap C^{(p,q)}$ for some $p,q\in[m]$ such that $p\neq q$, and $\{x_1^n\}_n$ and $\{x_2^n\}_n$ are strictly increasing sequences. For each $n$ there exists $y^n = (y_1^n,y_2^n)\in\ind\pi^{-1}\X\cap \A$ such that $\pi y^n\cong x^n$. Thus, there exists a sequence $\{|y_1^n,y_2^n|\}_n\subseteq \left(\pi^{-1}\X\cap \A\right)\cap C^{(p,q)}$ such that $\{y_1^n\}_n$ and $\{y_2^n\}_n$ are stricly increasing sequences. Since $\pi^{-1}\X\cap \A$ satisfies $(\PC 1)$, then there exist strictly decreasing sequences $\{z_1^n\}_n\subseteq\Z^{(p^+)}$ and $\{z_2^n\}_n\subseteq\Z^{(q^+)}$ such that $\{|z_1^n,z_2^n|\}_n \subseteq \pi^{-1}\X\cap \A$. As a consequence we have that $|p^+,q^+|\in\X$. This proves that $\X$ satisfies $(\ovl\PC 1)$.

It is straightforward to check that $\X$ satisfies $(\ovl\PC 3)$ and $(\ovl\PC3')$ if and only if $\pi^{-1}\X\cap\A$ satisfies $(\PC 3)$ and $(\PC 3')$. Moreover, if $\X$ satisfies $(\ovl\PC 2)$, $(\ovl\PC 2')$, $(\ovl\PC 3)$, and $(\ovl\PC 3')$ then $\pi^{-1}\X\cap\A$ satisfies $(\PC 2)$ and $(\PC 2')$. Finally, if $\pi^{-1}\X\cap\A$ satisfies $(\PC 2)$ and $(\PC 2')$ then $\X$ satisfies $(\ovl\PC 2)$ and $(\ovl\PC 2')$. We can conclude that $\X$ satisfies the $\ovl\PC$ conditions if and only if $\pi^{-1}\X\cap\A$ satisfies the $\PC$ conditions.
\end{proof} 

\begin{proposition}\label{proposition if preimage of Y almost precovering then Y precovering}
Let $\X$ be a subcategory of $\ovl\C_m$. If $\pi^{-1}\X\cap \A$ is a precovering subcategory of $\C_{2m}$ then $\X$ is a precovering subcategory of $\ovl\C_m$.
\end{proposition}
\begin{proof}
From Remark \ref{remark precovering at the level of indecomposables} it is enough to check that $\X$ is precovering at the level of the indecomposable objects. Consider $a\in\ind\ovl\C_m$, then there exists $a'\in\ind\C_{2m}$ such that $\pi a'\cong a$, and there exists $f\colon x\to a'$ a $\pi^{-1}\X\cap \A$-precover of $a'$. Consider $\pi f\colon \pi x\to a$, we show that $\pi f$ satisfies the condition of Remark \ref{remark precovering at the level of indecomposables}. First assume that $f$ does not factor through $\D$. Consider $b\in\ind\X$ and $g\colon b\to a$ in $\ovl\C_m$. Without loss of generality we can assume that $g\neq 0$. From Lemma \ref{lemma A aisle t-structure} and Lemma \ref{lemma hom sets A-cover}, there exists $b'\in\ind \pi^{-1}\X\cap \A$ such that $\pi b'\cong b$ and there exists a non-zero morphism $g'\colon b'\to a'$ in $\C_{2m}$ which does not factor through $\D$. Since the $\Hom$-spaces in $\ovl\C_m$ are at most one dimensional, we have that $g = \l\pi g'$ for some $\l\in\K^*$. Since $f\colon x\to a'$ is a $\pi^{-1}\X\cap \A$-precover of $a'$, there exists $h\colon b'\to x$ in $\C_{2m}$ such that $fh = g'$. We obtain that $\l\pi (f)\pi (h) = \l\pi(fh) = g$ in $\ovl\C_m$. This proves that $\pi f\colon \pi x\to a$ is an $\X$-precover of $a$.

Now we consider the case where $f$ factors through $\D$. We show that $\Hom_{\ovl\C_m}(b,a) = 0$ for all $b\in\ind\X$. Assume that there exists a non-zero morphism $g\colon b\to a$ in $\ovl\C_m$ for some $b\in\ind\X$, then as above there exists $b'\in\ind \pi^{-1}\X\cap \A$ such that $\pi b'\cong b$ and there exists a non-zero morphism $g'\colon b'\to a'$ in $\C_{2m}$ which does not factor through $\D$. Since $f\colon x\to a'$ is a $\pi^{-1}\X\cap\A$-precover of $a'$, there exists $h\colon b'\to x$ in $\C_{2m}$ such that $fh = g'$. Since $f$ factors through $\D$, we have that $g'$ factors through $\D$, giving a contradiction. We can conclude that if $f$ factors through $\D$ then $\Hom_{\ovl\C_m}(b,a) = 0$ for all $b\in\ind\X$. As a consequence, $\pi f = 0$ is an $\X$-precover of $a$.
\end{proof}

The following proposition is the analogue of \cite[Proposition 3.7]{GHJ} in $\ovl\C_m$ and its proof is similar. 

\begin{proposition}\label{proposition if precovering then completed PC conditions}
Let $\X$ be a subcategory of $\ovl\C_m$. If $\X$ is a precovering subcategory then it satisfies the $\ovl\PC$ conditions.
\end{proposition}
\begin{proof}
Assume that $\X$ is a precovering subcategory of $\ovl\C_m$, we show that is satisfies $(\ovl\PC 1)$. Assume that there is a sequence $\{x^n = (x_1^n,x_2^n)\}_n\subseteq\ind\X\cap C^{(p,q)}$ for some $p,q\in[m]$ with $p\neq q$ such that $\{x_1^n\}_n$ and $\{x_2^n\}_n$ are strictly increasing sequences. We show that $a = |p^+,q^+|\in\ind\X$. Consider $(f_1,\dots,f_k)\colon y_1\oplus\cdots\oplus y_k\to a$ an $\X$-precover of $a$ with $y_1,\dots,y_k\in\ind\X$. Note that $\Hom_{\ovl\C_m}(x^n, a)\cong \K$ for each $n$. Fix $n$ and consider a non-zero morphism $g^n\colon x^n\to a$. Then there exists $h^n = (h_1^n,\dots,h_k^n)^T\colon x^n\to y_1\oplus\cdots\oplus y_k$ such that $fh^n = g^n$. Then, for each $n$ there exists $l\in\{1,\dots,k\}$ such that $g^n$ factors through $f_l$. 
	
There exists $l\in\{1,\dots,k\}$ such that $g^n$ factors through $f_l$ for infinitely many $n\in\Z$. Indeed, if for each $l$ only finitely many of the $g^n$ factor through $f_l$, then there are only finitely many $g^n$'s and this contradicts the fact that the sequence $\{x^n\}_n$ is infinite. Now fix an $l$ such that $g^n$ factors through $f_l$ for infinitely many $n\in\Z$. Without loss of generality we can assume that for each $n\in\Z$ the morphism $g^n$ factors through $f_l$. Indeed, if this is not the case, we can extract an infinite subsequence of $\{x^n\}_n$ such that all $g^n\colon x^n\to a$ satisfy that property. From now on we denote the object $y_l$ as $y$, the morphism $f_l\colon y_l\to a$ as $f\colon y\to a$, and we denote by $h^n\colon x^n\to y$ the morphism such that $f h^n = g^n$.

Since $\Hom_{\ovl\C_m}(y,a)\cong\K$ and $\Hom_{\ovl\C_m}(x^n,y)\cong \K$ for all $n\in\Z$, we have that 
\[
y\in\left(\bigcap_{n\in\Z} \ovl{H}^+(x^n)\cup \ovl{H}^-(\S^2 x^n)\right) \cap \left( \ovl I^-(a)\cup \ovl I^+(\S^{-2} a)\right).
\]
We refer to Figure \ref{figure precovering then pc conditions} for an illustration. 

\begin{figure}[ht]
\centering
\includegraphics[height = 3.9cm]{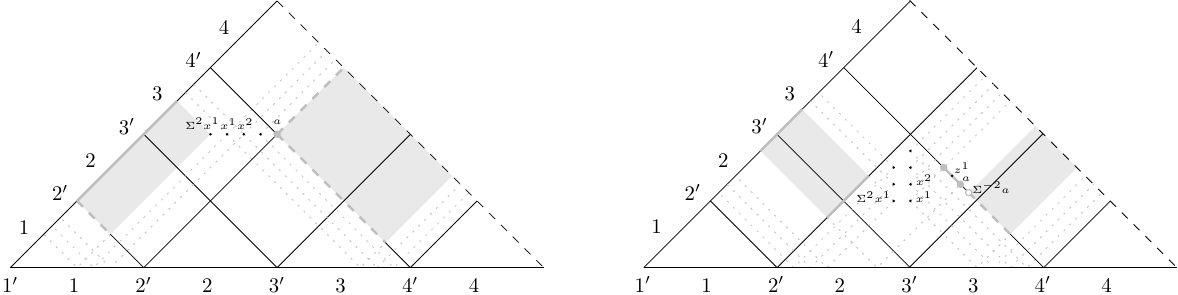}
\caption{On the left the argument for $(\ovl\PC 1)$, on the right for $(\ovl\PC 2)$. The object $y$ belongs to the grey areas.}
\label{figure precovering then pc conditions}
\end{figure}

We show that $y\cong a$. Assume that $y\not\cong a$, from Proposition \ref{proposition factorization properties} there exists a non-zero morphism $f'\colon a\to y$ such that $h^n = f'g^n$ for each $n\in\Z$. Since $fh^n = g^n\neq 0$, then $ff'g^n \neq 0$ and $ff'\colon a\to a$ is non-zero. Thus, $ff' = \l 1_a$ for some $\l\in\K^*$ and $a\cong y$, which contradicts our assumption. We obtain that $a\cong y\in\X$, and we conclude that $\X$ satisfies $(\ovl\PC 1)$.
	
Now we show that $\X$ satisfies $(\ovl\PC 2)$. Assume that there is a sequence $\{(x_1^n,x_2^n)\}\subseteq\ind\X\cap C^{(p,q)}$ for some $p,q\in[m]$ such that $\{x_1^n\}_n$ is strictly decreasing and $\{x_2^n\}_n$ is strictly increasing. We show that there is a strictly decreasing sequence $\{y_1^n\}_n\subseteq C^{(p)}$ such that $\{|y_1^n,q^+|\}_n\subseteq\X$.
	
Consider an object $a = |a_1,q^+|$ with $a_1\in\Z^{(p)}$ such that $x_1^1<a_1\leq x_2^1-2$. Then for each $n$ there exists a non-zero morphism $g^n\colon x^n\to a$. Consider an $\X$-precover $(f_1,\dots,f_k)\colon y_1\oplus\dots\oplus y_k\to a$ of $a$. With the same argument as above there exists $l\in\{1,\dots,k\}$ such that $g^n\colon x^n\to a$ factors through $f_l\colon y_l\to a$ for all $n$ (up to taking subsequences). Let $y = y_l$, proceeding similarly as above we obtain that $y\in\left\{|z,q^+|\middle| x_1^1\leq z\leq a_1^1\right\}$, see Figure \ref{figure precovering then pc conditions}. We define $z^1 = y$, which is the first element of our desired sequence. Now we consider $a' = |a_1',q^+|$ with $x_1^1<a_1'\leq x_1^2$. By repeating the same argument there exists $z^2\in\left\{|z,q^+|\middle| a_1'\leq z\leq a_1\right\}$ which is an object of $\X$. With this procedure we obtain our desired sequence $\{z^n\}_n$. This proves that $\X$ satisfies $(\ovl\PC 2)$. 

The argument of $(\ovl\PC 2')$ is similar to the argument of $(\ovl\PC 2)$, the arguments of $(\ovl\PC 3)$ and $(\ovl\PC 3')$ are similar to the argument of $(\ovl\PC 1)$. We can conclude that $\X$ satisfies the $\ovl\PC$ conditions.
\end{proof}

We now have our classification of the precovering subcategories of $\ovl\C_m$.

\begin{proof}[Proof of Theorem \ref{theorem classification precovering subcategories}]
The claim follows directly from Lemma \ref{lemma X completed PC conditions if and only if preimage of X intersected with A PC conditions}, Proposition \ref{proposition if preimage of Y almost precovering then Y precovering}, Proposition \ref{proposition if precovering then completed PC conditions}, and \cite[Theorem 3.1]{GHJ}.	
\end{proof}

\subsection{Preenveloping subcategories of $\overline{\mathcal{C}}_m$}\label{section preenveloping subcategories}

Here we discuss a characterization of the preenveloping subcategories of $\ovl\C_m$ dual to Theorem \ref{theorem classification precovering subcategories}. First we define an auxiliary subcategory $\B$, which is the dual version of $\A$. We recall that in Section \ref{section precovering subcategories of the completion} we fixed an integer $z^0\in\Z$, now we define $w^0 = z^0-1$. For each $p\in[m']$ we denote by $w_p^0\in\Zc_{2m}$ the copy of $w^0$ belonging to $\Z^{(p)}$.

\begin{definition}\label{definition B}
We define the subcategory $\B$ of $\C_{2m}$ as
\[
\B = \add\left\{(b_1,b_2)\in\ind\C_{2m}\middle| b_1,b_2\in\bigcup_{p\in[m']}[w_p^0,p^{++})\right\}.
\]
\end{definition}

\begin{figure}[ht]
\centering
\includegraphics[height = 4cm]{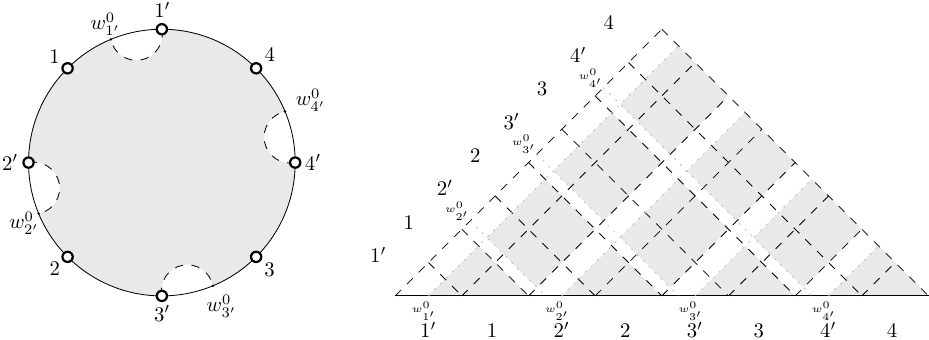}
\caption{The category $\B$.}
\label{figure the category B}
\end{figure}

Figure \ref{figure the category B} illustrates the subcategory $\B$.

For convenience of the reader, we record the duals of Definition \ref{definition completed pc conditions} and Theorem \ref{theorem classification precovering subcategories}.

\begin{definition}\label{definition completed pe conditions}
Let $\X$ be a subcategory of $\ovl\C_m$. We say that $\X$ satisfies the \emph{completed preenveloping conditions} ($\ovl\PE$ for short) if it satisfies the following combinatorial conditions.
\begin{enumerate}
\item[($\ovl\PE 1$)] If there exists a sequence $\{(x_1^n,x_2^n)\}_n\subseteq \X\cap C^{(p,q)}$ for some $p,q\in[m]$ such that $p\neq q$ and the sequences $\{x_1^n\}_n$ and $\{x_2^n\}_n$ are strictly decreasing, then $(p^-,q^-)\in\X$.
\item[($\ovl\PE 2$)] If there exists a sequence $\{(x_1^n,x_2^n)\}_n\subseteq \X\cap C^{(p,q)}$ for some $p,q\in[m]$ such that $p\neq q$ and the sequences $\{x_1^n\}_n$  and $\{x_2^n\}_n$ are respectively strictly increasing and strictly decreasing, then there exists a strictly increasing sequence $\{y_1^n\}_n\subseteq C^{(p)}$ such that $\{(y_1^n,q^-)\}_n\subseteq\X$.
\item[($\ovl\PE 2'$)]If there exists a sequence $\{(x_1^n,x_2^n)\}_n\subseteq \X\cap C^{(p,q)}$ for some $p,q\in[m]$ such that the sequences $\{x_1^n\}_n$ and $\{x_2^n\}_n$ are respectively strictly decreasing and strictly increasing, then there exists a strictly increasing sequence $\{y_2^n\}_n\subseteq C^{(q)}$ such that $\{(p^-,y_2^n)\}_n\subseteq \X$.
\item[($\ovl\PE 3$)] If there exists a sequence $\{(x_1,x_2^n)\}_n\subseteq \X\cap C^{(p,q)}$ for some $p\in[m']\cup[m]$ and $q\in[m]$ such that $p\neq q$, $p\neq q^-$ and the sequence $\{x_2^n\}_n$ is strictly decreasing, then $(x_1,q^-)\in\X$.
\item[($\ovl\PE 3'$)] If there exists a sequence $\{(x_1^n,x_2)\}_n\subseteq\X\cap C^{(p,q)}$ for some $p\in[m]$ and $q\in[m']\cup[m]$ such that the sequence $\{x_1^n\}_n$ is strictly decreasing, then $(p^-,x_2)\in\X$.
\end{enumerate}
\end{definition}

\begin{theorem}\label{theorem classification preenveloping subcategories}
Let $\X$ be a subcategory of $\ovl\C_m$. The following statements are equivalent.
\begin{enumerate}
\item $\X$ is preenveloping in $\ovl\C_m$.
\item $\pi^{-1}\X\cap\B$ is preenveloping in $\C_{2m}$.
\item $\X$ satisfies the $\ovl\PE$ conditions.
\end{enumerate}
\end{theorem}

The following lemma will be useful in Section \ref{section co-aisles t-structures} for computing the heart of a t-structure.

\begin{lemma}\label{lemma B co-aisle t-structure}
The following statements hold.
\begin{enumerate}
\item The category $\B$ is the co-aisle of a t-structure in $\C_{2m}$.
\item If $x\in\ind\C_{2m}$ and $x\notin\D$, then there exists $x'\in\ind\A\cap\S^{-1}\B\subseteq\ind\A\cap\S\B$ such that $\pi x'\cong x$.
\end{enumerate}
\end{lemma}
\begin{proof}
Statement (1) is the dual of statement (1) of Lemma \ref{lemma A aisle t-structure}. For statement (2), consider the $\A$-cover $x'\to x$ of $x$ as in statement (2) of Lemma \ref{lemma A aisle t-structure}. We have that $\pi(x')\cong x$, and it is straightforward to check that $x'\in\S^{-1}\B$. Moreover, since $\S^{-1}\B\subseteq\B\subseteq\S\B$, we have that $x'\in\S\B$. This concludes the proof.
\end{proof}

\section{Extension-closed subcategories}\label{section extension-closed subcategories}

In this section we classify the extension-closed subcategories of $\ovl\C_m$. To do so, we first show that the extension-closed subcategories of $\C_m$ are precisely those closed under extensions having indecomposable outer terms.

\subsection{Extension-closed subcategories of $\mathcal{C}_m$}

The extension-closed subcategories of $\C_m$ were classified in \cite[Theorem 7.2]{CP} for the case $m = 1$. The precovering extension-closed subcategories of $\C_m$, i.e. the torsion classes, were classified in \cite[Theorem 4.7]{GHJ}. Here we classify the subcategories of $\C_m$ which are just extension-closed for all $m\geq 1$.

We recall that we identify the indecomposable objects of $\C_m$ with the arcs of $\Zc_m$.

\begin{definition}
Let $a,b\in\ind\C_m$ be crossing arcs. The arcs of $\ind\C_m\setminus\{a,b\}$ which connect the endpoints of $a$ and $b$ are called \emph{Ptolemy arcs}. We say that a subcategory $\U$ of $\C_m$ satisfies the \emph{Ptolemy condition}, $\PT$ condition for short, if it is closed under taking Ptolemy arcs.
\end{definition}

\begin{figure}[ht]
\centering
\includegraphics[height = 4cm]{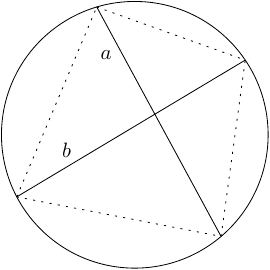}
\caption{The dotted arcs are the Ptolemy arcs of $a$ and $b$.}
\label{figure Ptolemy arcs}
\end{figure}

Figure \ref{figure Ptolemy arcs} provides an illustration of Ptolemy arcs.

Consider a non-split triangle of the form $a\lora e\lora b\lora \S a$ with $a,b\in\ind\C_m$ and $b\not\cong \S a$. The middle term $e$ is determined by the Ptolemy arcs of $a$ and $b$. More precisely 
\begin{itemize}
\item if $a_1<b_1<a_2<b_2$ then $e \cong e_1\oplus e_2$ with $e_1 = (a_1,b_2)$ and $e_2 = (b_1,a_2)$, and
\item if $b_1<a_1<b_2<a_2$ then $e\cong e_1'\oplus e_2'$ with $e_1' = (b_1,a_1)$ and $e_2' = (b_2,a_2)$.
\end{itemize}
In the first case, if $a_2 = b_2+1$ we interpret $(b_1,a_2)$ as the zero object. In the second case, if $a_1 = b_1+1$ we interpret $(b_1,a_1)$ as the zero object, and if $a_2 = b_2+1$ we interpret $(b_2,a_2)$ as the zero object. 

Now consider a triangle of the form $a\lora e\lora b\overset{h}{\lora} \S a$ with $a,b_1,\dots,b_n\in\ind\C_m$, such that the objects $b_1,\dots,b_n$ are non-isomorphic to $\S a$ and pairwise $\Hom$-orthogonal, i.e. $\Hom_{\C_m}(b_i,b_j) = 0$ for each $i\neq j$, and all the entries of $h = (h_1,\dots,h_n)$ are non-zero, cf. \cite[Lemma 3.2]{GZ}. The middle term $e$ of such a triangle was computed in \cite[Lemma 4.16]{GZ}. Their result generalises \cite[Proposition 4.12]{CP} for the case $m = 1$. Now we show that computing the middle term of a triangle of that form is enough to obtain the middle term of a triangle of the form $a\lora e\lora b\lora \S a$ with $a\in\ind\C_m$. This follows from Proposition \ref{proposition decomposition of a triangle in triangulated setting}, which holds in a more general framework.

\begin{remark}\label{remark factorization free hom orthogonal}
Let $f_1\colon x_1\to y$ and $f_2\colon x_2\to y$ be non-zero morphisms in $\C_m$ with $x_1,x_2,y\in\ind\C_m$. By applying the same argument of \cite[Lemma 4.6]{CP}, we obtain that $f_1$ and $f_2$ are factorization free if and only if $x_1$ and $x_2$ are $\Hom$-orthogonal. We observe that the ``if" implication holds in the setting of additive categories, while the ``only if" implication is specific to the category $\C_m$.
\end{remark}

\begin{corollary}[of Proposition \ref{proposition decomposition of a triangle in triangulated setting}]\label{corollary decomposition of a triangle}
Let $a\lora e \lora b\overset{h}{\lora} \S a$ be a triangle in $\C_m$ with $a,b_1,\dots,b_n\in\ind\C_m$, $b = \bigoplus_{i = 1}^n b_i$, and $h = (h_1,\dots,h_n)$. Then there exist $b_1',\dots,b_k'\in\ind\C_m$ and a morphism $h' = (h_1'\dots,h_k')\colon b' = \bigoplus_{i = 1}^k b_i'\to \S a$ such that $b'$ is a direct summand of $b$, $b_1',\dots,b_k'$ are non-isomorphic to $\S a$ and are pairwise $\Hom$-orthogonal, all the entries of $h'$ are non-zero, and there is the following isomorphism of triangles.
\[
\begin{tikzcd}
a\ar[r]\ar[d,"1"] & e'\oplus b''\ar[r]\ar[d,"\wr"] & b'\oplus b''\ar[r, "{(h', 0)}"]\ar[d,"\wr"]  & \S a\ar[d,"1"]\\
a\ar[r] & e\ar[r]& b\ar[r,"h"]& \S a
\end{tikzcd}
\]
\end{corollary}
\begin{proof} 
By Proposition \ref{proposition decomposition of a triangle in triangulated setting} there exists $b_1',\dots,b_k'\in\ind\C_m$ such that $b' = \bigoplus_{i = 1}^k b_i'$ is a direct summand of $b$, all the entries of $h' = (h_1',\dots,h_k')\colon b'\to \S a$ are pairwise factorization free, and there exists an isomorphism of triangles as in the statement. By Remark \ref{remark consequence lemma decomposition of a triangle in triangulated setting} and Remark \ref{remark factorization free hom orthogonal} we conclude that $b_i'\not\cong\S a$ and $h_i'\neq 0$ for each $i$, and $b_1',\dots,b_k'$ are pairwise $\Hom$-orthogonal.
\end{proof}

From \cite[Theorem 4.1]{CP} we know that the middle terms of arbitrary triangles of $\C_m$ can be computed iteratively when $m = 1$. It is straightforward to check that the same holds for $m\geq 2$. 

Let $\U$ be a subcategory of $\C_m$, and consider a triangle $a\lora e\lora b\lora \S a$ in $\C_m$ with $a,b\in\U$. From \cite[Lemma 3.4]{GZ} the coordinates of the indecomposable summands of $e$ belong to the set of coordinates of the indecomposable summands of $a$ and $b$. Note that in general this does not imply that $e\in\U$. Now we discuss a necessary and sufficient condition for $\U$ to be extension-closed. First, we need the following lemma.

\begin{lemma}\label{lemma subcategory closed under extensions with indecomposable outer terms implies closed under extensions}
Let $\U$ be a subcategory of $\C_m$. Assume that $\U$ is closed under extensions of the form $a\lora e\lora b\lora \S a$ with $a,b\in\ind\U$. Then $\U$ is closed under extensions of the form $a\lora e\lora b\overset{h}{\lora} \S a$ with $a\in\ind\C_m$, $b = \bigoplus_{i = 1}^n b_i$ where $b_1,\dots,b_n\in\ind\C_m$ are non-isomorphic to $\S a$ and pairwise $\Hom$-orthogonal, and all the entries of $h = (h_1,\dots,h_n)$ are non-zero.
\end{lemma}
\begin{proof}
Consider a triangle $a\lora e\lora b\overset{h}{\lora} \S a$ with $a\in\ind\C_m$, $b = \bigoplus_{i = 1}^n b_i$ where $b_1,\dots,b_n\in\ind\C_m$ are non-isomorphic to $\S a$ and pairwise $\Hom$-orthogonal, and all the entries of $h = (h_1,\dots,h_n)$ are non-zero. We prove that $e\in \U$. We proceed by induction on $n$. If $n = 1$ then $e\in\U$ by assumption. Now assume that $n\geq 2$. For each $i\in\{1,\dots,n\}$ we have that $b_i$ and $a$ cross, i.e. $b_i\in H^+(\S^{-1}a)\cup H^-(\S a)$. We have the following possibilities: $b_i\in H^-(\S a)$ for each $i$, or there exists $i$ such that $b_i\in H^+(\S^{-1} a)$. In the first case, we rename $b_1,\dots,b_n$ in such a way that the first coordinate of $b_n$ is the minimum of the first coordinates of $b_1,\dots,b_n$. In the second case, we rename $b_1,\dots,b_n$ in such a way that the first coordinate of $b_n$ is the maximum of the first coordinates of $b_1,\dots,b_n$. We refer to Figure \ref{figure a crossing bn} for an illustration.
\begin{figure}[ht]
\centering
\includegraphics[height = 4cm]{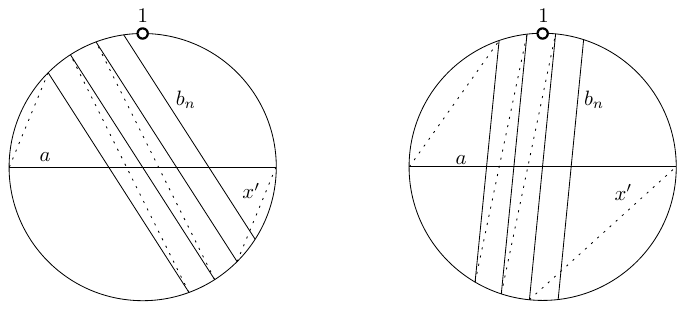}
\caption{On the left when $b_i\in H^-(\S a)$ for each $i\in\{1,\dots,n\}$, on the right when there exists $i\in\{1,\dots,n\}$ such that $b_i\in H^+(\S^{-1} a)$. The dotted arcs represent the indecomposable direct summands of $x$.}
\label{figure a crossing bn}
\end{figure}
We write $b = b'\oplus b_n$, where $b' = \bigoplus_{i = 1}^{n-1}b_i$. Consider the following Octahedral Axiom diagram.
\[
\begin{tikzcd}[ampersand replacement=\&]
\& b'\ar[r,"1"]\ar[d,"\left(\begin{smallmatrix}1\\0\end{smallmatrix}\right)"'] \& b'\ar[d, "{h'}", dashed] \& \\
e\ar[r]\ar[d,"1"'] \& b'\oplus b_n\ar[r,"h"]\ar[d,"{\left(\begin{smallmatrix} 0 & 1\end{smallmatrix}\right)}"'] \& \S a\ar[r]\ar[d,dashed] \& \S e\ar[d,"1"'] \\
e\ar[r] \& b_n\ar[r,"g"]\ar[d,"0"'] \& \S x\ar[r]\ar[d,dashed] \& \S e \\ \& \S b'\ar[r,"1"] \& \S b' \& 
\end{tikzcd}
\]
Consider the triangle $a\lora x\lora b'\overset{h'}{\lora}\S a$, and note that all the entries of $h'$ are non-zero. Thus, by induction hypothesis we have that $x\in\U$. By \cite[Lemma 4.16]{GZ} we can compute the indecomposable direct summands of $x$. Moreover, there is precisely one indecomposable direct summand of $x$, denoted by $x'$, such that $\Hom_{\C_m}(b_n,x')\neq 0$, see Figure \ref{figure a crossing bn}. We write $x = x'\oplus x''$ and $g = \left(\begin{smallmatrix}g'\\0\end{smallmatrix}\right)\colon b_n\to \S x'\oplus \S x''$. Since $x\in\U$, we have that $x',x''\in\U$. Consider the triangle $x'\lora e'\lora b_n\overset{g'}{\lora}\S x'$ and note that $e'\in\U$ because $x',b_n\in\ind\U$. We have the following isomorphism of triangles.
\[
\begin{tikzcd}
x'\oplus x''\ar[r]\ar[d,"1"'] & e\ar[r]\ar[d,"\wr"'] & b_n\ar[r, "\left(\begin{smallmatrix} g' \\ 0 \end{smallmatrix}\right)"]\ar[d,"1"'] & \S x'\oplus \S x''\ar[d,"1"']\\
x'\oplus x''\ar[r] & e'\oplus x''\ar[r] & b_n\ar[r,"\left(\begin{smallmatrix} g' \\ 0 \end{smallmatrix}\right)
"] & \S x'\oplus \S x''
\end{tikzcd}
\]
Therefore, we conclude that  $e\cong e'\oplus x''\in\U$.
\end{proof}

\begin{proposition}\label{proposition extension-closed subategories of IT category}
Let $\U$ be a subcategory of $\C_m$. Then the following statements are equivalent.
\begin{enumerate}
\item The subcategory $\U$ satisfies the $\PT$ condition.
\item The subcategory $\U$ is closed under extensions of the form $a\lora e\lora b\lora \S a$ where $a,b\in\ind\C_m$.		
\item The subcategory $\U$ is closed under extensions.
\end{enumerate}
\end{proposition}
\begin{proof}
The equivalence of statements (1) and (2) follows from the computation of the middle term of an extension having indecomposable outer terms. We prove that (2) implies (3), the other implication is trivial. Since $\U$ is closed under extensions having indecomposable outer terms, by  Lemma \ref{lemma subcategory closed under extensions with indecomposable outer terms implies closed under extensions} $\U$ is closed under extensions of the form $a\lora e\lora b\lora \S a$ with $a\in\ind\C_m$, $b = \bigoplus_{i = 1}^n b_i$ where $b_1,\dots,b_n\in\ind\C_m$ are non-isomorphic to $\S a$ and are pairwise $\Hom$-orthogonal, and all the entries of $h = (h_1,\dots,h_n)$ are non-zero. By Remark \ref{remark factorization free hom orthogonal}, in particular $h_1,\dots,h_n$ are pairwise factorization free. Then, by Proposition \ref{proposition closure under extensions}, $\U$ is closed under extensions.
\end{proof}

\subsection{Extension-closed subcategories of $\overline{\mathcal{C}}_m$}

Here we characterise the extension-closed subcategories of $\ovl\C_m$. First we introduce the completed version of the $\PT$ condition. We refer to Proposition \ref{proposition Hom sets completion} for the computation of the $\Hom$-spaces of $\ovl\C_m$.

Recall that we identify the indecomposable objects of $\ovl\C_m$ with the arcs of $\ovl\Zc_m$.

\begin{definition}
Let $x,y\in\ind\ovl\C_m$ be such that $\Hom_{\ovl\C_m}(x,\S y)\cong \K$. The arcs of $\ind\ovl\C_m\setminus\{x,y\}$ which connect the endpoints of $x$ and $y$ are called \emph{Ptolemy arcs} of $x$ and $y$. We say that a subcategory $\X$ of $\ovl\C_m$ satisfies the \emph{completed Ptolemy condition}, $\ovl\PT$ condition for short, if it is closed under taking Ptolemy arcs.
\end{definition}

\begin{figure}[ht]
\centering
\includegraphics[height = 4cm]{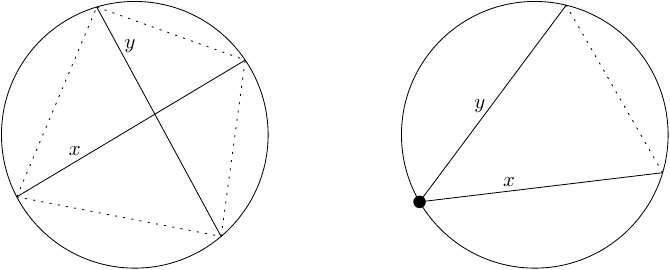}
\caption{The dotted arcs are the Ptolemy arcs of $x$ and $y$. On the left $x$ and $y$ cross, on the right they share one endpoint which is an accumulation point.}
\label{figure Ptolemy arcs completion}
\end{figure}

Figure \ref{figure Ptolemy arcs completion} provides an illustration of the Ptolemy arcs in $\ovl\C_m$.

The middle term of a non-split extension in $\ovl{\C}_m$ having indecomposable outer terms was computed in \cite[Section 3]{PY}.

\begin{proposition}\label{proposition extension closed subcategories completion}
Let $\X$ be a subcategory of $\ovl\C_m$. The following statements are equivalent.
\begin{enumerate}
\item The subcategory $\X$ satisfies the $\ovl\PT$ condition.
\item The subcategory $\X$ is closed under extensions of the form $x_1\lora c\lora x_2\lora \S x_1$ with $x_1,x_2\in\ind\ovl\C_m$.
\item The subcategory $\X$ is closed under extensions.
\end{enumerate}
\end{proposition}
\begin{proof}
The proof of the equivalence of (1) and (2) is straightforward and follows from \cite[Section 3]{PY}. The fact that (3) implies (2) is trivial. We prove that (2) implies (3). To this end, first we show that $\pi^{-1}\X$ is closed under extensions, and then that $\X$ is closed under extensions.

Assume that $\X$ is closed under extensions with indecomposable outer terms. Consider the preimage $\pi^{-1}\X$ in $\C_{2m}$. It is straightforward to check that $\pi^{-1}\X$ is an additive subcategory of $\C_{2m}$. We show that $\pi^{-1}\X$ is closed under extensions having indecomposable outer terms. Consider a triangle $a\lora e\lora b\lora \S a$ in $\C_{2m}$ with $a,b\in\ind\pi^{-1}\X$. Then $\pi a\lora \pi e\lora \pi b\lora \pi\S a$ is a triangle in $\ovl\C_m$, see \cite[Lemma 4.3.1]{K}. Moreover, from \cite[Proposition 3.10]{PY} it follows that $\pi a$ and $\pi b$ are either indecomposable objects or zero. From (2) we obtain that $\pi e\in\X$, i.e. $e\in\pi^{-1}\X$. This proves that $\pi^{-1}\X$ is closed under extensions having indecomposable outer terms. By Lemma \ref{lemma subcategory closed under extensions with indecomposable outer terms implies closed under extensions} we obtain that $\pi^{-1}\X$ is closed under extensions. 

Now we show that $\X$ is closed under extensions. Consider a triangle $x_1\lora c\lora x_2\lora \S x_1$ in $\ovl\C_m$ with $x_1,x_2\in\X$. Then there exists a triangle $a\lora e\lora b\lora \S a$ in $\C_{2m}$ whose image after $\pi$ is isomorhic in $\ovl\C_m$ to $x_1\lora c\lora x_2\lora \S x_1$. Thus $\pi a,\pi b\in\X$, i.e. $a,b\in\pi^{-1}\X$. Since $\pi^{-1}\X$ is closed under extensions, we have $e\in\pi^{-1}\X$ and then $c \cong \pi e \in\X$. This completes the proof.
\end{proof}

\section{Torsion pairs}\label{section torsion pairs}

Torsion pairs in $\C_m$ were classified in \cite[Section 4]{GHJ}, we classify the torsion pairs in $\ovl\C_m$. We recall from Proposition \ref{proposition iyama-yoshino torsion pairs} that a torsion pair $(\X,\Y)$ is uniquely determined by its torsion class $\X$, and therefore it is enough to classify the torsion classes.

\begin{theorem}\label{theorem classification of torsion pairs}
Let $\X$ be a subcategory of $\ovl\C_m$. Then $\X$ is a torsion class in $\ovl\C_m$ if and only if $\X$ satisfies the $\ovl{\mathrm{PC}}$ conditions and the $\ovl{\mathrm{PT}}$ condition. Moreover, there is an inclusion preserving bijection.	
\begin{align*}
\left\{\parbox{4cm}{\centering Torsion-classes in $\ovl\C_m$}\right\} & \longleftrightarrow\left\{\parbox{7.5cm}{\centering Extension-closed subcategories $\U\subseteq\C_{2m}$ such that $\D\subseteq\U$ and $\U\cap\A$ is precovering} \right\}\\
\X & \longmapsto \pi^{-1}\X\\
\pi\U & \longmapsfrom \U
\end{align*}
\end{theorem}
\begin{proof}
The first statement follows directly from Proposition \ref{proposition iyama-yoshino torsion pairs}, Theorem \ref{theorem classification precovering subcategories}, and Proposition \ref{proposition extension closed subcategories completion}. The bijection follows from Proposition \ref{proposition bijection between extension closed subcategories in the verdier quotient} and Theorem \ref{theorem classification precovering subcategories}.
\end{proof}

We have the following corollaries of Theorem \ref{theorem classification of torsion pairs} and which will be useful in Section \ref{section t-structures} and Section \ref{section co-t-structures} for classifying t-structures and co-t-structures.

\begin{corollary}\label{corollary t-structures}
Let $\X$ be a subcategory of $\ovl\C_m$. Then $\X$ is the aisle of a t-structure in $\ovl\C_m$ if and only if $\X$ satisfies the $\ovl{\mathrm{PC}}$ conditions, the $\ovl{\mathrm{PT}}$ condition, and $\X$ is closed under clockwise rotations. Moreover, there is an inclusion preserving bijection.	
\begin{align*}
\left\{\parbox{5cm}{\centering Aisles of t-structures in $\ovl\C_m$}\right\} & \longleftrightarrow\left\{\parbox{7.5cm}{\centering Suspended subcategories $\U\subseteq\C_{2m}$ such that $\D\subseteq\U$ and $\U\cap\A$ is precovering} \right\}\\
\X & \longmapsto \pi^{-1}\X\\
\pi\U & \longmapsfrom \U
\end{align*}
\end{corollary}

\begin{corollary}\label{corollary co-t-structures}
Let $\X$ be a subcategory of $\ovl\C_m$. Then $\X$ is the aisle of a co-t-structure in $\ovl\C_m$ if and only if $\X$ satisfies the $\ovl{\mathrm{PC}}$ conditions, the $\ovl{\mathrm{PT}}$ condition, and $\X$ is closed under anticlockwise rotations. Moreover, there is an inclusion preserving bijection.	
\begin{align*}
\left\{\parbox{5.5cm}{\centering Aisles of co-t-structures in $\ovl\C_m$}\right\} & \longleftrightarrow\left\{\parbox{7.5cm}{\centering Co-suspended subcategories $\U\subseteq\C_{2m}$ such that $\D\subseteq\U$ and $\U\cap\A$ is precovering} \right\}\\
\X & \longmapsto \pi^{-1}\X\\
\pi\U & \longmapsfrom \U
\end{align*}
\end{corollary}

\section{t-structures}\label{section t-structures}

The t-structures in $\C_m$ were classified in \cite{GZ} for $m\geq 1$, and in \cite{N} and \cite{CZZ} for the cases $m = 1$ and $m = 2$. Here we classify the t-structures in $\ovl\C_m$. We start by classifying the aisles, then we compute the co-aisles and the hearts. Finally, we classify the bounded and non-degenerate t-structures. 

\subsection{Aisles of t-structures}\label{section aisles of t-structures}

We recall the classification of the aisles of the t-structures in $\C_m$ from \cite[Section 4]{GZ} in terms of decorated non-crossing partitions.

\begin{definition}[{\cite[Definition 4.5]{GZ}}]\label{definition GZ decorated non-crossing partition}
A \emph{decorated non-crossing partition} of $[m]$ is a pair $(\P,X)$ given by a non-crossing partition $\P$ of $[m]$ and an $m$-tuple $X = (x_p)_{p\in[m]}$ where for each $p\in[m]$
\[
x_p \in 
\begin{cases}
[p,p^+) & \text{if $\{p\}\in\P$},\\
(p,p^+] & \text{if $p,p^+\in B$ for some block $B\in\P$},\\
(p,p^+) & \text{otherwise}.
\end{cases}
\]
\end{definition}

The following theorem gives a classification of the aisles of the t-structures in $\C_m$.

\begin{theorem}[{\cite[Theorem 4.6]{GZ}}]\label{theorem GZ classification aisles t-structures}
There is a bijection between the decorated non-crossing partitions of $[m]$ and the aisles of t-structures in $\C_m$.
\end{theorem}

In order to classify the aisles of the $t$-structures in $\ovl\C_m$ we need to adapt Definition \ref{definition GZ decorated non-crossing partition} in our setting.

\begin{definition}\label{definition half-decorated non-crossing partition}
A \emph{half-decorated non-crossing partition} of $[m']\cup[m]$ is a pair $(\P,X)$ given by a non-crossing partition $\P$ of $[m']\cup[m]$ and an $m$-tuple $X = (x_p)_{p\in[m]}$ such that for each $p\in[m]$ 
\[
x_p\in
\begin{cases}
[p,p^+) & \text{if $\{p\}\in\P$},\\
(p,p^+] & \text{if $p,p^+\in B$ for some block $B\in\P$,}\\
(p,p^+) & \text{otherwise.}
\end{cases}
\]
\end{definition}

\begin{figure}[ht]
\centering
\includegraphics[height = 4cm]{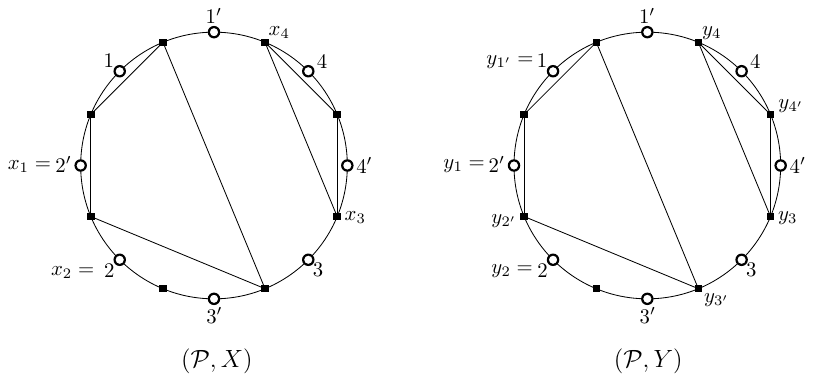}
\caption{On the left $(\P,X)$ is a half-decorated non-crossing partition of $[4']\cup[4]$, and on the right $(\P,Y)$  is a decorated non-crossing partition of $[4']\cup[4]$, where $\P = \{\{1',1,2',3'\},\{2\},\{3,4',4\}\}$. Note that $(\P,Y)$ does not correspond to $(\P,X)$ according to Remark \ref{remark half-decorated non-crossing partitions and decorated non-crossing partitions}.}
\label{figure half-decorated non-crossing partition non-crossing partition}
\end{figure}

Note that a decorated non-crossing partition $(\P,X)$ of $[m']\cup[m]$ is not a half-decorated non-crossing partition of $[m']\cup[m]$ because $X$ is a $2m$-tuple, not an $m$-tuple. Figure \ref{figure half-decorated non-crossing partition non-crossing partition} gives examples of a half-decorated non-crossing partition and a decorated non-crossing partition of $[m']\cup[m]$.

\begin{remark}\label{remark half-decorated non-crossing partitions and decorated non-crossing partitions}
Given a half-decorated non-crossing partition of $[m']\cup[m]$, we can obtain a decorated non-crossing partition by adding the decoration $z_p^0$ for each $p\in[m']$. Conversely, given such a decorated non-crossing partition we obtain a half-decorated non-crossing partition by removing the decorations $x_p$ for each $p\in[m']$. These assignments are mutually inverse. Therefore, we have a bijection  between the half-decorated non-crossing partitions of $[m']\cup[m]$ and the decorated non-crossing partitions of $[m']\cup[m]$ of the form $(\P,X)$ with $X = (x_p)_{p\in[m']\cup[m]}$ such that $x_p = z_p^0$ for each $p\in[m']$.
\end{remark}

The main result of this section is the following analogue of Theorem \ref{theorem GZ classification aisles t-structures}. The notation employed in the statement will be defined in Definition \ref{definition almost aisle of t-structure} and Definition \ref{definition from almost aisle of t-structure to half-decorated non-crossing partition}.

\begin{theorem}\label{theorem t-structures}
The following is a bijection. 
\begin{align*}
\left\{\parbox{5cm}{\centering Half-decorated non-crossing partitions of $[m']\cup[m]$}\right\} & \longleftrightarrow \left\{\parbox{5cm}{\centering Aisles of t-structures in $\ovl\C_m$} \right\}\\
(\P,X) & \longmapsto \pi\U_{(\P,X)}\\
(\P_{\pi^{-1}\X},X_{\pi^{-1}\X}) & \longmapsfrom \X
\end{align*}
\end{theorem}

From Remark \ref{remark half-decorated non-crossing partitions and decorated non-crossing partitions} and Theorem \ref{theorem t-structures} it follows that the aisles of the t-structures in $\ovl\C_m$ are in bijection with certain aisles of t-structures in $\C_{2m}$ (namely those corresponding to the decorated non-crossing partitions $(\P,X)$ of $[m']\cup[m]$ with $X = (x_p)_{p\in[m']\cup[m]}$ such that $x_p = z_p^0$ for each $p\in[m']$).

\begin{figure}[ht]
\centering
\includegraphics[height = 5cm]{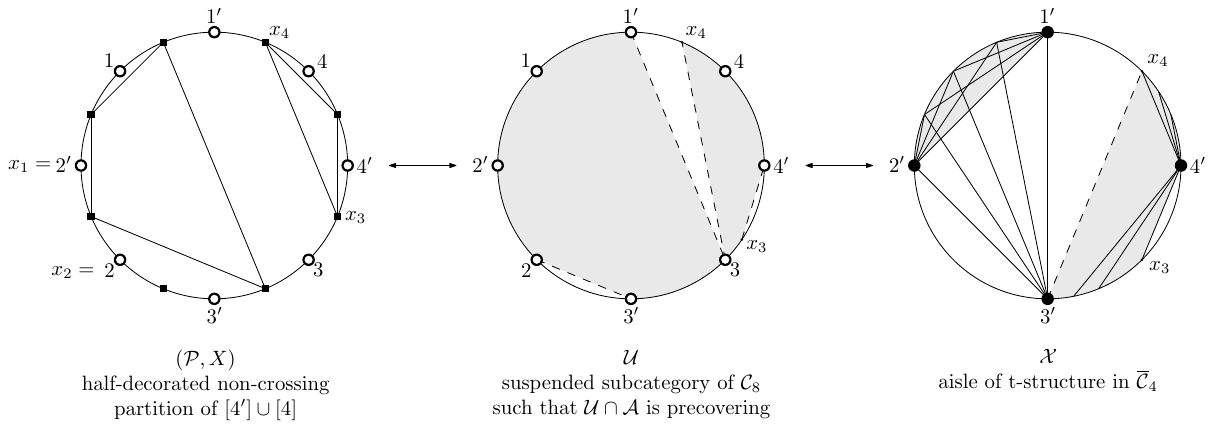}
\caption{Illustration of how to obtain the aisle of a t-structure of $\ovl\C_m$ from a half-decorated non-crossing partition of $[m']\cup[m]$.}
\label{figure half decorated non-crossing partition aisle t-structure}
\end{figure}

To prove Theorem \ref{theorem t-structures} we take an intermediate step through $\C_{2m}$. From Corollary \ref{corollary t-structures} the aisles of t-structures in $\ovl\C_m$ are in bijection with the suspended subcategories $\U$ of $\C_{2m}$ such that $\D\subseteq\U$ and $\U\cap\A$ is precovering. These can be regarded as ``almost aisles" of t-structures in $\C_{2m}$ and are classified in terms of half-decorated non-crossing partitions of $[m']\cup[m]$, see Proposition \ref{proposition classification t-structures}. The aisles of the t-structures in $\ovl\C_m$ are then obtained by localising the ``almost aisles" in $\C_{2m}$. Figure \ref{figure half decorated non-crossing partition aisle t-structure} illustrates this process.

The following proposition classifies the ``almost aisles" of t-structures in $\C_{2m}$.

\begin{proposition}\label{proposition classification t-structures}
The following is a bijection.
\begin{align*}
\left\{\parbox{5cm}{\centering Half-decorated non-crossing partitions of $[m']\cup[m]$}\right\} & \longleftrightarrow \left\{\parbox{7.5cm}{\centering Suspended subcategories $\U\subseteq\C_{2m}$ such that $\D\subseteq\U$ and $\U\cap\A$ is precovering} \right\}\\
\a\colon (\P,X) & \longmapsto \U_{(\P,X)}\\
(\P_{\U},X_{\U}) & \longmapsfrom \U\colon\b
\end{align*}
\end{proposition}

The rest of this section is devoted to prove Proposition \ref{proposition classification t-structures}. We start by defining the assignments of the maps $\a$ and $\b$. 

\begin{definition}\label{definition almost aisle of t-structure}
Let $(\P,X)$ be a half-decorated non-crossing partition of $[m']\cup[m]$. We define
\[
\U_{(\P,X)} = \add\bigsqcup_{B\in\P} \left\{(u_1,u_2)\in\ind\C_{2m}\middle| u_1,u_2\in\left(\bigcup_{p\in B\cap[m]}(p,x_p]\right)\cup \left(\bigcup_{p\in B\cap[m']}\Z^{(p)}\right)\right\},
\]
where we use the following convention: for $p\in[m]$, if $x_p = p$ then $(p,x_p] = \emptyset$, and if $x_p = p^+$ then $(p,x_p] = \Z^{(p)}$.
\end{definition}

We check that the map $\a$ is well defined.

\begin{proposition}\label{proposition from half decorated non-crossing partition to almost aisle of t-structure}
Let $(\P,X)$ be a half-decorated non-crossing partition of $[m']\cup[m]$. Then $\U_{(\P,X)}$ is a suspended subcategory of $\C_{2m}$ such that $\D\subseteq\U_{(\P,X)}$ and $\U_{(\P,X)}\cap\A$ is precovering.
\end{proposition}
\begin{proof}
We show that $\D\subseteq \U_{(\P,X)}$. Consider $d = (d_1,d_2)\in\ind\D$, then $d_1,d_2\in\Z^{(p)}$ for some $p\in[m']$. Since $p\in B$ for some block $B\in\P$, then $\Z^{(p,p)}\subseteq\ind\U_{(\P,X)}$, and then $d\in\U_{(\P,X)}$. Moreover, it is straightforward to check that $\S\U_{(\P,X)}\subseteq\U_{(\P,X)}$. For showing that $\U_{(\P,X)}$ is extension-closed, we can proceed as in the argument of \cite[Proposition 4.8]{GZ}. 
	
Now we show that $\U_{(\P,X)}\cap\A$ is precovering. Let $\wt{X}$ be the $2m$-tuple $\wt{X} = (\wt x_p)_{p\in[m']\cup[m]}$ where for each $p\in[m']\cup[m]$
\[
\wt x_p = 
\begin{cases}
z_p^0 & \text{if $p\in[m']$,}\\
x_p & \text{if $p\in[m]$.}
\end{cases}
\]
By Remark \ref{remark half-decorated non-crossing partitions and decorated non-crossing partitions} $(\P,\wt X)$ is a decorated non-crossing partition of $[m']\cup[m]$ and we can associate to it the aisle of a t-structure $\U'$ in $\C_{2m}$, see Theorem \ref{theorem GZ classification aisles t-structures}. By \cite[p. 986]{GZ} we have that
\[
\U' = \add\bigsqcup_{B\in\P}\left\{(u_1,u_2)\in\ind\C_{2m}\middle| u_,u_2\in\bigcup_{p\in B}(p,\wt x_p]\right\}.
\]
It is straightforward to check that $\U' = \U_{(\P,X)}\cap\A$. Thus, $\U_{(\P,X)}\cap\A$ is the aisle of a t-structure in $\C_{2m}$, and in particular it is precovering. This concludes the argument.
\end{proof}  

Now we define the map $\b$. To this end, given an ``almost aisle" in $\C_{2m}$, we define an equivalence relation $\sim_{\U}$ on the set $[m']\cup[m]$ in the same way as in \cite[Section 4.1]{GZ}. The same argument of \cite[Lemma 4.10]{GZ} shows that $\sim_{\U}$ is an equivalence relation.

\begin{definition}\label{definition equivalence relation}
Let $\U$ be a suspended subcategory of $\C_{2m}$ such that $\D\subseteq\U$ and $\U\cap\A$ is precovering. The relation $\sim_{\U}$ on the set $[m']\cup[m]$ is defined as follows: for any $p,q\in[m']\cup[m]$ we have that $p\sim_{\U}q$ if and only if $p = q$ or there exists an arc of $\U$ with an endpoint in $\Z^{(p)}$ and the other in $\Z^{(q)}$. 
\end{definition}

\begin{definition}\label{definition from almost aisle of t-structure to half-decorated non-crossing partition}
Keeping the assumptions and notation of Definition \ref{definition equivalence relation}, we define $\P_{\U}$ to be the partition of $[m']\cup[m]$ given by the equivalence classes of $\sim_{\U}$. For each $p\in[m]$ we define 
\[
x_p = \sup\{z\in\Z^{(p)}\mid\text{there exists an arc of $\U$ with an endpoint equal to $z$}\}.
\]	
We denote by $X_{\U}$ the $m$-tuple $X_{\U} = (x_p)_{p\in[m]}$.
\end{definition}

With Proposition \ref{proposition from almost aisle of t-structure to half decorated non-crossing partition} we will show that $(\P_{\U},X_{\U})$ is a half-decorated non-crossing partition of $[m']\cup[m]$. The following remark and lemmas are useful for that purpose. 

\begin{remark}
Consider a suspended subategory $\U$ of $\C_{2m}$ such that $\D\subseteq\U$ and $\U\cap\A$ is precovering. We observe that $\U\cap\A$ is the aisle of a t-structure in $\C_{2m}$. We denote by $(\P_{\U\cap\A},X_{\U\cap\A})$ the decorated non-crossing partition associated to $\U\cap\A$, as defined in \cite[Definition 4.11]{GZ}. The following lemmas relate $(\P_{\U\cap\A},X_{\U\cap\A})$, defined in \cite{GZ}, and $(\P_{\U},X_{\U})$, defined above. We recall that $\P_{\U\cap\A}$ is the set of equivalence classes of $[m']\cup[m]$ under the equivalence relation $\sim_{\U\cap\A}$. 
\end{remark}

\begin{lemma}\label{lemma partitions}
Let $\U$ be a suspended subcategory of $\C_{2m}$ such that $\D\subseteq\U$ and $\U\cap\A$ is precovering. Let $(\P_{\U\cap\A},X_{\U\cap\A})$ be the decorated non-crossing partition associated to $\U\cap\A$. Then $\P_{\U} = \P_{\U\cap\A}$.
\end{lemma}
\begin{proof}
We show that for any $p,q\in[m']\cup[m]$ we have that $p\sim_{\U}q$ if and only if $p\sim_{\U\cap\A}q$. It is straightforward to check that if $p\sim_{\U\cap\A}q$ then $p\sim_{\U}q$. Assume that $p\sim_{\U}q$. If $p = q$ then the claim is trivial. If $p\neq q$ then there exists $u\in\ind\U$ having one endpoint in $\Z^{(p)}$ and the other endpoint in $\Z^{(q)}$. Note that there exists $n\geq 0$ such that $\S^n u\in\A$ and then, since $\S^n\U\subseteq\U$, we obtain that $u\in\U\cap\A$. Then we have that $p\sim_{\U\cap\A}q$. This concludes the argument.
\end{proof}

\begin{lemma}\label{lemma decorations}
Let $\U$ be a subcategory of $\C_{2m}$ as in Lemma \ref{lemma partitions}, $(\P_{\U\cap\A}, X_{\U\cap\A})$ be the decorated non-crossing partition associated to $\U\cap\A$. Denote $X_{\U} = (x_p)_{p\in[m]}$ and $X_{\U\cap\A} = (\wt x_p)_{p\in[m']\cup[m]}$. Then for each $p\in[m']\cup[m]$ we have that
\[
\wt x_p = 
\begin{cases}
z_p^0 & \text{if $p\in[m']$,}\\
x_p & \text{if $p\in[m]$.}
\end{cases}
\]
\end{lemma}
\begin{proof}
Let $p\in[m']\cup[m]$. We recall that by construction, see \cite[Section 4.1]{GZ}, we have that
\[
\wt x_p = \sup\{z\in\Z^{(p)}\mid \text{there exists an arc of $\U\cap\A$ with an endpoint equal to $z$}\}.
\]
If $p\in[m']$, there exists an arc of $\U\cap\A$ with an endpoint equal to $z_p^0$. Let $z\in\Z^{(p)}$ be such that $z>z_p^0$. There is no arc of $\A$, and therefore no arc of $\U\cap\A$, with an endpoint equal to $z$. Thus, $\wt x_p = z_p^0$. Now consider $p\in[m]$, we show that $\wt x_p = x_p$. We divide the argument into claims.
	
\emph{Claim 1.} Let $z\in\Z^{(p)}$. If there exists an arc of $\U$ with an endpoint equal to $z$, then there exists an arc of $\U\cap\A$ with an endpoint equal to $z$.
	
Assume that there exists $u\in\ind\U$ having an endpoint equal to $z$. If $u\in\U\cap\A$, then we have the claim. Now assume that $u\in\U$ and $u\notin\A$. We denote $u = (u_1,u_2)$. We assume that $u_1 = z$, the other case is analogous. Since $u\notin\A$, we have that $u_2\in\Z^{(q)}$ for some $q\in[m']$ and $u_2>z_p^0$. Then we are in the situation of Figure \ref{figure lemma decorations}.
\begin{figure}[ht]
\centering
\includegraphics[height = 4cm]{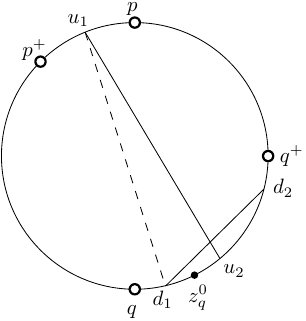}
\caption{Illustration of the argument of Claim 1.}
\label{figure lemma decorations}
\end{figure}
	
Consider $d = (d_1,d_2)\in\Z^{(q,q)}$ with $d_1\leq z_q^0<u_2<d_2$. Since $d\in\D\subseteq\U$, $u$ and $d$ are crossing, and $\U$ is extension-closed, we obtain that $(z,d_1) = (u_1,d_1)\in\U$. Moreover, $(z,d_1)\in\A$. This concludes the argument of Claim 1.
	
\emph{Claim 2.} If $\wt x_p = p$ then $x_p = p$.
	
Assume that $\wt x_p = p$, i.e. there is no $z\in\Z^{(p)}$ such that there is an arc of $\U\cap\A$ with an endpoint equal to $z$. By Claim 1 we have that there is no $z\in\Z^{(p)}$ such that there is an arc of $\U$ with an endpoint equal to $z$, i.e. $x_p = p$. This concludes the argument of Claim 2.
	
\emph{Claim 3.} If $\wt x_p = p^+$ then $x_p = p^+$.
	
The proof is straightforward.
	
\emph{Claim 4.} $\wt x_p = x_p$.
	
If $\wt x_p = p$ or $\wt x_p = p^+$ then the claim follows from Claim 2 and Claim 3. Assume that there exists $z\in\Z^{(p)}$ such that $\wt x_p = z$. As a consequence, there is an arc of $\U\cap\A$ with an endpoint equal to $z$, and then there is an arc of $\U$ with an endpoint equal to $z$. Moreover, for any $z'\in\Z^{(p)}$ such that $z'>z$ there is no arc of $\U\cap\A$ with an endpoint equal to $z'$. Then, by Claim 1, for any $z'\in\Z^{(p)}$ such that $z'>z$ there is no arc of $\U$ with an endpoint equal to $z'$. Thus, $x_p = z$. This concludes the argument of Claim 4.
	
We can conclude that $\wt x_p = z_p^0$ for each $p\in[m]'$, and $\wt x_p = x_p$ for each $p\in[m]$.
\end{proof}

Now we can prove that the map $\b$ of Proposition \ref{proposition classification t-structures} is well defined.

\begin{proposition}\label{proposition from almost aisle of t-structure to half decorated non-crossing partition}
Let $\U$ be a suspended subcategory of $\C_{2m}$ such that $\D\subseteq\U$ and $\U\cap\A$ is precovering. Then $(\P_{\U},X_{\U})$ is a half-decorated non-crossing partition.
\end{proposition}
\begin{proof}
We check that $(\P_{\U},X_{\U})$ satisfies the conditions of Definition \ref{definition half-decorated non-crossing partition}. Consider $\U\cap\A$, which is the aisle of a t-structure, and its associated decorated non-crossing partition $(\P_{\U\cap\A},X_{\U\cap\A})$, as defined in \cite{GZ}. We recall that $\P_{\U}$ is a partition of $[m']\cup[m]$ and that, from Lemma \ref{lemma partitions}, $\P_{\U\cap\A} = \P_{\U}$. As a consequence, $\P_{\U}$ is a non-crossing partition of $[m']\cup[m]$. Now, we denote the decorations by $X_{\U} = (x_p)_{p\in[m]}$ and $X_{\U\cap\A} = (\wt x_p)_{p\in[m']\cup[m]}$. By Lemma \ref{lemma decorations} we have that $x_p = \wt x_p$ for each $p\in[m]$. We conclude that $(\P_{\U},X_{\U})$ is a half-decorated non-crossing partition of $[m']\cup[m]$.
\end{proof}

Given $\U$ a suspended subcategory of $\C_{2m}$ such that $\D\subseteq\U$ and $\U\cap\A$ is precovering, the following lemma shows that any shift of $\U$ has the same properties of $\U$. This fact will be useful in the proof of Proposition \ref{proposition classification t-structures}.

\begin{lemma}\label{lemma Sn U}
Let $\U$ be a suspended subcategory of $\C_{2m}$ such that $\D\subseteq\U$ and $\U\cap\A$ is precovering. Consider the associated half-decorated non-crossing partition $(\P_{\U},X_{\U})$ with $X_{\U} = (x_p)_{p\in[m]}$. The following statements hold. 
\begin{enumerate}
\item For any $n\in\Z$ the subcategory $\S^n\U$ of $\C_{2m}$ is suspended, $\D\subseteq\S^n\U$, and $\S^n\U\cap\A$ is precovering. 
\item Consider $(\P_{\S^n\U},X_{\S^n\U})$. Then $\P_{\S^n\U} = \P_{\U}$ and $X_{\S^n\U} = (x_p-n)_{p\in[m]}$.
\end{enumerate}
\end{lemma}
\begin{proof}
First we prove statement (1), statement (2) follows by construction, see Definition \ref{definition almost aisle of t-structure}. It is straightforward to check that $\S^n\U$ is extension-closed and contains $\D$, we show that $\S^n\U\cap\A$ is precovering. By Proposition \ref{proposition bijection between extension closed subcategories in the verdier quotient} we have that $\pi^{-1}\pi\U = \U$, and, since $\pi^{-1}\pi\U\cap\A = \U\cap\A$ is precovering, by Theorem \ref{theorem classification precovering subcategories} we have that $\pi\U$ is precovering in $\ovl\C_m$. Now fix $n\in\Z$. Since $\pi\U$ is precovering, then $\S^n\pi\U$ is precovering in $\ovl\C_m$. As a consequence, again by Theorem \ref{theorem classification precovering subcategories}, we have that $\pi^{-1}\S^n\pi\U\cap\A$ is precovering. Since $\pi^{-1}\S^n\pi\U\cap\A = \pi^{-1}\pi\S^n\U\cap\A = \S^n\U\cap\A$, we obtain that $\S^n\U\cap\A$ is precovering.
\end{proof}

Finally, we can prove Proposition \ref{proposition classification t-structures}.

\begin{proof}[Proof of Proposition \ref{proposition classification t-structures}]
By Proposition \ref{proposition from half decorated non-crossing partition to almost aisle of t-structure} and Proposition \ref{proposition from almost aisle of t-structure to half decorated non-crossing partition} the maps $\a$ and $\b$ are well defined, we prove that they are mutually inverse. We divide the proof into steps.
	
\emph{Step 1.} The map $\b$ is injective.
	
Let $\U$ and $\U'$ be suspended subcategories of $\C_{2m}$ such that $\D\subseteq\U$, $\D\subseteq\U'$, with $\U\cap\A$ and $\U'\cap\A$ precovering. Assume that $(\P_{\U},X_{\U}) = (\P_{\U'}, X_{\U'})$, we show that $\U = \U'$.
	
First we show that $\S^n\U\cap\A = \S^n\U'\cap\A$ for each $n\in\Z$. By Lemma \ref{lemma partitions} and Lemma \ref{lemma decorations} and Lemma \ref{lemma Sn U}, $(\P_{\S^n\U\cap\A}, X_{\S^n\U\cap\A}) = (\P_{\S^n\U'\cap\A},X_{\S^n\U'\cap\A})$. By Theorem \ref{theorem GZ classification aisles t-structures} we obtain that $\S^n\U\cap\A = \S^n\U'\cap\A$.
	
Now we show that $\U\subseteq\U'$, the other inclusion can be obtained in the same way. Consider $u\in\ind\U$, we have that $\S^n u\in\A$ for some $n\geq 0$. Then $\S^n u\in\S^n\U\cap\A$. Since $\S^n\U\cap\A = \S^n\U'\cap\A$, we have that $\S^n u\in\S^n\U'\cap\A$ and then $u\in\U'$. We obtain that $\U = \U'$ and this concludes the argument of Step 1.
	
\emph{Step 2.} We show that $\b\a = \id$.
	
Let $(\P,X)$ be a half-decorated non-crossing partition of $[m']\cup[m]$. Let $\U_{(\P,X)}$ be the associated subcategory of $\C_{2m}$, which we denoted by $\U$. Let $(\P_{\U},X_{\U})$ be the half-decorated non-crossing partition associated to $\U$. We show that $(\P,X) = (\P_{\U},X_{\U})$. 
	
Showing the equality $\P = \P_{\U}$ is equivalent to show that for any $p,q\in[m']\cup[m]$ we have that $p\sim_{\U}q$ if and only if $p,q\in B$ for some block $B\in\P$. This follows directly from Definition \ref{definition almost aisle of t-structure} and Definition \ref{definition equivalence relation}. Now we show that $X = X_{\U}$. We denote $X = (x_p)_{p\in[m]}$ and $X_{\U} = (y_p)_{p\in[m]}$. By construction, see Definition \ref{definition almost aisle of t-structure} and Definition \ref{definition from almost aisle of t-structure to half-decorated non-crossing partition}, we have the following equalities for each $p\in[m]$.
\[
y_p = \sup\{z\in\Z^{(p)}\mid\text{there exists an arc of $\U = \U_{(\P,X)}$ with an endpoint equal to $z$}\} = x_p
\]
Therefore we have that $(\P,X) = (\P_{\U},X_{\U})$. This concludes the argument of Step 2.
	
We can conclude that the maps $\a$ and $\b$ are mutually inverse.
\end{proof}

\subsection{Co-aisles of t-structures}\label{section co-aisles t-structures}

\begin{figure}[ht]
\centering
\includegraphics[height = 5cm]{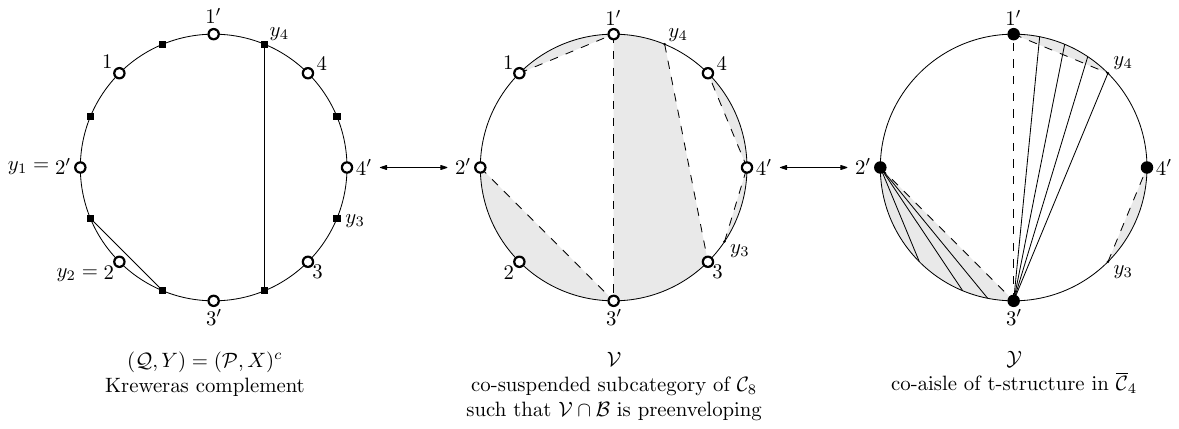}
\caption{Illustration of how to obtain the co-aisle of the aisle of Figure \ref{figure half decorated non-crossing partition aisle t-structure}.}
\label{figure complement half-decorated non-crossing partition co-aisle t-structure}
\end{figure}

From Theorem \ref{theorem t-structures} we have a classification of the aisles of the t-structures in $\ovl\C_m$, now we compute the corresponding co-aisles in terms of non-crossing partitions. As before, we take an intermediate step through $\C_{2m}$. Given a half-decorated non-crossing partition $(\P,X)$ of $[m']\cup[m]$, we consider its complement $(\P,X)^c = (\Q,Y)$, where $\Q = \P^c$ is the Kreweras complement of $\P$, see Section \ref{section non-crossing partitions}. With a computation similar to \cite[Section 4.2]{GZ}, $(\Q,Y)$ corresponds to a subcategory $\V$ of $\C_{2m}$. This is a co-suspended subcategory of $\C_{2m}$ such that $\D\subseteq\V$ and $\V\cap\B$ is preenveloping, therefore $\V$ can be thought as an ``almost co-aisle" in $\C_{2m}$. From such $\V$ we obtain the corresponding co-aisle in $\ovl\C_m$ after localising. Figure \ref{figure complement half-decorated non-crossing partition co-aisle t-structure} illustrates this process. 

\begin{definition}\label{definition complement half decorated non-crossing partition}
Let $(\P,X)$ be a half-decorated non-crossing partition of $[m']\cup[m]$ with $X = (x_p)_{p\in[m]}$. We define the \emph{complement}, $(\P,X)^c$, of $(\P,X)$ to be the pair $(\Q,Y)$ where $\Q = \P^c$ is the Kreweras complement of $\P$, and $Y$ is the $m$-tuple $Y = X-1 = (x_p-1)_{p\in[m]}$ for each $p\in[m]$. 
\end{definition}

We describe how to obtain an ``almost co-aisle" of t-structure in $\C_{2m}$ from the complement of a half-decorated non-crossing partition of $[m']\cup[m]$.

\begin{definition}\label{definition almost co-aisle of t-structure}
Let $(\P,X)$ be a half-decorated non-crossing partition of $[m']\cup[m]$ and let $(\Q, Y) = (\P,X)^c$. We define
\[
\V_{(\Q,Y)} = \add\bigsqcup_{B\in\Q}\left\{(v_1,v_2)\in\ind\C_{2m}\middle| v_1,v_2\in\left(\bigcup_{p\in B\cap[m]} [y_{p},p^+)\right)\cup \left(\bigcup_{p\in B\cap[m']}\Z^{(p)}\right)\right\}.
\]
\end{definition}

Consider the complement $(\Q,Y)$ of a half-decorated non-crossing partition of $[m']\cup[m]$. The following lemmas and remark establish some properties of the subcategory $\V_{(\Q,Y)}$ of $\C_{2m}$. The first is analogous to Proposition \ref{proposition from half decorated non-crossing partition to almost aisle of t-structure}.

\begin{lemma}\label{lemma V is extension-closed and contains D}
Let $(\P,X)$ be a half-decorated non-crossing partition of $[m']\cup[m]$ and let $(\Q,Y) = (\P,X)^c$. Then $\V_{(\Q,Y)}$ is co-suspended and contains $\D$.
\end{lemma}
\begin{proof}
The proof is analogous to the argument of Proposition \ref{proposition from half decorated non-crossing partition to almost aisle of t-structure}.
\end{proof}

\begin{lemma}\label{lemma V intersected with B is equal to the co-aisle of U intersected with A}
Let $(\P,X)$ be a half-decorated non-crossing partition of $[m']\cup[m]$ and let $(\Q,Y) = (\P,X)^c$. Then $\V_{(\Q,Y)}\cap\B = \left(\U_{(\P,X)}\cap\A\right)^{\perp}$. 
\end{lemma}
\begin{proof}
For each $p\in[m']\cup[m]$ we define
\[
\wt y_p = 
\begin{cases}
y_p & \text{ if $p\in[m]$,}\\
w_p^0 & \text{ if $p\in[m']$}
\end{cases}
\]
where we recall from Section \ref{section preenveloping subcategories} that $w_p^0 = z_p^0-1$. It is straightforward to check that
\[
\V_{(\Q,Y)}\cap \B= \add\bigsqcup_{B\in\Q}\left\{(v_1,v_2)\in\ind\C_{2m}\middle| v_1,v_2\in\bigcup_{p\in B}[\wt{y_p},p^+)\right\}.
\]
Moreover, from \cite[Corollary 4.14]{GZ} the right hand side is equal to $\left(\U_{(\P,X)}\cap\A\right)^{\perp}$.
\end{proof}

\begin{remark}\label{remark U intersected with A and V intersected with B t-structure}
Let $\U_{(\P,X)}$ and $\V_{(\Q,Y)}$ be as in Lemma \ref{lemma V intersected with B is equal to the co-aisle of U intersected with A}. Since $\U_{(\P,X)}\cap\A$ is precovering and suspended, by Proposition \ref{proposition iyama-yoshino torsion pairs} $\left(\U_{(\P,X)}\cap\A, \V_{(\Q,Y)}\cap\B\right)$ is a t-structure in $\C_{2m}$.
\end{remark}

The following lemma and proposition show that an ``almost co-aisle" in $\C_{2m}$ is, after localising, the co-aisle of a t-structure in $\ovl\C_m$.

\begin{lemma}\label{lemma pi V is included in X perp}
Let $\X$ be the aisle of a t-structure in $\ovl\C_m$, let $(\P,X)$ be the half decorated non-crossing partition associated to $\X$, and let $(\Q,Y) = (\P,X)^c$. Then $\pi\V_{(\Q,Y)}\subseteq \X^{^{\perp}}$.
\end{lemma}
\begin{proof}
Assume that there exist $x\in\ind\X$ and $y\in\ind\pi\V_{(\Q,Y)}$ such that $\Hom_{\ovl\C_m}(x,y)\cong\mathbb{K}$. Note that there exists $y'\in\ind\B$ such that $\pi(y')\cong y$. Then $y'\in\pi^{-1}\pi\V_{(\Q,Y)}$ and, by Proposition \ref{proposition bijection between extension closed subcategories in the verdier quotient} and Lemma \ref{lemma V is extension-closed and contains D}, we have that $y'\in\ind\V_{(\Q,Y)}\cap\B$. We define $\U = \pi^{-1}\X$, we have that $\U = \U_{(\P,X)}$. Now, by Lemma \ref{lemma A aisle t-structure} and Lemma \ref{lemma hom sets A-cover} we have that there exists $x'\in\ind\A$ such that $\pi(x')\cong x$, and then $x'\in\ind\U\cap\A$, and $\Hom_{\C_{2m}}(x',y')\cong\mathbb{K}$. Since $x\in\ind\U\cap\A$ and $y'\in\ind\V_{(\Q,Y)}\cap\B$, this gives a contradiction with Lemma \ref{lemma V intersected with B is equal to the co-aisle of U intersected with A}. Then we can conclude that $\Hom_{\ovl\C_m}(\X,\pi\V_{(\Q,Y)}) = 0$.
\end{proof}

\begin{proposition}\label{proposition compute the co-aisle of a t-structure}
Let $(\X,\Y)$ be a t-structure in $\ovl\C_m$, $\U = \pi^{-1}\X$, $(\P,X)$ be its associated half-decorated non-crossing partition, and $(\Q,Y) = (\P,X)^c$. Then 
\[
\Y = \pi\V_{(\Q,Y)} = \pi\left(\V_{(\Q,Y)}\cap\B\right) = \pi\left((\U\cap\A)^{^{\perp}}\right).
\]
\end{proposition}
\begin{proof}
First we show that $\pi\V_{(\Q,Y)} = \pi\left(\V_{(\Q,Y)}\cap\B\right)$. The inclusion $\pi\left(\V_{(\Q,Y)}\cap\B\right)\subseteq\pi\V_{(\Q,Y)}$ is straightforward. We show the other inclusion. Consider $y\in\ind\pi\V_{(\Q,Y)}$, then there exists $y'\in\ind\B$ such that $\pi(y')\cong y$. Since $y\in\pi\V_{(\Q,Y)}$, we have that $y'\in\pi^{-1}\pi\V_{(\Q,Y)}$. By Proposition \ref{proposition bijection between extension closed subcategories in the verdier quotient} and Lemma \ref{lemma V is extension-closed and contains D} we have that $\pi^{-1}\pi\V_{(\Q,Y)} = \V_{(\Q,Y)}$. Thus, $y'\in\ind\V_{(\Q,Y)}\cap\B$ and then $y\cong \pi(y')\in\pi\left(\V_{(\Q,Y)}\cap\B\right)$. By Lemma \ref{lemma V intersected with B is equal to the co-aisle of U intersected with A} we also have the equality $\pi\left(\V_{(\Q,Y)}\cap\B\right) = \pi\left((\U\cap\A)^{^{\perp}}\right)$. It remains to show the equality $\Y = \pi\V_{(\Q,Y)}$, to do so we check that $\left(\X,\pi\V_{(\Q,Y)}\right)$ is a torsion pair.
	
By Lemma \ref{lemma pi V is included in X perp} we have that $\pi\V_{(\Q,Y)}\subseteq\X^{\perp}$, we show that $\X\ast\pi\V_{(\Q,Y)} = \ovl\C_m$. Since $\X = \pi(\U\cap\A)$ and $\pi\V_{(\Q,Y)} = \pi\left((\U\cap\A)^{\perp}\right)$, it is equivalent to show that $\pi\left(\U\cap\A\right) \ast\pi\left((\U\cap\A)^{\perp}\right)= \ovl\C_m$. 
	
Let $a\in\ovl\C_m$, there exists $a'\in\C_{2m}$ such that $\pi(a')\cong a$. Since $(\U\cap\A)\ast (\U\cap\A)^{^\perp} = \C_{2m}$, there exists a triangle $u\lora a'\lora v\lora \S a$ in $\C_{2m}$ with $u\in\U\cap\A$ and $v\in(\U\cap\A)^{^\perp}$. After localising we obtain the triangle $\pi(u)\lora a\lora \pi(v)\lora \S\pi(u)$ in $\ovl\C_m$. Note that $\pi(u)\in\pi(\U\cap\A)$ and $\pi(v)\in\pi\left((\U\cap\A)^{\perp}\right)$, thus we have that $a\in \pi\left(\U\cap\A\right) \ast\pi\left((\U\cap\A)^{\perp}\right)$. We can conclude that $\left(\X,\pi\V_{(\Q,Y)}\right)$ is a torsion pair, and as a consequence $\Y = \pi\V_{(\Q,Y)}$.
\end{proof}

\subsection{Hearts}\label{hearts t-structures}

\begin{figure}[ht]
\centering
\includegraphics[height = 5cm]{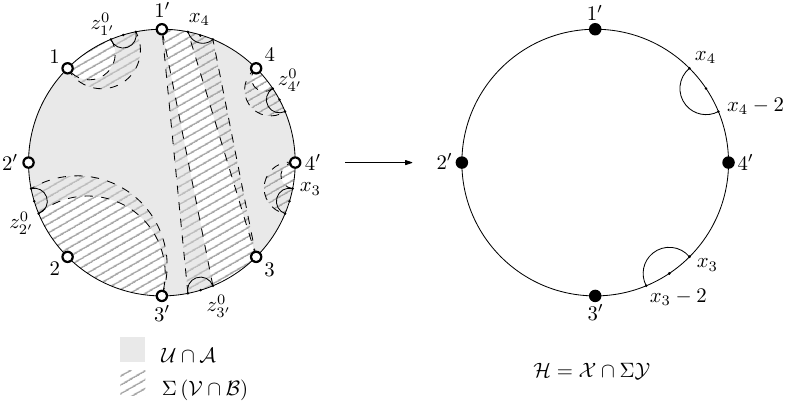}
\caption{The heart of the t-structure $(\X,\Y)$ of Figure \ref{figure half decorated non-crossing partition aisle t-structure} and Figure \ref{figure complement half-decorated non-crossing partition co-aisle t-structure}.}
\label{figure heart t-structure}
\end{figure}

With Theorem \ref{theorem t-structures} we classified the aisles of t-sturctures in $\ovl\C_m$, and with Proposition \ref{proposition compute the co-aisle of a t-structure} we computed the corresponding co-aisles. Now we can compute the heart of a t-structure $(\X,\Y)$ in $\ovl\C_m$. We first consider the preimage of $(\X,\Y)$ under $\pi$, which we denote by $(\U,\V)$. Note that $(\U,\V)$ is not a t-structure of $\C_{2m}$, but $(\U\cap\A, \V\cap\B)$ is. We can compute the heart of $(\U\cap\A,\V\cap\B)$ as in \cite[Corollary 4.15]{GZ}, and then obtain the heart of $(\X,\Y)$ by localising. Figure \ref{figure heart t-structure} illustrates this process.

\begin{corollary}\label{corollary heart t-structure}
Let $(\X,\Y)$ be a t-structure in $\ovl\C_m$. Consider its associated decorated non-crossing partition $(\P,X)$ of $[m']\cup[m]$ with $X = (x_p)_{p\in[m]}$. Then the heart $\mathcal{H} = \X\cap\S \Y$ is given by
\[
\mathcal{H} = \add \{ (x_p-2,x_p)\mid p\in[m] \text{ and } x_p\in\Zc_m\}.
\]
\end{corollary}
\begin{proof}
Let $\U = \pi^{-1}\X$, $\V = \pi^{-1}\Y$, $\U' = \U\cap\A$, and $\V' = \V\cap\B$. By Lemma \ref{lemma V intersected with B is equal to the co-aisle of U intersected with A}, the pair $\left(\U',\V'\right)$ is a t-structure in $\C_{2m}$. Consider the heart $\mathcal{H}' = \U'\cap\S\mathcal{V}'$, we show that $\pi\mathcal{H}' = \mathcal{H}$. Then the claim follows directly from \cite[Corollary 4.15]{GZ}.
	
First we show the inclusion $\pi\mathcal{H}'\subseteq\mathcal{H}$. Consider $h'\in\ind\mathcal{H}'$. Since $h'\in \U' \subseteq \U$, we have that $\pi h'\in \pi\U$ and from Proposition \ref{proposition bijection between extension closed subcategories in the verdier quotient} we have that $\pi\U = \X$. Similarly, since $h' \in \S\V' \subseteq\S\V$, we obtain that $\pi h' \in \pi\S\V = \S\pi\V = \S\Y$. Thus, $\pi h' \in \X\cap\S\Y = \mathcal{H}$. 
	
Now we show the inclusion $\mathcal{H}\subseteq\pi\mathcal{H}'$. Let $h\in\ind\mathcal{H}$, by Lemma \ref{lemma B co-aisle t-structure} there exists $h'\in\ind\A\cap\S\B$ such that $\pi h'\cong h$. Since $h\in\X$, $h'\in\pi^{-1}\X = \U$. Moreover, since $h\in\S\Y$, then $h'\in\pi^{-1}\S\Y = \S\V$. Thus, $h'\in\U\cap\A$ and $h'\in\S\V\cap\S\B$. We obtain that $h'\in\U'\cap\S\V' = \mathcal{H}'$, and then $h\cong\pi h'\in\pi\mathcal{H}'$. We can conclude that $\mathcal{H} = \pi\mathcal{H}'$. 
\end{proof}

\subsection{Boundedness}\label{section bounded t-structures}

The bounded t-structures in $\C_m$ were classified in \cite[Section 4.4]{GZ}. In $\C_m$ there exist bounded t-structures only if $m = 1$, see \cite[Remark 4.20, Corollary 4.22]{GZ}. Here we classify the bounded t-structures in $\ovl\C_m$, and we obtain that for each $m\geq 1$ there are no bounded t-structures in $\ovl\C_m$.

\begin{proposition}\label{proposition bounded above t-structures}
Let $(\X,\Y)$ be a t-structure in $\ovl\C_m$, let $(\P,X)$ be its associated half-decorated non-crossing partition of $[m']\cup[m]$, $\U = \pi^{-1}\X$ and $\V = \pi^{-1}\Y$. The following statements are equivalent.
\begin{enumerate}
\item The t-structure $(\X,\Y)$ is left bounded in $\ovl\C_m$.
\item The t-structure $(\U\cap\A,\V\cap\B)$ is left bounded in $\C_{2m}$.
\item The non-crossing partition $\P$ has as unique block $\{1',1,\dots,m',m\}$.
\end{enumerate}
\end{proposition}
\begin{proof}
We prove the equivalence of statements (1) and (2), for the equivalence between (2) and (3) we refer to \cite[Proposition 4.21]{GZ}. Assume that (1) holds, we check the inclusion $\C_{2m}\subseteq \bigcup_{n\in\Z}\S^n(\U\cap\A)$, the other inclusion is trivial. Consider $a\in\ind\C_{2m}$. Note that there exists $k\in\Z$ such that $a\in\S^k\A$. Moreover, since $\pi(a)\in\ovl\C_m$, there exists $l\in\Z$ such that $\pi(a)\in\S^l\X$, and  then $a\in\pi^{-1}\S^l\X = \S^l\U$. Thus, we have that $a\in\S^l\U\cap\S^k\A$. Let $n = \min\{k,l\}$, then $a\in\S^n(\U\cap\A)$.
	
Now we assume that (2) holds, we check the inclusion $\ovl\C_m\subseteq\bigcup_{n\in\Z}\S^n\X$, the other inclusion is trivial. Consider $a\in\ind\ovl\C_m$, then there exists $a'\in\ind\C_{2m}$ such that $\pi(a')\cong a$. Then there exists $n\in\Z$ such that $a'\in\S^n(\U\cap\A)\subseteq\S^n\U$, and then $a\cong\pi(a')\in\pi\S^n\U = \S^n\X$. This concludes the proof.
\end{proof}

Dually, we have the following proposition.

\begin{proposition}\label{proposition bounded below t-structures}
Keeping the assumptions and notation of Proposition \ref{proposition bounded above t-structures}, the following statements are equivalent.
\begin{enumerate}
\item The t-structure $(\X,\Y)$ is right bounded in $\ovl\C_m$.
\item The t-structure $(\U\cap\A,\V\cap\B)$ is right bounded in $\C_{2m}$.
\item The non-crossing partition $\P$ has as blocks $\{1'\}$, $\{1\}$, \dots, $\{m'\}$, $\{m\}$.
\end{enumerate}
\end{proposition}

We have the following corollary of Proposition \ref{proposition bounded above t-structures} and Proposition \ref{proposition bounded below t-structures}.

\begin{corollary}
For each $m\geq 1$ there are no bounded t-structures in $\ovl\C_m$.
\end{corollary}

\subsection{Non-degeneracy}

We classify the non-degenerate t-structures in $\ovl\C_m$. We refer to \cite[Corollary 4.19]{GZ} for the classification of the non-degenerate t-structures in $\C_m$.

\begin{proposition}\label{proposition left non-degenerate t-structures}
Let $(\X,\Y)$ be a t-structure in $\ovl\C_m$, let $(\P,X)$ be its associated half-decorated non-crossing partition with $X = (x_p)_{p\in[m]}$, and $\U = \pi^{-1}\X$. The following statements are equivalent.
\begin{enumerate}
\item The t-structure $(\X,\Y)$ is left non-degenerate in $\ovl\C_m$.
\item We have that $\bigcap_{n\in\Z}\S^n\U = \D$.
\item For each $p\in[m]$ we have that $x_p\neq p^+$, and for each $p,q\in [m']$ if $p,q\in B$ for some block $B\in\P$, then $p = q$.
\end{enumerate}
\end{proposition}
\begin{proof}
First we show the equivalence between the statements (1) and (2). Assume that $(\X,\Y)$ is left non-degenerate, i.e. $\bigcap_{n\in\Z}\S^n\X = 0$. The inclusion $\D\subseteq\bigcap_{n\in\Z}\S^n\U$ is straightforward, we show the other inclusion. Consider $u\in\ind\C_{2m}$ such that $u\in\S^n\U$ for all $n\in\Z$, then $\pi(u)\in\pi\S^n\U = \S^n\X$ for all $n\in\Z$. As a consequence, $\pi(u) \cong 0$ and then $u\in\D$. Now assume that $\bigcap_{n\in\Z}\S^n\U = \D$, we show that $\bigcap_{n\in\Z}\S^n\X = 0$. Assume that there exists $x\in\ind\ovl\C_m$ such that $x\in\S^n\X$ for all $n\in\Z$. Then there exists $x'\in\ind\C_{2m}$ such that $\pi(x')\cong x$, and $x'\in\pi^{-1}(\S^n\X) = \S^n\U$ for all $n\in\Z$. Then $x'\in\D$ and $x\cong\pi(x') \cong 0$, contradicting the fact that $x\in\ind\ovl\C_m$. This proves the equivalence between (1) and (2).
	
Now we prove the equivalence between statements (2) and (3). Assume that $\bigcap_{n\in\Z}\S^n\U = \D$ and that there exists $p\in[m]$ such that $x_p = p^+$, then $\Z^{(p,p)}\subseteq\ind\U$. Let $u\in\Z^{(p,p)}$, then $x\in\S^n\U$ for each $n\in\Z$. As a consequence $x\in\D$, and this contradicts the fact that $p\in[m]$. This proves that $x_p\neq p^+$. Now consider $p,q\in[m']$ such that $p,q\in B$ for some block $B\in\P$, then $\U$ contains all arcs having one endpoint in $\Z^{(p)}$ and the other in $\Z^{(q)}$. Consider such $u$, then $u\in\S^n\U$ for each $n\in\Z$. As a consequence, $u\in\D$ and then $p = q$. This proves that (2) implies (3).
	
Now assume that statement (3) holds, we show that $\bigcap_{n\in\Z}\S^n\U = \D$. The inclusion $\D\subseteq\bigcap_{n\in\Z}\S^n\U$ is straightforward, we prove the other inclusion. Let $u\in\ind\bigcap_{n\in\Z}\S^n\U$, we show that $u\in\D$. Assume that $u$ has an endpoint $z\in\Z^{(p)}$ for some $p\in[m]$. Since $u\in\ind\U$, then $z\in(p,x_p]$. Moreover, since $x_p\neq p^+$, there exists $n\in\Z$ such that $\S^n u\notin\U$, and this contradicts the fact that $u\in\bigcap_{n\in\Z}\S^n\U$. Thus, $u\in\Z^{(p,q)}$ for some $p,q\in[m']$. Then $p,q\in B$ for some $B\in\P$ and as a consequence $p = q$, i.e. $u\in\D$. This concludes the argument.
\end{proof}

Dually, we have the following proposition.

\begin{proposition}\label{proposition right non-degenerate t-structures}
Keeping the assumptions and notation of Proposition \ref{proposition left non-degenerate t-structures}, let $(\Q,Y) = (\P,X)^c$ and $\V = \pi^{-1}\Y$. The following statements are equivalent.
\begin{enumerate}
\item The t-structure $(\X,\Y)$ is right non-degenerate.
\item We have that $\bigcap_{n\in\Z}\S^n\V = \D$.
\item For each $p\in[m]$ we have that $x_p\neq p$, and for each $p,q\in [m']$ if $p,q\in C$ for some block $C\in\Q$, then $p = q$.
\end{enumerate}
\end{proposition}

Combining Proposition \ref{proposition left non-degenerate t-structures} and Proposition \ref{proposition right non-degenerate t-structures} we obtain the following corollary.

\begin{corollary}\label{corollary non-degenerate t-structures}
Keeping the assumptions and notation of Proposition \ref{proposition left non-degenerate t-structures} and Proposition \ref{proposition right non-degenerate t-structures}, the following statements are equivalent.
\begin{enumerate}
\item The t-structure $(\X,\Y)$ is non-degenerate.
\item We have that $\bigcap_{n\in\Z}\S^n\U = \D = \bigcap_{n\in\Z}\S^n\V$.
\item For each $p\in[m]$ we have that $x_p\in\Z^{(p)}$, and for each $p,q\in [m']$ if $p,q\in B$ for some block $B\in\P$, or $p,q\in C$ for some block $C\in\Q$, then $p = q$.
\end{enumerate}
\end{corollary}

With the following example we show that there exist half-decorated non-crossing partitions of $[m']\cup[m]$ satisfying condition (3) of Corollary \ref{corollary non-degenerate t-structures}.

\begin{example}\label{example non-degenerate t-structure}
Let $\P$ be the non-crossing partition $\{\{1',1\}, \{2',2\},\dots,\{m',m\}\}$ of $[m']\cup[m]$, and let $X = (x_p)_{p\in[m]}$ with $x_p\in\Z^{(p)}$ for each $p\in[m]$. Then $(\P,X)$ is a half-decorated non-crossing partition of $[m']\cup[m]$ and $\P^c = \{\{1'\},\{2'\},\dots,\{m'\},\{1,2,\dots,m\}\}$, see Section \ref{section non-crossing partitions}. Note that $(\P,X)$ satisfies condition (3) of Corollary \ref{corollary non-degenerate t-structures}.
\end{example}

As a consequence, we have the following corollary.

\begin{corollary}\label{corollary existence non-degenerate t-structures}
Non-degenerate t-structures in $\ovl\C_m$ exist for each $m\geq 1$.
\end{corollary}

\section{Co-t-structures}\label{section co-t-structures}

\begin{figure}[ht]
\centering
\includegraphics[height = 4cm]{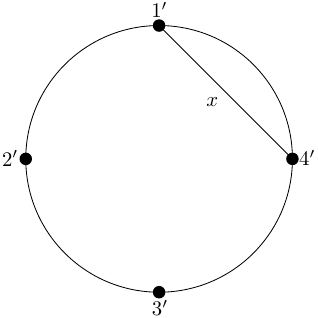}
\caption{The subcategory $\add\{x\}$ of $\ovl\C_4$ is the aisle of a co-t-structure.}
\label{figure example co-t-structure in completion}
\end{figure}

From \cite[Proposition 4.6]{ZZ} we know that in the category $\C_{m}$ the only co-t-structures are $(\C_m,0)$ and $(0,\C_m)$. In $\ovl\C_m$ this is not the case, Figure \ref{figure example co-t-structure in completion} gives an example of a non-trivial co-t-structure in $\ovl\C_m$. In this section we classify the aisles of the co-t-structures, we compute the co-aisles and co-hearts. We also classify the bounded and non-degenerate co-t-structures, and the co-t-structures having a left or right adjacent t-structure. Moreover, from the classification of the co-t-structures, we can easily obtain the classification of the recollements of $\ovl\C_m$.

\subsection{Aisles of co-t-structures}

Here we classify the aisles of the co-t-structures in $\ovl\C_m$ in a way similar to the classification of the aisles of the t-structures in Section \ref{section aisles of t-structures}. The definition below is the analogue of Definition \ref{definition half-decorated non-crossing partition}.

\begin{definition}\label{definition alternating non-crossing partition}
An \emph{alternating non-crossing partition} of $[m']\cup[m]$ is a pair $(\P,X)$ given by a non-crossing partition $\P$ of $[m']$ and an $m$-tuple $X = (x_p)_{p\in[m]}$ such that $x_p\in [p,p^+]$ for each $p\in[m]$.
\end{definition}

The main result of this section is the following analogue of Theorem \ref{theorem t-structures}. The notation employed in the statement will be defined in Definition \ref{definition almost aisle of co-t-structure} and Definition \ref{definition from almost aisle of co-t-structure to h-decorated h-nc partition}.

\begin{theorem}\label{theorem classification aisles co-t-structures}
The following is a bijection. 
\begin{align*}
\left\{\parbox{4.5cm}{\centering Alternating non-crossing partitions of $[m']\cup[m]$}\right\} & \longleftrightarrow \left\{\parbox{5.2cm}{\centering Aisles of co-t-structures in $\ovl\C_m$} \right\}\\
(\P,X) & \longmapsto \pi\U_{(\P,X)}\\
(\P_{\pi^{-1}\X},X_{\pi^{-1}\X}) & \longmapsfrom \X
\end{align*}
\end{theorem}

\begin{figure}[ht]
\centering
\includegraphics[height = 5cm]{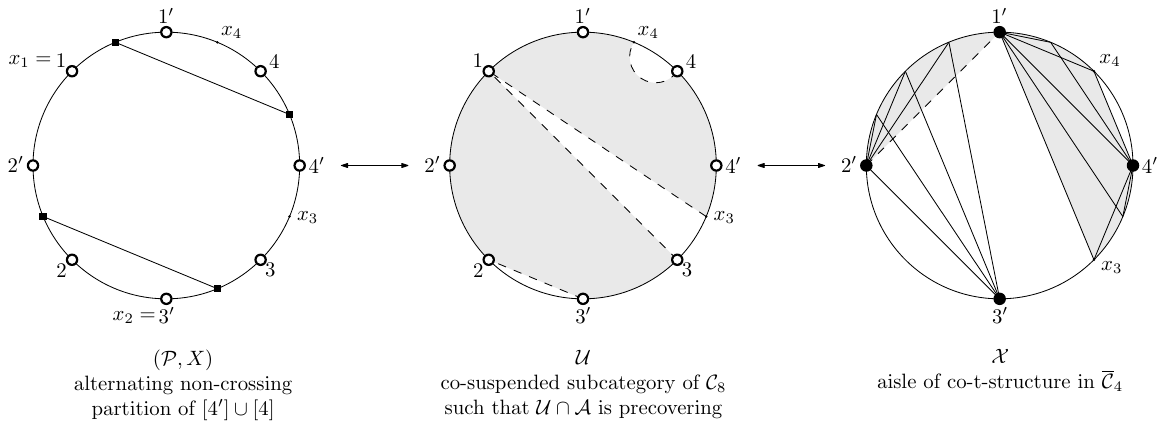}
\caption{Illustration of how to obtain the aisle of a co-t-structure in $\ovl\C_m$ from an alternating non-crossing partition of $[m']\cup[m]$.}
\label{figure alternating non-crossing partition aisle co-t-structure}
\end{figure}

To prove this result, we proceed as in Section \ref{section aisles of t-structures} by taking an intermediate step through $\C_{2m}$. From Corollary \ref{corollary co-t-structures} the aisles of co-t-structures in $\ovl\C_m$ are in bijection with certain subcategories of $\C_{2m}$, which can be regared as ``almost aisles" of co-t-structures. These are co-suspended subcategories $\U$ of $\C_{2m}$ such that $\D\subseteq\U$ and $\U\cap\A$ is precovering, and are classified in Proposition \ref{proposition classification co-t-structures} in terms of alternating non-crossing partitions of $[m']\cup[m]$. The aisles of co-t-structures in $\ovl\C_m$ are then obtained after localising the ``almost aisles" of co-t-sturctures in $\C_{2m}$. Figure \ref{figure alternating non-crossing partition aisle co-t-structure} illustrates this process.

The rest of this section is devoted to prove the following proposition. The assignments of the maps $\a$ and $\b$ will be defined in Definition \ref{definition almost aisle of co-t-structure} and Definition \ref{definition from almost aisle of co-t-structure to h-decorated h-nc partition}.

\begin{proposition}\label{proposition classification co-t-structures}
The following is a bijection.
\begin{align*}
\left\{\parbox{6cm}{\centering  Alternating non-crossing partitions of $[m']\cup[m]$}\right\} & \longleftrightarrow \left\{\parbox{7.5cm}{\centering Co-suspended subcategories $\U\subseteq\D$ such that $\D\subseteq\U$ and $\U\cap\A$ is precovering} \right\}\\
\a\colon (\P,X) & \longmapsto \U_{(\P,X)}\\		
(\P_{\U},X_{\U}) & \longmapsfrom \U\colon \b
\end{align*}
\end{proposition}

The following definition and lemma define the map $\a$ and show that it is well defined.

\begin{definition}\label{definition almost aisle of co-t-structure}
Let $(\P,X)$ be an alternating non-crossing partition of $[m']\cup[m]$. We define 
\[
\U_{(\P,X)} = \add\bigsqcup_{B\in\P} \left\{(u_1,u_2)\in\ind\C_{2m}\middle| u_1,u_2\in\bigcup_{p\in B}[x_{p^-},p^+)\right\},
\]
where  we use the following convention: if $x_{p^-} = p^-$ then $[x_{p^-},p^+) = \Z^{(p^-)}\sqcup\Z^{(p)}$, and if $x_{p^-} = p$ then $[x_{p^-},p^+) = \Z^{(p)}$.
\end{definition}

\begin{proposition}\label{proposition from alternating non-crossing partition to almost aisle of co-t-structure}
Let $(\P,X)$ be an alternating non-crossing partition of $[m']\cup[m]$. Then $\U_{(\P,X)}$ is co-suspended subcategory of $\C_{2m}$ such that $\D\subseteq\U$ and $\U\cap\A$ is precovering.
\end{proposition}
\begin{proof}
In order to show that $\U_{(\P,X)}$ is extension-closed, contains $\D$, and is closed under $\S^{-1}$, we can proceed similarly to the argument of Proposition \ref{proposition from half decorated non-crossing partition to almost aisle of t-structure}. We show that $\U_{(\P,X)}\cap\A$ is precovering. From \cite[Theorem 3.1]{GHJ} we know that this is equivalent to show that $\U_{(\P,X)}$ satisfies the $\PC$ conditions, see Definition \ref{definition pc conditions}. 

We check that $\U_{(\P,X)}$ satisfies (PC1), the other conditions are analogous. Assume that there exists a sequence $\{(x_1^n,x_2^n)\}_n\subseteq\U_{(\P,X)}\cap\A\cap\Z^{(p,q)}$ for some $p,q\in[m']\cup[m]$ such that $p\neq q$ and the sequences $\{x_1^n\}_n$ and $\{x_2^n\}_n$ are strictly increasing. Then $p,q\in[m]$. By Definition \ref{definition almost aisle of co-t-structure} we have that $p^+, q^+\in B$ for some block $B\in\P$. As a consequence, there exist strictly decreasing sequences $\{y_1^n\}_n\subseteq\Z^{(p^+)}$ and $\{y_2^n\}_n\subseteq\Z^{(q^+)}$ such that $\{|y_1^n,y_2^n|\}_n\subseteq\U_{(\P,X)}\cap\A$. This proves that (PC1) holds and concludes the argument.
\end{proof}

With the following definition and proposition we define the map $\b$ of Proposition \ref{proposition classification co-t-structures} and we check that it is well defined. Given a co-suspended subcategory $\U$ of $\C_{2m}$ such that $\D\subseteq\U$ and $\U\cap\A$ is precovering, we define the equivalence relation $\sim_{\U}$ on the set $[m']$ as in Definition \ref{definition equivalence relation}.

\begin{definition}\label{definition from almost aisle of co-t-structure to h-decorated h-nc partition}
Let $\U$ be a co-suspended subcategory of $\C_{2m}$ such that $\D\subseteq\U$ and $\U\cap\A$ is precovering. We define $\P_{\U}$ to be the partition of $[m']$ given by the equivalence classes of $\sim_{\U}$. For each $p\in[m]$ we define 
\[
x_p = \inf\{z\in\Z^{(p)}\mid\text{there exists $u\in\U$ with an endpoint equal to $z$}\}.
\]
We denote by $X_{\U}$ the $m$-tuple $X_{\U} = (x_p)_{p\in[m]}$.
\end{definition}

\begin{proposition}
Keeping the notation of Definition \ref{definition from almost aisle of co-t-structure to h-decorated h-nc partition}, the pair $(\P_{\U},X_{\U})$ is an alternating non-crossing partition of $[m']\cup[m]$.
\end{proposition}
\begin{proof}
We already know that $\P_{\U}$ is a partition of $[m']$, we only need to check that $\P_{\U}$ non-crossing. To this end, we can apply the same argument of \cite[Lemma 4.12]{GZ}.
\end{proof}

The following lemma is useful for the argument of Proposition \ref{proposition classification co-t-structures}.

\begin{lemma}\label{lemma intervals contained in U}
Let $\U$ be a co-suspended subcategory of $\C_{2m}$ such that $\D\subseteq \U$ and $\U\cap\A$ is precovering. Consider the alternating non-crossing partition $(\P_{\U},X_{\U})$ with $X_{\U} = (x_{p})_{p\in[m]}$. Let $p,q\in[m']$ be such that $p,q\in B$ for some block $B\in\P_{\U}$. Then any arc of $\C_{2m}$ having one endpoint in $[x_{p^-},p^+)$ and the other in $[x_{q^-},q^+)$ is an arc of $\U$.
\end{lemma}
\begin{proof}
In order to simplify the notation we assume that $q\neq 1'$, if $q = 1'$ we can proceed analogously. We denote by $[x_{p^-},p^+)\times [x_{q^-},q^+)$ the set of arcs $a = (a_1,a_2)\in\ind\C_{2m}$ such that $a_1\in[x_{p^-},p^+)$ and $a_2\in[x_{q^-},q^+)$. We show that $[x_{p^-},p^+)\times [x_{q^-},q^+)\subseteq \ind\U$. We have the equality
\[
[x_{p^-},p^+)\times [x_{q^-},q^+) = \Z^{(p,q)}\sqcup \left([x_{p^-},p)\times \Z^{(q)}\right)\sqcup \left(\Z^{(p)}\times [x_{q^-},q)\right)\sqcup \left([x_{p^-},p)\times [x_{q^-},q)\right).
\]
We assume that $x_{p^-}\neq p$ and $x_{q^-}\neq q$, the other cases are analogous. We divide the proof into steps.

\emph{Step 1.} We show that $\Z^{(p,q)}\subseteq\ind\U$.

If $p = q$ we have the claim from the fact that $\D\subseteq\U$. Now assume that $p\neq q$. Since $p,q\in B$ for some block $B\in\P_{\U}$, there exists $u\in\ind\U$ such that $u\in\Z^{(p,q)}$. We show that $\S^n u\in\ind\U$ for each $n\in\Z$, then, using the fact that $\U$ is extension-closed, it is straightforward to check that $\U$ contains any arc of $\Z^{(p,q)}$. Since $\S^{-1}\U\subseteq\U$, we already know that $\S^n u\in\U$ for each $n\leq 0$, it remains to check that $\S^n u\in\U$ for each $n\geq 1$. Consider the arcs $a = (u_1-1,u_1+1)\in\Z^{(p,p)}$ and $b = (u_2-1,u_2+1)\in\Z^{(q,q)}$. The arcs $u$ and $a$ are crossing, and then, since $\U$ satisfies the $\PT$ condition, $u' = (u_1-1,u_2)\in\U$. Moreover, the arcs $u'$ and $b$ are crossing and $u'' = (u_1-1,u_2-1) = \S u\in\U$. Repeating this argument we obtain that $\S^n u\in\U$ for each $n\geq 1$. This concludes the argument of Step 1.

\emph{Step 2.} We show that $[x_{p^-},p)\times \Z^{(q)}\subseteq\ind\U$.

Let $a = (a_1,a_2)\in\ind\C_{2m}$ with $a_1\in[x_{p^-},p)$ and $a_2\in\Z^{(q)}$, we show that $a\in\U$. First we show that there exists an arc of $\U$ with an endpoint equal to $a_1$. If there is not such arc, then there is no arc $u\in\ind\U$ with an endpoint in $[x_{p^-},a_1]$, otherwise $\S^{n} u\in\U$ has an endpoint equal to $a_1$ for some $n\leq 0$. Since 
\[
x_{p^-} = \inf\{z\in\Z^{(p^-)}\mid\text{there exists an arc of $\U$ with an endpoint equal to $z$}\}
\] 
this gives a contradiction, and therefore there exists an arc of $\U$ with an endpoint equal to $a_1$. Let $u'$ be such arc, then $\S^{n} u'\in\U$ for each $n\leq 0$. Moreover, since $\U$ satisfies the $\PT$ condition, we obtain that $(a_1,a_1+2), (a_1,a_1+3), \dots\in\U$. Note that these arcs are also in $\A$ because they belong to $\Z^{(p^-)}$. Therefore, we have a sequence $\{(a_1,a_1+2+n)\}_{n\geq 0}\subseteq\ind\U\cap\A$ such that $\{a_2+2+n\}_{n\geq 0}$ is strictly increasing. Since $\U\cap\A$ is precovering and satisfies condition $(\PC 3)$, it follows that there exists an arc $v = (a_1,v_2)\in\Z^{(p^-,p)}\cap\U$ such that $v_2\leq a_2$. Consider an arc of the form $z = (z_1,a_2)\in\Z^{(p,q)}$ with $p<z_1<v_2$. Since the arcs $v$ and $z$ cross and $\Z^{(p,q)}\subseteq\U$, then $a = (a_1,a_2)\in\U$. This concludes the argument of Step 2.

\emph{Step 3.} Analogously as in Step 2 we have that $\Z^{(p)}\times[x_{q^-},q)\subseteq\ind\U$.

\emph{Step 4.} We show that $[x_{p^-},p)\times [x_{q^-},q)\subseteq\ind\U$.

Let $a = (a_1,a_2)\in[x_{p^-},p)\times [x_{q^-},q)$, we show that $a\in\U$. If $p = q$, consider the sequence $\{(a_1,a_1+2+n)\}_{n\geq 0}\subseteq\U\cap\A$ of Step 2. Since $(a_1,a_2) = (a_1,a_1+2+n)\in\U$ for some $n\geq 0$, we have that $a\in\U$. Now assume that $p\neq q$. We consider the arc $v = (a_1,v_2)\in\Z^{(p^-,p)}$ of Step 2, and an arc $z = (z_1,a_2)\in\Z^{(p,q^-)}$ with $p<z_1<v_2$.  Since the arcs $v$ and $z$ cross and $z\in\Z^{(p)}\times[x_{q^-},q)\subseteq\U$, by Step 3 $a = (a_1,a_2)\in\U$. This concludes the argument of Step 4.

We can conclude that $[x_{p^-},p^+)\times [x_{q^-},q^+)\subseteq\ind\U$.
\end{proof}

Finally, we can prove Proposition \ref{proposition classification co-t-structures}.

\begin{proof}[Proof of Proposition \ref{proposition classification co-t-structures}]
From Definition \ref{definition almost aisle of co-t-structure} and Proposition \ref{proposition from alternating non-crossing partition to almost aisle of co-t-structure} we have that the maps are well defined, we prove that they are mutually inverse. We divide the proof into steps.

\emph{Step 1.} The map $\a$ is injective.

Let $(\P,X)$ and $(\Q,Y)$ be two alternating non-crossing partitions of $[m']\cup[m]$ such that $\U_{(\P,X)} = \U_{(\Q,Y)}$, we show that $(\P,X) = (\Q,Y)$. Assume that $\P\neq \Q$. Then there exist $p,q\in[m']$ with $p\neq q$ such that
$p$ and $q$ belong to the same block of $\P$ and to distinct blocks of $\Q$, or vice versa $p$ and $q$ belong to the same block of $\Q$ and to distinct blocks of $\P$. In the first case, there exists an arc of $\U_{(\P,X)}$ with an endpoint in $\Z^{(p)}$ and the other in $\Z^{(q)}$, while there is no such arc in $\U_{(\Q,Y)}$. As a consequence $\U_{(\P,X)}\neq \U_{(\Q,Y)}$, giving a contradiction. In the second case the role of $\P$ and $\Q$ exchange and we obtain the same contradiction. Thus we have that $\P = \Q$.

Now we show that $X = Y$. We denote $X = (x_p)_{p\in[m]}$ and $Y = (y_p)_{p\in[m]}$, and we assume that $X\neq Y$. Let $p\in[m]$ be such that $x_p\neq y_p$, then either $x_p<y_p$ or $x_q<y_p$. Assume that $x_p<y_p$, the other case is analogous. Since $p\leq x_p<y_p$, there is an arc of $\U_{(\Q,Y)}$ with an endpoint greater that $x_p$, while there is no such arc in $\U_{(\P,X)}$. We obtain that $\U_{(\P,X)}\neq \U_{(\Q,Y)}$, giving a contradiction. This concludes the argument of Step 1.

\emph{Step 2.} We show that $\a\b = \id$.

Consider $\U$ a co-suspended subcategory of $\C_{2m}$ such that $\D\subseteq\U$ and $\U\cap\A$ is precovering. We show that $\U = \U_{(\P_{\U},X_{\U})}$. First we show the inclusion $\U_{(\P_{\U},X_{\U})}\subseteq \U$. Consider $u = (u_1,u_2)\in\ind\U_{(\P_{\U},X_{\U})}$, then there exist a block $B\in\P_{\U}$ and $p,q\in B$ such that $u_1\in[x_{p^-},p^+)$ and $u_2\in[x_{q^-},q^+)$, where $X_{\U} = (x_p)_{p\in[m]}$. Then $u\in\ind\U$ by Lemma \ref{lemma intervals contained in U}.

Now we show the inclusion $\U\subseteq\U_{(\P_{\U}, X_{\U})}$. Consider $u = (u_1,u_2)\in\ind\U$, then there exist $p,q\in[m']$ such that $u_1\in[x_{p^-},p^+)$ and $u_2\in[x_{q^-},q^+)$ where 
\[x_{p^-} = \inf\{z\in\Z^{(p^-)}\mid \text{there exists an arc of $\U$ with an endpoint equal to $z$}\}\]
and $x_{q^-}$ is defined similarly. We show that $p,q\in B$ for some block $B\in\P_{\U}$, i.e. that there exists an arc of $\U$ with an endpoint in $\Z^{(p)}$ and the other in $\Z^{(q)}$. Then we can conclude that $u_1,u_2\in \bigcup_{p\in B} [x_{p^-},p^+)$, and then $u\in\U_{(\P_{\U},X_{\U})}$. If $p = q$ the claim is straightforward, we assume that $p\neq q$. We can write $[x_{p^-},p^+) = [x_{p^-},p)\sqcup\Z^{(p)}$ and $[x_{q^-},q^+) = [x_{q^-},q)\sqcup\Z^{(q)}$. If $u_1\in\Z^{(p)}$ and $u_2\in\Z^{(q)}$ then the claim follows directly. We assume that $u_1\notin\Z^{(p)}$ or $u_2\notin\Z^{(q)}$. 

Assume that $u_1\in[x_{p^-},p)$ and $u_2\in[x_{q^-},q)$. We consider the sequence $\{\S^{-n}u = (u_1+n,u_2+n)\}_{n\geq 0}\subseteq\ind\U$. This sequence is also in $\A$ because it is contained in $\Z^{(p^-,q^-)}$. Since $\U\cap\A$ is precovering and satisfies condition $(\PC 1)$, there exists an arc of $\U$ with an endpoint in $\Z^{(p)}$ and the other in $\Z^{(q)}$. This gives the claim.

Now assume that $u_1\in\Z^{(p)}$ and $u_2\in[x_{q^-},q)$. We consider the sequence $\{\S^{-n}u = (u_1+n,u_2+n)\}_{n\geq 0}\subseteq\ind\U$. The sequence $\{(u_2,u_2+2+n)\}_{n\geq 0}\subseteq\Z^{(q^-,q^-)}$ is obtained from the crossings of the sequence $\{\S^{-n}u\}_{n\geq 0}$. We have that $\{(u_2,u_2+2+n)\}_{n\geq 0}\subseteq\ind\U\cap\A$ because this sequence is contained in $\Z^{(q^-,q^-)}$ and $\U$ satisfies the $\PT$ condition. Since $\U\cap\A$ is precovering and satisfies condition $(\PC 3)$, there exists an arc $x\in\ind\U$ with an endpoint equal to $u_2$ and the other in $\Z^{(q)}$. The arcs $\S^{-1}u$ and $x$ cross, and from this crossing we obtain an arc $u'\in\ind\U$ with an endpoint in $\Z^{(p)}$ and the other endpoint in $\Z^{(q)}$. The case where $u_1\in[x_{p^-},p)$ and $u_2\in\Z^{(q)}$ is analogous, therefore we have the claim. This concludes the argument of Step 2.

We can conclude that the two maps of the claim are mutually inverse.
\end{proof}

\subsection{Co-aisles of co-t-structures}\label{section co-aisles of co-t-structures}

\begin{figure}[ht]
\centering
\includegraphics[height = 5cm]{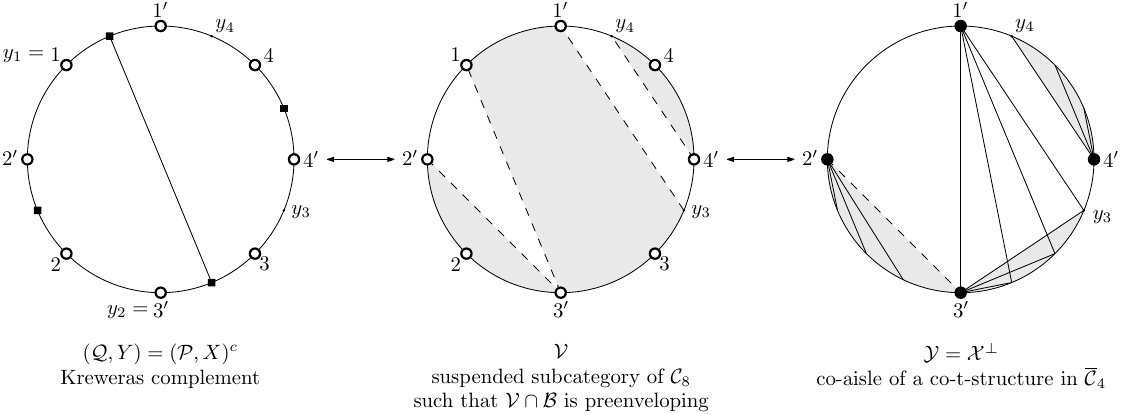}
\caption{Illustration of how to obtain the co-aisle of the aisle of Figure \ref{figure alternating non-crossing partition aisle co-t-structure}}
\label{figure complement alternating non-crossing partition co-aisle co-t-structure}
\end{figure}

We compute the co-aisles of co-t-structures in $\ovl\C_m$ using a method similar to Section \ref{section co-aisles t-structures}. From an alternating non-crossing partition $(\P,X)$ of $[m']\cup[m]$, we consider its complement $(\P,X)^c$, obtained from the Kreweras complement $\P^c$ of $\P$, see Section \ref{section non-crossing partitions}. This corresponds to a subcategory $\V$ of $\C_{2m}$, which can be thought as an ``almost co-aisle" of a co-t-structure in $\C_{2m}$. This is a suspended subcategory $\V$ of $\C_{2m}$ such that $\D\subseteq\V$ and $\V\cap\B$ is preenveloping. The subcategory $\V$ gives a co-aisle of a co-t-structure in $\ovl\C_m$ after localising. Figure \ref{figure complement alternating non-crossing partition co-aisle co-t-structure} illustrates this process.  

\begin{definition}\label{definition complement co-aisle co-t-structure}
Let $(\P,X)$ be an alternating non-crossing partition of $[m']\cup[m]$ with $X = (x_p)_{p\in[m]}$. We define the \emph{complement} $(\P,X)^c$ of $(\P,X)$ to be the pair $(\Q,Y)$ where $\Q = \P^c$ is the Kreweras complement of $\P$, and $Y$ is the $m$-tuple $Y = X-1 = (x_p-1)_{p\in[m]}$ for each $p\in[m]$.
\end{definition}

From the complement of an alternating non-crossing partition of $[m']\cup[m]$ we obtain an ``almost co-aisle" of co-t-structure in $\C_{2m}$. The following definition is similar to Definition \ref{definition almost co-aisle of t-structure}.

\begin{definition}\label{definition almost co-aisle of co-t-structure}
Let $(\P,X)$ be an alternating non-crossing partition of $[m']\cup[m]$ and let $(\Q, Y) = (\P,X)^c$. We define
\[
\V_{(\Q,Y)} = \add\bigsqcup_{B\in\Q}\left\{(v_1,v_2)\in\ind\C_{2m}\middle| v_1,v_2\in\bigcup_{p\in B} (p,y_{p^+}]\right\}.
\]
\end{definition}

The proof of the following lemma is analogous to the proof of Proposition \ref{proposition from half decorated non-crossing partition to almost aisle of t-structure}.

\begin{lemma}
Let $(\P,X)$ be an alternating non-crossing partition of $[m']\cup[m]$ and let $(\Q,Y) = (\P,X)^c$. Then $\V_{(\Q,Y)}$ is suspended and contains $\D$.
\end{lemma}

Consider the complement $(\Q,Y)$ of an alternating non-crossing partition of $[m']\cup[m]$. The following lemmas and remark describe some properties of the subcategory $\V_{(\Q,Y)}$.

\begin{lemma}\label{lemma V intersected with B is equal to the co-aisle of U intersected with A for co-t-structures}
Let $(\P,X)$ be an alternating non-crossing partition of $[m']\cup[m]$ and let $(\Q,Y) = (\P,X)^c$. Then $\V_{(\Q,Y)}\cap\B = (\U_{(\P,X)}\cap\A)^{{\perp}}$.
\end{lemma}
\begin{proof}
It is straightforward to check that
\[
\V_{(\Q,Y)}\cap \B= \add\bigsqcup_{B\in\Q}\left\{(v_1,v_2)\in\ind\C_{2m}\middle| v_1,v_2\in\bigcup_{p\in B}[w^0_p,y_{p^+}]\right\}
\]
where we recall from Section \ref{section precovering subcategories of the completion} that $w_p^0 = z_p^0-1$ for each $p\in[m']$. We denote the right hand side of the equality by $\mathcal{W}$. Proceeding analogously as in the argument of \cite[Corollary 4.14]{GZ}, it is straightforward to check that $\S^{-1}\mathcal{W}$ consists precisely of all the arcs of $\C_{2m}$ which do not cross $\U\cap\A$. As a consequence $\V_{(\Q,Y)}\cap\B = (\U\cap\A)^{^\perp}$.
\end{proof}

\begin{remark}\label{remark U intersected with A and V intersected with B torsion pair}
Let $\U_{(\P,X)}$ and $\V_{(\Q,Y)}$ be as in Lemma \ref{lemma V intersected with B is equal to the co-aisle of U intersected with A for co-t-structures}. Since $\U_{(\P,X)}\cap\A$ is precovering and extension-closed, by Proposition \ref{proposition iyama-yoshino torsion pairs} $\left(\U_{(\P,X)}\cap\A, \V_{(\Q,Y)}\cap\B\right)$ is a torsion pair. It is not a t-structure nor a co-t-structure because in general $\U_{(\P,X)}\cap\A$ is not closed under $\S$ or $\S^{-1}$, cf. Remark \ref{remark U intersected with A and V intersected with B t-structure}.
\end{remark}

Let $(\Q,Y)$ be the complement of an alternating non-crossing partition of $[m']\cup[m]$. With the following proposition we prove that by localising $\V_{(\Q,Y)}$ we obtain the co-aisle of a co-t-structure in $\ovl\C_m$. The argument is the same of Proposition \ref{proposition compute the co-aisle of a t-structure}.

\begin{proposition}\label{proposition compute the co-aisle of a co-t-structure}
Let $(\X,\Y)$ be a co-t-structure in $\ovl\C_m$, $\U = \pi^{-1}\X$, $(\P,X)$ be its associated alternating non-crossing partition, and $(\Q,Y) = (\P,X)^c$. Then
\[
\Y = \pi\V_{(\Q,Y)} = \pi\left(\V_{(\Q,Y)}\cap\B\right) = \pi\left((\U\cap\A)^{\perp}\right).
\]
\end{proposition}

\subsection{Co-hearts}

We classified the aisles of co-t-structures in $\ovl\C_m$ in Theorem \ref{theorem classification aisles co-t-structures}, and we computed the co-aisle of a co-t-structure in Proposition \ref{proposition compute the co-aisle of a co-t-structure}. Here we compute the co-heart of a co-t-structure in $\ovl\C_m$.

First we introduce some notation. Let $(\P,X)$ be an alternating non-crossing partition of $[m']\cup[m]$. Consider $p,q\in[m']\cup[m]$, we write $q = p^{+_B}$ if 
\begin{itemize}
\item $p,q\in[m']$, and
\item $p,q\in B$ for some block $B\in\P$, and 
\item $q$ is the next element of $[m']\cap B$ we meet while moving from $p$ along $S^1$ in the anticlockwise direction. 
\end{itemize}
If $B = \{p\}$, then by convention $p^{+_B} = p$.

\begin{figure}[ht]
\centering
\includegraphics[height = 5cm]{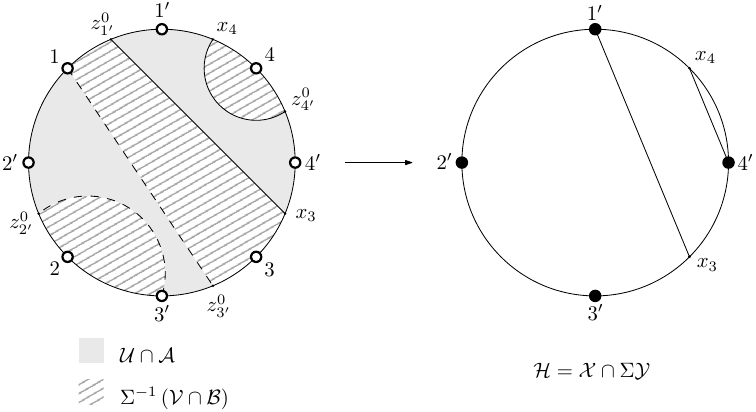}
\caption{Illustration of how to obtain the co-heart of the co-t-structure of Figure \ref{figure alternating non-crossing partition aisle co-t-structure} and Figure \ref{figure complement alternating non-crossing partition co-aisle co-t-structure}.}
\label{figure co-heart}
\end{figure}

Now let $(\X,\Y)$ be a co-t-structure in $\ovl\C_m$. We consider the preimage $(\pi^{-1}\X,\pi^{-1}\Y)$ of $(\X,\Y)$, which we denote by $(\U,\V)$. The pair $(\U,\V)$ is not a torsion pair, but $(\U\cap\A,\V\cap\B)$ is. Moreover, $(\U\cap\A,\V\cap\B)$ is not a co-t-structure, but we can still compute $\mathcal{S}' = (\U\cap\A)\cap\S^{-1}(\V\cap\B)$ similarly. The co-heart of $(\X,\Y)$ is obtained by localising $\mathcal{S}'$. Figure \ref{figure co-heart} illustrates this process.

In the proposition below we recall that $|z_{p}^0,x_{q^-}|$ is equal to $(z_p^0,x_{q^-})$ if $p<q^-$ and is equal to $(x_{q^-},z_p^0)$ if $q^-<p$, see Section \ref{section infinity-gons}.

\begin{proposition}\label{proposition almost co-heart of co-t-structure}
Let $(\X,\Y)$ be a co-t-structure in $\ovl\C_m$, $(\P,X)$ be its associated alternating non-crossing partition of $[m']\cup[m]$ with $X = (x_p)_{p\in[m]}$, $\U = \pi^{-1}\X$, and $\V = \pi^{-1}\Y$. Then
\[(\U\cap\A)\cap\S^{-1}(\V\cap\B) = \add\left\{|z_p^0,x_{q^-}|\middle|\text{$p,q\in[m']$, $q = p^{+_B}$ for some $B\in\P$, and $x_{q^-}\in\Z^{(q^-)}$}\right\}.\]
\end{proposition}
\begin{proof}
First we show that arcs of the form $a = |z_p^0,x_{q^-}|$, where $q = p^{+_B}$ for some block $B\in\P$ and $x_{q^-}\in\Z^{(q^-)}$, belong to $(\U\cap\A)\cap\S^{-1}(\V\cap\B)$. From Definition \ref{definition A} and Definition \ref{definition almost aisle of co-t-structure} we have that $a\in\ind\U\cap\A$, we check that $a\in\ind\S^{-1}(\V\cap\B)$. From Lemma \ref{lemma V intersected with B is equal to the co-aisle of U intersected with A for co-t-structures} this is equivalent to check that $a$ does not cross any arc $u\in\ind\U\cap\A$. Note that $z_p^0\in[z_p^0,x_{p^+}]$ and $x_{q^-}\in[z_{q^{--}}^0,x_{q^-}]$. Moreover, since $q = p^{+_B}$ for some block $B\in\P$, we have that $p,q^{--}\in C$ for some block $C\in\P^c$, see Section \ref{section non-crossing partitions}. From Definition \ref{definition B} and Definition \ref{definition almost co-aisle of co-t-structure} this implies that $a = |z_p^0,x_{q^-}|\in\ind\S^{-1}\V\cap\B$.

Now we show that any arc $a\in\ind(\U\cap\A)\cap\S^{-1}(\V\cap\B)$, provided that it exists, is of the form $a = |z_p^0,x_{q^-}|$ with $p,q\in[m']$ such that $q = p^{+_B}$ for some block $B\in\P$, and $x_{q^-}\in\Z^{(q^-)}$. We divide the argument into steps.

\emph{Step 1.} Let $z$ be an endpoint of $a$. We show that $z = z_0^p$ for some $p\in[m']$, or $z = x_{p^-}$ for some $p\in[m']$ such that $x_{p^-}\in\Z^{(p^-)}$. 

Since $a\in\ind\U\cap\A$, then $z\in[x_{p^-},z_p^0]$ for some $p\in[m']$, and since $a\in\ind\V\cap\B$, then $z\in[z_q^0,x_{q^+}]$ for some $q\in[m']$. If $z\in(p,z_p^0]$ then $q = p$ and $z = z_p^0$. If $z\in[x_{p^-},p)$ then $q = p^{--}$ and $z = x_{p^-}$. Therefore we have the claim.

\emph{Step 2.} Let $p,q\in[m']$ be such that one endpoint of $a$ is of the form $x_{p^-}$ or $z_p^0$, and the other endpoint is of the form $x_{q^-}$ or $z_q^0$.  We show that $a\not\cong|z_p^0,z_q^0|$ and $a\not\cong|x_{p^-},x_{q^-}|$.
 
If $a\cong |z_p^0,z_q^0|$ then $a$ and $\S a$ are crossing and, since $\S a = |z_p^0-1,z_q^0-1|\in\ind\U\cap\A$, this gives a contradiction. Similarly, if $a\cong|x_{p^-},x_{q^-}|$ then $a$ is crossed by $\S^{-1}a = |x_{p^-}+1,x_{q^-}+1|\in\ind\U\cap\A$ and we obtain again a contradiction.

\emph{Step 3.} We know that $a \cong |z_p^0,x_{q^-}|$ for some $p,q\in[m']$ such that $p,q\in B$ for some block $B\in\P$ and $x_{q^-}\in\Z^{(q^-)}$. We show that $q = p^{+_B}$. 

Assume that $q\neq p^{+_B}$. Then $B\neq\{p\}$, otherwise $p = q$ and $q = p^{+_B}$. If $p = q$ consider $r\in B\setminus\{p\}$. Then there exists an arc in $\ind\U\cap\A$ with an endpoint in $(p,z_p^0]$ and the other endpoint in $(r,z_0^r]$ which crosses $a$, and this gives a contradiction. Now assume that $p\neq q$, then there exists $r\in B\setminus\{p,q\}$ such that $p,r,q$ are in cyclic order. The arc $|z_0^r,z_0^q|\in\U\cap\A$ crosses $a$, and this gives a contradiction. This concludes the argument.
\end{proof}

The following corollary can be proved with the same argument of Corollary \ref{corollary heart t-structure}.

\begin{corollary}
Let $(\X,\Y)$ be a co-t-structure in $\ovl\C_m$. Consider $(\P,X)$ its associated alternating non-crossing partition of $[m']\cup[m]$ with $X = (x_p)_{p\in[m]}$. Then the co-heart $\mathcal{S} = \X\cap\S^{-1}\Y$ is given by
\[\mathcal{S} = \add\left\{|p,x_{q^-}|\middle|\text{$p,q\in[m']$, $q = p^{+_B}$ for some $B\in\P$, and $x_{q^-}\in\Z^{(q^-)}$}\right\}.\]
\end{corollary}

\subsection{Boundedness}\label{section bounded co-t-structures}

We study the bounded co-t-structures in $\ovl\C_m$. We find that for $m\geq 2$ there are no bounded co-t-structures.

\begin{proposition}\label{proposition left bounded co-t-structures}
Let $(\X,\Y)$ be a co-t-structure in $\ovl\C_m$, $(\P,X)$ be its associated alternating non-crossing partition with $X = (x_p)_{p\in[m]}$, and $\U = \pi^{-1}\X$. The following statements are equivalent.
\begin{enumerate}
\item The co-t-structure $(\X,\Y)$ is left bounded in $\ovl\C_m$.
\item We have that $\bigcup_{n\in\Z}\S^n\U = \C_{2m}$.
\item The non-crossing partition $\P$ has as unique block $\{1',\dots,m'\}$ and $x_p\neq p^+$ for each $p\in[m]$.
\end{enumerate}
\end{proposition}
\begin{proof}
The equivalence between the statements (1) and (2) is straightforward, we show the equivalence between (2) and (3). Assume that $\P = \{\{1',\dots,m'\}\}$ and $x_p\neq p^+$ for each $p\in[m]$. Let $a = (a_1,a_2)\in\ind\C_{2m}$, we check that $a\in\S^n\U$ for some $n\in\Z$. There exists $n\geq 0$ such that $a_1+n\in[x_{p^-},p^+)$ and $a_2+n\in[x_{q^-},q^+)$ for some $p,q\in[m']$. Since $p$ and $q$ belong to the same block of $\P$, we have that $\S^{-n}a = (a_1+n,a_2+n)\in\U$, and then $a\in\S^{n}\U$. 

Now assume that $\bigcup_{n\in\Z}\S^n\U = \C_{2m}$, we check that (3) holds. Let $p,q\in[m']$, and consider $a\in\ind\C_{2m}$ with an endpoint in $\Z^{(p)}$ and the other in $\Z^{(q)}$. By assumption there exists $n\in\Z$ such that $a\in\S^n\U$, and then $\S^{-n}a\in\U$. Since the endpoints of $\S^{-n}a$ still belong to $\Z^{(p)}$ and $\Z^{(q)}$, we have that $p,q\in B$ for some block $B\in\P$. This means that any two elements of $[m']$ belong to the same block of $\P$, i.e. $\P = \{\{1',\dots,m'\}\}$. Now, assume that $x_p = p^+$ for some $p\in[m]$. Consider an arc $a\in\Z^{(p,p)}$, we observe that $\S^n a\notin\U$ for each $n\in\Z$, and this gives a contradiction. This concludes the argument.
\end{proof}

Dually, we have the following proposition.

\begin{proposition}\label{proposition right bounded co-t-structures}
Let $(\X,\Y)$ be a co-t-structure in $\ovl\C_m$, let $(\P,X)$ be its associated alternating non-crossing partition with $X = (x_p)_{p\in[m]}$, and $\V = \pi^{-1}\X$. The following statements are equivalent.
\begin{enumerate}
\item The co-t-structure $(\X,\Y)$ is right bounded in $\ovl\C_m$.
\item We have that $\bigcup_{n\in\Z}\S^n\V = \C_{2m}$.
\item The non-crossing partition $\P$ has as blocks $\{1'\}$, \dots, $\{m'\}$, and $x_p\neq p$ for each $p\in[m]$.
\end{enumerate}
\end{proposition}

\begin{corollary}
For each $m\geq 2$ there are no bounded co-t-structures in $\ovl\C_m$.
\end{corollary}
\begin{proof}
Assume that $m\geq 2$. If there exists a bounded co-t-structure in $\ovl\C_m$, then, by Proposition \ref{proposition left bounded co-t-structures} and Proposition \ref{proposition right bounded co-t-structures}, its associated alternating non-crossing partition $(\P,X)$ of $[m']\cup[m]$ is such that $\P = \{1',\dots,m'\} = \{\{1'\},\dots,\{m'\}\}$, giving a contradiction. Therefore, there are no bounded co-t-structures in $\ovl\C_m$ if $m\geq 2$.
\end{proof}

\subsection{Non-degeneracy}\label{section non-degenerate co-t-structures}

We classify the non-degenerate co-t-structures in $\ovl\C_m$. We find that for $m\geq 2$ there are no non-degenerate co-t-structures. In general it is straightforward to check that left or right bounded co-t-structures are also right or left non-degenerate respectively. We will see that also the converse holds in $\ovl\C_m$.

\begin{proposition}\label{proposition left non-degenerate co-t-structures}
Let $(\X,\Y)$ be a co-t-structure in $\ovl\C_m$, $(\P,X)$ be its associated alternating non-crossing partition with $X = (x_p)_{p\in[m]}$, and $\U = \pi^{-1}\X$. The following statements are equivalent.
\begin{enumerate}
\item The co-t-structure $(\X,\Y)$ is left non-degenerate.
\item We have that $\bigcap_{n\in\Z}\S^n\U = \D$.
\item The non-crossing partition $\P$ has blocks $\{1'\}$, \dots, $\{m'\}$, and $x_p\neq p$ for each $p\in[m]$.
\end{enumerate}
\end{proposition}
\begin{proof}
For the equivalence between the statements (1) and (2) we can use the same argument of Proposition \ref{proposition left non-degenerate t-structures}. We prove the equivalence between (2) and (3). Assume that $\bigcap_{n\in\Z}\S^n\U = \D$ and that there exist $p,q\in[m']$ such that $p,q\in B$ for some $B\in\P$. Then $\U$ contains any arc having one endpoint in $\Z^{(p)}$ and the other endpoint in $\Z^{(q)}$. Consider such arc $u$, then $\S^n u\in \U$ for each $n\in\Z$, i.e. $u\in\bigcap_{n\in\Z}\S^n\U$. Then $u\in\D$ and $p = q$. Now assume that there exists $p\in[m]$ such that $x_p = p$, then $\U$ contains any arc $u\in\Z^{(p,p)}$. Thus, $u\in\bigcap_{n\in\Z}\S^n\U = \D$, and then $u\in\Z^{(q,q)}$ for some $q\in[m']$ and this contradicts the fact that $p\in[m]$. This proves that (3) holds.
	
Now assume that statement (3) holds, we check that $\bigcap_{n\in\Z}\S^n\U\subseteq\D$, the other inclusion is straightforward. Let $u\in\ind\bigcap_{n\in\Z}\S^n\U$, then $u\in\U$ and there exist $p,q\in[m']$ such that $u$ has one endpoint in $[x_{p^-},p^+)$ and the other endpoint in $[x_{q^-},q^+)$. Then $p,q\in B$ for some block $B\in\P$, and as a consequence $p = q$ and $u$ has both endpoints in $[x_{p^-},p^+)$. Assume that $u$ has an endpoint in $[x_{p^-},p)$, then, since $x_p\neq p$, there exists $n\in\Z$ such that $\S^n u\notin\U$, i.e. $u\notin\bigcap_{n\in\Z}\S^n\U$. Then $u\in\Z^{(p,p)}$, and as a consequence $u\in\D$. This concludes the argument.
\end{proof}

Dually, we have the following proposition.

\begin{proposition}\label{proposition right non-degenerate co-t-structures}
Let $(\X,\Y)$ be a co-t-structure of $\ovl\C_m$, $(\P,X)$ be its associated alternating non-crossing partition with $X = (x_p)_{p\in[m]}$, and $\V = \pi^{-1}\X$. The following statements are equivalent.
\begin{enumerate}
\item The co-t-structure $(\X,\Y)$ is right non-degenerate.
\item We have that $\bigcap_{n\in\Z}\S^n\V = \D$.
\item The non-crossing partition $\P$ has as unique block $\{1',\dots,m'\}$ and $x_p\neq p^+$ for each $p\in[m]$.
\end{enumerate}
\end{proposition}

\begin{corollary}
For each $m\geq 2$ there are no non-degenerate co-t-structures in $\ovl\C_m$.
\end{corollary}

We also have the following corollary, which combines these results with those in Section \ref{section bounded co-t-structures}.

\begin{corollary}
Let $(\X,\Y)$ be a co-t-structure in $\ovl\C_m$. Then $(\X,\Y)$ is left bounded if and only if it is right non-degenerate, and $(\X,\Y)$ is right-bounded if and only if it is left non-degenerate.
\end{corollary}

\subsection{Adjacent triples}\label{section adjacent triples and functorially finite co-t-structures}

We classify the co-t-structures in $\ovl\C_m$ having a left adjacent or right adjacent t-structure.

\begin{theorem}\label{theorem functorially finite co-t-structures}
Let $(\X,\Y)$ be a co-t-structure in $\ovl\C_m$ and $(\P,X)$ be its associated alternating non-crossing partition with $X = (x_p)_{p\in[m]}$. The following statements hold.
\begin{enumerate}
\item The co-t-structure $(\X,\Y)$ has a right adjacent t-structure if and only if for each $p\in[m]$ if $x_p = p^+$ then $\{p^+\}\in\P$.
\item The co-t-structure $(\X,\Y)$ has a left adjacent t-structure if and only if for each $p\in[m]$ if $x_p = p$ then $p^-,p^+\in B$ for some block $B\in\P$.
\end{enumerate}
\end{theorem}
\begin{proof}
We prove statement (1), statement (2) is dual. Let $\V = \pi^{-1}\Y$ and $(\Q,Y) = (\P,X)^c$ with $Y = (y_p)_{p\in[m]}$. If $(\X,\Y)$ has a right adjacent t-structure, then $\Y$ is precovering and $\V\cap\A$ satisfies the $\PC$ conditions, see Theorem \ref{theorem classification precovering subcategories}. Let $p\in[m]$ be such that $x_p = p^+$, then $y_p = p^+$. We show that $p^-,p^+\in C$ for some block $C\in\Q$. Since $y_p = p^+$, $\V\cap\A$ contains all the arcs having one endpoints in $(p^-,z_{p^-}^0]$ and the other in $\Z^{(p)}$, see Definition \ref{definition A} and Definition \ref{definition almost co-aisle of co-t-structure}. By $(\PC 3)$ or $(\PC 3')$ there exists an arc of $\V$ with an endpoint in $\Z^{(p^-)}$ and the other in $\Z^{(p^+)}$. Thus, $p^-,p^+\in C$ for some block $C\in\Q$. Since $\Q = \P^c$, this is equivalent to $\{p^+\}\in\P$.

Now assume that (2) holds, i.e. if $y_p = p^+$ then $p^-,p^+\in C$ for some block $C\in\Q$. We show that $\Y$ is precovering, i.e. that $\V\cap\A$ is precovering. We check that $\V\cap\A$ satisfies $(\PC 1)$ the other conditions are analogous. Assume that there exists a sequence $\{(v_1^n,v_2^n)\}_n\subseteq\V\cap\A\cap\Z^{(p,q)}$ for some $p,q\in[m']\cup[m]$ such that $p\neq q$ with $\{v_1^n\}_n$ and $\{v_2^n\}_n$ strictly increasing. Then $p,q\in [m]$ and, since there exist arcs of $\V$ in $\Z^{(p,q)}$, $p^-,q^-\in C$ for some $C\in\Q$. Moreover, $y_p = p^+$ and $y_q = q^+$, and then by assumption $p^-,p^+,q^-,q^+\in C$. Then, $\V\cap\A$ contains any arc having one endpoint in $(p^+,z_{p^+}^0]$ and the other endpoint in $(q^+,z_{q^+}^0]$. In particular, there exist strictly decreasing sequences $\{w_1^n\}_n\subseteq\Z^{(p^+)}$ and $\{w_2^n\}_n\subseteq\Z^{(q^+)}$ such that $\{|w_1^n,w_2^n|\}_n\subseteq\V\cap\A$. This concludes the argument.
\end{proof}

\subsection{Recollements}\label{section recollements}

We recall that in a triangulated category recollements are in bijection with TTF triples, which are triples $(\X,\Y,\Zc)$ such that $(\X,\Y)$ and $(\Y,\Zc)$ are t-structures, we refer to \cite[Section 2.2]{NS} for more details. Since $\ovl\C_m$ is $\Hom$-finite and Krull--Schmidt, by Proposition \ref{proposition iyama-yoshino torsion pairs}, TTF triples are in bijection with functorially finite thick subcategories, which we classify here. Thick subcategories of $\C_m$ and $\ovl\C_m$ were classified in \cite{GZ} and \cite{M} respectively. By \cite[Proposition 4.6]{ZZ}, $\C_m$ and $0$ are the only precovering or preenveloping thick subcategories of $\C_m$, but this is no longer the case in $\ovl\C_m$, see Figure \ref{figure example co-t-structure in completion} for an example.

The following theorem follows directly from Theorem \ref{theorem classification aisles co-t-structures}.

\begin{theorem}\label{theorem precovering and preenveloping thick subcategories}
Let $(\X,\Y)$ be a co-t-structure of $\ovl\C_m$ and $(\P,X)$ be its associated alternating non-crossing partition with $X = (x_p)_{p\in[m]}$. The following statements are equivalent.
\begin{enumerate}
\item $\X$ is a precovering thick subcategory.
\item $\Y$ is a preenveloping thick subcategory.
\item For each $p\in[m]$ either $x_p = p$ or $x_p = p^+$.
\end{enumerate}
\end{theorem}

The following corollary combines Theorem \ref{theorem functorially finite co-t-structures} and Theorem \ref{theorem precovering and preenveloping thick subcategories}. 

\begin{corollary}\label{corollary functorially finite thick subcategories}
Let $\X$ be a subcategory of $\ovl\C_m$. The following statements are equivalent.
\begin{enumerate}
\item $\X$ is a functorially finite thick subcategory.
\item The alternating non-crossing partition of $[m']\cup[m]$ associated to the co-t-structure $(\X,\X^\perp)$, which we denote by $(\P,X)$ with $X = (x_p)_{p\in[m]}$, satisfies the following condition: for each $p\in[m]$ either $x_p = p$ or $x_p = p^+$, and if $x_p = p$ then $p^-,p^+\in B$ for some block $B\in\P$.
\end{enumerate}
\end{corollary}

\section{Lattice structures}

The t-structures and co-t-structures in a triangulated category form a partially ordered sets under inclusion of aisles. In this section we prove that both t-structures and co-t-structures in $\ovl\C_m$ form a lattice similarly as the t-structures in $\C_m$, see \cite[Section 5]{GZ}. In Section \ref{section t-structures} and Section \ref{section co-t-structures} we proved that the t-structures and the co-t-structures in $\ovl\C_m$ are in bijection with, respectively, the half-decorated non-crossing partitions and the alternating non-crossing partitions of $[m']\cup[m]$. We prove that these sets have lattice structures. We recall that the non-crossing partitions forms a lattice under refinement, see Section \ref{section non-crossing partitions}.

We introduce some notation. Given $m$-tuples $X = (x_p)_{p\in[m]}$ and $X' = (x'_p)_{p\in[m]}$, we write $X\leq X'$ if $x_p\leq x'_p$ for each $p\in[m]$. We write 
\[
\min\{X,X'\} = (\min\{x_p,x'_p\})_{p\in[m]}\text{ and }\max\{X,X'\} = (\max\{x_p,x'_p\})_{p\in[m]}.
\] 
For $(\P,X)$ and $(\P',X')$ half-decorated non-crossing partitions of $[m']\cup[m]$, we define $(\P,X)\leq (\P',X')$ if $\P\leq\P'$ and $X\leq X'$, as in \cite[Section 5]{GZ}. If $(\P,X)$ and $(\P',X')$ are alternating non-crossing partitions of $[m']\cup[m]$, then we define $(\P,X)\leq (\P',X')$ if $\P\leq\P'$ and $X'\leq X$.

\begin{lemma}\label{lemma lattice structures half decorated non-crossing partitions and alternating non-crossing partitions}
The following statements hold.
\begin{enumerate}
\item The half-decorated non-crossing partitions of $[m']\cup[m]$ form a lattice where, for each $(\P,X)$ and $(\P',X')$, we have that
\begin{align*}
(\P,X)\wedge(\P',X') & = (\P\wedge\P',\min\{X,X'\})\text{ and}\\
(\P,X)\vee(\P',X') & = (\P\vee\P',\max\{X,X'\}).
\end{align*}
\item The alternating non-crossing partitions of $[m']\cup[m]$ form a lattice where, for each $(\P,X)$ and $(\P',X')$, we have that
\begin{align*}
(\P,X)\wedge(\P',X') & = (\P\wedge\P',\max\{X,X'\})\text{ and}\\
(\P,X)\vee(\P',X') & = (\P\vee\P',\min\{X,X'\}).
\end{align*}
\end{enumerate}
\end{lemma}
\begin{proof}
We prove statement (1). Let $(\P,X)$ and $(\P',X')$ be half-decorated non-crossing partitions of $[m']\cup[m]$. It is straightforward to check that $(\P\wedge\P',\min\{X,X'\})$ and $(\P\vee\P',\max\{X,X'\})$ are respectively the greater lower bound and the least upper bound of $(\P,X)$ and $(\P',X')$, provided that they are well defined. We check that they both are half-decorated non-crossing partitions of $[m']\cup[m]$.
	
Let $(\P,\wt X)$ and $(\P,\wt X')$ be the decorated non-crossing partitions of $[m']\cup[m]$ obtained respectively from $(\P,X)$ and $(\P',X')$ by adding the decoration $z_p^0$ for each $p\in[m']$, see Remark \ref{remark half-decorated non-crossing partitions and decorated non-crossing partitions}. The set of decorated non-crossing partitions has a partial order defined in \cite[Section 5]{GZ}. By \cite[Theorem 5.2]{GZ} the meet and join of $(\P,\wt X)$ and $(\P',\wt X')$ are respectively $(\P\wedge\P',\min\{\wt X,\wt X'\})$ and $(\P\vee\P',\max\{\wt X,\wt X'\})$. Note that $(\P\wedge\P',\min\{\wt X,\wt X'\})$ and $(\P\vee\P',\max\{\wt X,\wt X'\})$ have $z_p^0$ as decoration for each $p\in[m']$. Thus,  $(\P\wedge\P',\min\{X,X'\})$ and $(\P\vee\P',\max\{X,X'\})$ are their correspondent half-decorated non-crossing partitions under the bijection of Remark \ref{remark half-decorated non-crossing partitions and decorated non-crossing partitions}. Therefore, $(\P\wedge\P',\min\{X,X'\})$ and $(\P\vee\P',\max\{X,X'\})$ are well defined.

Now we prove (2). Let $(\P,X)$ and $(\P',X')$ be alternating non-crossing partitions of $[m']\cup[m]$. Since $\P\wedge\P'$ and $\P\vee\P'$ are non-crossing partitions of $[m']$, then, by Definition \ref{definition alternating non-crossing partition}, $(\P\wedge\P',\max\{X,X'\})$ and  $(\P\vee\P',\min\{X,X'\})$ are alternating non-crossing partitions of $[m']\cup[m]$. Moreover, it is straightforward to check that these are respectively the greater lower bound and the least upper bound of $(\P,X)$ and $(\P',X')$.
\end{proof}

The following lemma is the analogue of \cite[Proposition 5.3]{GZ}. We recall that a co-aisle of a t-structure or of a co-t-structure is related to the Kreweras complement of the corresponding non-crossing partition, see Section \ref{section co-aisles t-structures} and Section \ref{section co-aisles of co-t-structures}. The complements of half-decorated non-crossing partitions and of alternating non-crossing partitions are defined in Definition \ref{definition complement half decorated non-crossing partition} and Definition \ref{definition complement co-aisle co-t-structure}.

\begin{lemma}\label{lemma computation intersection aisles}
The following statements hold.
\begin{enumerate}
\item Let $(\X,\Y)$ and $(\X',\Y')$ be t-structures in $\ovl\C_m$, $(\P,X)$ and $(\P',X')$ be their correspondent half-decorated non-crossing partitions of $[m']\cup[m]$, and $(\Q,Y) = (\P,X)^c$ and $(\Q',Y') = (\P',X')^c$ be their complements. Then
\[
\X\cap\X' = \pi\U_{(\P\wedge\P',\min\{X,X'\})} \text{ and } \Y\cap\Y' = \pi\V_{(\Q\wedge\Q',\max\{Y,Y'\})}.
\]
\item  Let $(\X,\Y)$ and $(\X',\Y')$ be co-t-structures in $\ovl\C_m$, $(\P,X)$ and $(\P',X')$ be their correspondent alternating non-crossing partitions of $[m']\cup[m]$, and $(\Q,Y) = (\P,X)^c$ and $(\Q',Y') = (\P',X')^c$ be their complements. Then
\[
\X\cap\X' = \pi\U_{(\P\wedge\P',\max\{X,X'\})} \text{ and } \Y\cap\Y' = \pi\V_{(\Q\wedge\Q',\min\{Y,Y'\})}.
\]
\end{enumerate}
\end{lemma}
\begin{proof}
We prove (1), statement (2) is analogous. We check that $\X\cap\X' = \pi\U_{(\P\wedge\P',\min\{X,X'\})}$, for the equality $\Y\cap\Y' = \pi\V_{(\Q\wedge\Q',\max\{Y,Y'\})}$ we can proceed analogously. By Theorem \ref{theorem t-structures}, $\X = \pi\U_{(\P,X)}$ and $\X' = \pi\U_{(\P',X')}$. We denote $\U = \U_{(\P,X)}$ and $\U' = \U_{(\P',X')}$. By applying a similar argument to the one of \cite[Proposition 5.3]{GZ}, we obtain that $\U\cap\U' = \U_{(\P\wedge\P',\min\{X,X'\})}$. We prove that $\pi\U\cap\pi\U' = \pi(\U\cap\U')$.
	
Let $x\in\pi\U\cap\pi\U'$, then $x\cong\pi(u)$ and $x\cong\pi(u')$ for some $u\in\U$ and $u'\in\U'$. Thus, $\pi(u)\cong\pi(u')$ and, by Lemma \ref{lemma Q(X) closed under isomorphisms}, $u,u'\in\U\cap\U'$. Then, $x\in\pi(\U\cap\U')$. This implies that $\pi\U\cap\pi\U'\subseteq\pi(\U\cap\U')$, the other inclusion in straightforward. It follows that $\X\cap\X' = \pi(\U\cap\U') = \pi\U_{(\P\wedge\P',\min\{X,X'\})}$.
\end{proof}

We can now describe the lattice structures of the t-structures and co-t-structures in $\ovl\C_m$. We refer to \cite[Theorem 5.2, Proposition 5.3]{GZ} for the non-completed case.

\begin{theorem}\label{theorem lattice structures t-structures and co-t-structures}
The t-structures and the co-t-structures in $\ovl\C_m$ have lattice structures under inclusion of aisles. For each pair of t-structures, or of co-t-structures, $(\X,\Y)$ and $(\X',\Y')$, we have that
\begin{align*}
(\X,\Y)\wedge(\X',\Y') &  = (\X\cap\X',(\X\cap\X')^\perp)\text{ and}\\
(\X,\Y)\vee(\X',\Y') & = (^\perp(\Y\cap\Y'),\Y\cap\Y').
\end{align*}
\end{theorem}
\begin{proof}
We prove the statement for the t-structures in $\ovl\C_m$. For the co-t-structures we can proceed similarly. We divide the proof into steps.
	
\emph{Step 1.} We prove that the t-structures form a lattice under inclusion of aisles.

By Corollary \ref{corollary t-structures} there is an inclusion preserving bijection between the aisles of t-structures in $\ovl\C_m$ and the suspended subcategories $\U$ of $\C_{2m}$ such that $\D\subseteq\U$ and $\U\cap\A$ is precovering. Moreover, by Proposition \ref{proposition classification t-structures} this set of subcategories of $\C_{2m}$ is in bijection with the set of half-decorated non-crossing partitions of $[m']\cup[m]$, which is a lattice by Lemma \ref{lemma lattice structures half decorated non-crossing partitions and alternating non-crossing partitions}. We prove that the latter bijection is order preserving. Let $\U$ and $\U'$ be suspended subcategories of $\C_{2m}$ as above,  and let $(\P,X)$ and $(\P',X')$ be their associated half-decorated non-crossing partitions. We prove that $\U\subseteq\U'$ if and only if $(\P,X)\leq (\P',X')$.

Assume that $\U\subseteq\U'$ and consider $p,q\in[m']\cup[m]$ such that $p,q\in B$ for some block $B$ of $\P$. We prove that $p$ and $q$ belong to the same block of $\P'$. There is an arc $u\in\ind\U$ with an endpoint in $\Z^{(p)}$ and the other in $\Z^{(q)}$. Since $u\in\U'$, we obtain that $p,q\in B'$ for some block $B'$ of $\P'$. Thus, $\P\leq\P'$. Moreover, since $x_p'\geq z'$ for each $z'\in\Z^{(p)}$ which is an endpoint of an arc of $\U'$, then  $x_p'\geq z$ for each $z\in\Z^{(p)}$ which is an endpoint of an arc of $\U$. Therefore, $x_p\leq x_p'$ because $x_p$ is the least upper bound of the set of such elements $z\in\Z^{(p)}$. We obtain that $(\P,X)\leq (\P',X')$.

Now assume that $(\P,X)\leq(\P',X')$ and consider $u\in\ind\U$, we show that $u\in\ind\U'$. There exists a block $B$ of $\P$ such that each endpoint of $u$ belongs to $(p,x_p]$ for some $p\in B\cap[m]$, or to $\Z^{(p)}$ for some $p\in B\cap[m']$. Since $B\subseteq B'$ and $(p,x_p]\subseteq(p,x_p']$ for $p\in[m]$, each endpoint of $u$ belongs to $(p,x_p']$ for some $p\in B'\cap[m]$, or to $\Z^{(p)}$ for some $p\in B'\cap[m']$. Therefore, $u'\in\ind\U'$. We conclude that $\U\subseteq\U'$ if and only if $(\P,X)\leq(\P',X')$.

\emph{Step 2.} We compute the meet of two t-structures.

Consider $(\X,\Y)$ and $(\X',\Y')$ t-structures in $\ovl\C_m$ and their corresponding half decorated non-crossing partitions $(\P,X)$ and $(\P',X')$. The meet of $(\X,\Y)$ and $(\X',\Y')$ corresponds to the meet of $(\P,X)$ and $(\P',X')$, which is equal to $(\P\wedge\P',\min\{X,X'\})$ by Lemma \ref{lemma lattice structures half decorated non-crossing partitions and alternating non-crossing partitions}. Thus, the aisle of $(\X,\Y)\wedge(\X',\Y')$ is equal to $\pi\U_{(\P\wedge\P',\min\{X,X'\})}$. By Lemma \ref{lemma computation intersection aisles} we obtain that $(\X,\Y)\wedge(\X',\Y')  = (\X\cap\X',(\X\cap\X')^\perp)$.

\emph{Step 3.} We compute the join of two t-structures.
	
The join of $(\X,\Y)$ and $(\X',\Y')$ corresponds to the join of $(\P,X)$ and $(\P',X')$, which is equal to $(\P\vee\P',\max\{X,X'\})$ by Lemma \ref{lemma lattice structures half decorated non-crossing partitions and alternating non-crossing partitions}. We denote $(\Q,Y) = (\P,X)^c$, $(\Q',Y') = (\P',X')^c$, and $(\R,Z) = (\P\vee\P',\max\{X,X'\})^c$. By Remark \ref{remark join and meet kreweras complement}, $\R = (\P\vee\P')^c = \Q\wedge\Q'$, and then $(\R,Z) = (\Q\wedge\Q',\max\{Y,Y'\})$. Since the co-aisle of $(\X,\Y)\vee(\X,\Y)$ is equal to $\pi\V_{(\R,Z)}$, by Lemma \ref{lemma computation intersection aisles} we have that $\pi\V_{(\R,Z)} = \Y\cap\Y'$.
\end{proof}

The lattice structures described above restrict to certain classes of t-structures and of co-t-structures in $\ovl\C_m$. We recall that a \emph{sublattice} of a lattice is a subposet which contains the join and meet of any pair of elements.

\begin{corollary}\label{corollary sublattices t-structures and co-t-structures}
The following statements hold.
\begin{enumerate}
\item The left bounded and the right bounded t-structures form  sublattices of the lattice of t-structures in $\ovl\C_m$.
\item The left bounded, the right bounded, the left non-degenerate, and the right non-degenerate co-t-structures form sublattices of the lattice of co-t-structures in $\ovl\C_m$.
\item The co-t-structures having a left adjacent t-structure, and  the co-t-structures having a right adjacent t-structure form sublattices of the lattice of co-t-structures in $\ovl\C_m$.
\end{enumerate}
\end{corollary}
\begin{proof}
The statements (1) and (2) are straightforward and follow directly from the combinatorial descriptions of Proposition \ref{proposition bounded above t-structures}, Proposition \ref{proposition bounded below t-structures}, Proposition \ref{proposition left bounded co-t-structures}, Proposition \ref{proposition right bounded co-t-structures}, Proposition \ref{proposition left non-degenerate co-t-structures}, and Proposition \ref{proposition right non-degenerate co-t-structures}. We prove (3).
	
Let $(\X,\Y)$ and $(\X',\Y')$ be co-t-structures in $\ovl\C_m$ admitting left adjacent t-structures, we check that their join and meet admit left adjacent t-structures. If $(\X,\Y)$ and $(\X',\Y')$ admit right adjacent t-structures the proof is similar. Let $(\P,X)$ and $(\P',X')$ be the corresponding alternating non-crossing partitions of $[m']\cup[m]$ with $X = (x_p)_{p\in[m]}$ and $X' = (x_p')_{p\in[m]}$. By Lemma \ref{lemma lattice structures half decorated non-crossing partitions and alternating non-crossing partitions}, $(\P,X)\wedge (\P',X') = (\P\wedge\P',\max\{X,X'\})$ and $(\P,X)\vee (\P',X')  = (\P\vee\P',\min\{X,X'\})$. Therefore, by Theorem \ref{theorem functorially finite co-t-structures} it is enough to prove that: for each $p\in[m]$, if $\max\{x_p,x_p'\} = p$ then $p^-,p^+\in C$ for some block $C$ of $\P\wedge\P'$, and if $\min\{x_p,x_p'\} = p$ then $p^-,p^+\in\C$ for some block $C$ of $\P\vee\P'$.

Let $p\in[m]$. If $\max\{x_p,x_p'\} = p$ then $x_p = x_p' = p$. Thus, $p^-,p^+\in B$ for some block $B$ of $\P$, and $p^-,p^+\in B'$ for some block $B'$ of $\P'$. Therefore, $p^-,p^+\in B\cap B'$ which is a block of $\P\wedge\P'$. Now assume that $\min\{x_p,x_p'\} = p$, and assume that $x_p = p$, for the case $x_p' = p$ we can proceed in the same way. We have that $p^-,p^+\in B$ for some block $B$ of $\P$. Since $\P\leq\P\vee\P'$, there exists a block $C$ of $\P\vee\P'$ such that $B\subseteq C$ and then $p^-,p^+\in C$. We conclude that the join and meet of $(\X,\Y)$ and $(\X',\Y')$ admit left adjacent t-structures.
\end{proof}

We observe that the non-degererate t-structures in $\ovl\C_m$ do not form a sublattice, below we have a counterexample. The same holds for the left non-degenerate and the right non-degenerate t-structures.

\begin{example}
Let $\P = \{\{1',1\},\{2',2\}\}$ and $\P' = \{\{1'\},\{2'\},\{1,2\}\}$ be non-crossing partitions of $[2']\cup[2]$. Let $X = (x_p)_{p\in[2]}$ and $X' = (x_p')_{p\in[2]}$ be such that $x_p,x_p'\in\Z^{(p)}$ for each $p\in[2]$. Then $(\P,X)$ and $(\P,X')$ are half-decorated non-crossing  partitions of $[2']\cup[2]$ such that their corresponding t-structures are non-degenerate, see Corollary \ref{corollary non-degenerate t-structures}. Note that $\P\vee\P'$ has as unique block $\{1',1,2',2\}$ and then the t-structures associated to $(\P,X)\vee(\P',X')$ is not left non-degenerate. Moreover, $\P\wedge\P'$ has as blocks $\{1'\},\{1\},\{2'\},\{2\}$ and then the t-structure associated to $(\P,X)\wedge(\P',X')$ is not right non-degenerate.
\end{example}

We know from \cite[Theorem 3.7]{GZ} and \cite[Theorem 4.9]{M} that the thick subcategories of $\C_m$ and of $\ovl\C_m$ form lattice structures. We observe that the functorially finite thick subcategories of $\ovl\C_m$ also form a lattice.

\begin{corollary}\label{corollary functorially finite thick subcategories lattice structure}
The functorially finite thick subcategories of $\ovl\C_m$ form a lattice under inclusion.
\end{corollary}
\begin{proof}
A functorially finite thick subcategory $\X$ of $\ovl\C_m$ can be regarded as the aisle of the co-t-structure $(\X,\X^\perp)$ which admits a left adjacent t-structure and is such that $\S\X\subseteq\X$. By the combinatorial description of Theorem \ref{theorem precovering and preenveloping thick subcategories} and Corollary \ref{corollary sublattices t-structures and co-t-structures}, both these classes of co-t-structures are closed under taking finite joins and meets. We conclude that the lattice structure of the co-t-structures in $\ovl\C_m$ restricts to the functorially finite thick subcategories.
\end{proof}

\end{document}